\documentclass[11pt]{article}

\usepackage{latexsym}
\usepackage{amsfonts}
\usepackage{amsmath}

\usepackage[dvips]{color}

\newtheorem{thrm}{Theorem}[section]
\newtheorem{lemma}[thrm]{Lemma}
\newtheorem{prop}[thrm]{Proposition}

\newtheorem{cor}[thrm]{Corollary}

\numberwithin{equation}{section}

\usepackage{enumerate}

\usepackage{amssymb}
\usepackage{mathtools}
\mathtoolsset{showonlyrefs}
\usepackage{mathrsfs}
\usepackage{comment}
\usepackage{hyperref}
\usepackage{xcolor}

\def\E{\mathbb{E} }
\def\P{\mathbb{P} }
\def\Q{\mathbb{Q} }
\def\R{\mathbb{R} }

\makeatletter
\begin{document}
\allowdisplaybreaks

\title{\Large \bf{ $1$-stable fluctuation of the derivative martingale of branching random walk
}\footnote{The research of this project is supported
		by the National Key R\&D Program of China (No. 2020YFA0712900).}}
\author{ \bf  Haojie Hou \hspace{1mm}\hspace{1mm}
	Yan-Xia Ren\footnote{Corresponding author.
The research of this author is supported by NSFC
	(Grant Nos.  12071011 and 12231002) and The Fundamental Research Funds for the Central Universities, Peking University LMEQF.
	} \hspace{1mm}\hspace{1mm} and \hspace{1mm}\hspace{1mm}
	Renming Song\thanks{Research supported in part by a grant from the Simons
	Foundation(\#960480, Renming Song).}
	\hspace{1mm} }
\date{}
\maketitle

\begin{abstract}
In this paper, we study the functional convergence in law of the fluctuations of the derivative martingale of branching random walk on the real line. Our main result strengthens the
results of  Buraczewski et. al. [Ann. Probab., 2021]
and is the branching random walk counterpart of the main result of  Maillard and Pain [Ann. Probab., 2019] for branching Brownian motion.
\end{abstract}

\medskip


\medskip


\section{Introduction}

Consider the following branching Brownian motion on $\R$: initially there is a particle at 0, it moves according to a standard Brownian motion with drift $1$. After an exponentially distributed time with parameter $\beta>0$, it dies and splits into a random number of offspring with law $\{p_k:k \geq 0\}$.
The offspring repeat the parent's behavior independently from where they were born.
We will use $N(t)$ to denote the set of particles alive at time $t$ and for $u\in N(t)$, we will use $X_u(t)$ to denote the position of $u$.  Without loss of generality, assume that
\[
\beta \left(\sum_{k=1}^\infty k p_k -1\right)=\frac{1}{2},
\]
which implies that for any $t\geq 0$ (for example, see \cite[(1.2) and (1.3)]{MP}),
\begin{align}
	\mathbb{E}\bigg(\sum_{u\in N(t)} e^{-X_u(t)}\bigg)	=1,\quad \mathbb{E}\bigg(
	\sum_{u\in N(t)} X_u(t)e^{-X_u(t)}	\bigg)	=0\quad\mbox{and}\quad \mathbb{E}\bigg(	\sum_{u\in N(t)} (X_u(t))^2 e^{-X_u(t)}	\bigg)	=t.
\end{align}
The derivative martingale of the branching Brownian motion is defined as
\[
Z_t: = \sum_{u\in N(t)} X_u(t)e^{-X_u(t)}.
\]
It was proved in \cite{Ky04, YR} that $Z_t$ converges almost surely to
a non-degenerate non-negative
limit $Z_\infty$ if and only if $\sum_{k=1}^\infty k (\log k)^2 p_k <\infty$.
Maillard and Pain
\cite{MP} studied the fluctuation of $Z_\infty -Z_t$. They showed that, under the assumption  $\sum_{k=1}^\infty k(\log k)^3 p_k<\infty$,
\begin{align}
	\bigg(\sqrt{t}\bigg(Z_\infty- Z_{at}+\frac{\log t}{\sqrt{2\pi at}}Z_\infty\bigg)_{a\geq 1}, \P(\cdot \big| 	\mathcal{F}_t^B) \bigg)\stackrel{f.d.d.}{\longrightarrow}  \left((S_{a^{-1/2}Z_\infty })_{a\geq 1}, \P\left(\cdot \big| Z_\infty\right) \right)
\end{align}
in probability,
where $S_t$ is a spectrally positive 1-stable L\'{e}vy process independent of $Z_\infty$ and  $\{\mathcal{F}_t^B\}_{t\ge 0}$ is the filtration of the branching Brownian motion.
More precisely, they showed that, for all
$m \geq 1, a_{1}, \ldots, a_{m} \in[1, \infty)$ and bounded continuous $f: \mathbb{R}^{m} \rightarrow \mathbb{R}$,
\begin{align}\label{FDD-BBM}
	&\mathbb{E}\bigg(f\bigg(\sqrt{t}\bigg(Z_\infty- Z_{a_k t}+\frac{\log t}{\sqrt{2\pi a_k t}}Z_\infty\bigg),
	1 \leq k \leq m\bigg)\Big| 	\mathcal{F}_t^B  \bigg)	\nonumber\\
	&\underset{t \rightarrow \infty}{\longrightarrow} \mathbb{E}\left(f\left(S_{a_k^{-1/2}Z_{\infty} }, 1 \leq k \leq m\right) \Big| Z_{\infty}\right), \quad \text { in probability. }
\end{align}

Now we turn to branching random walks.  A branching random walk on $\R$ is defined as follows. At generation 0, there is a particle at the origin. At generation $n=1$, this particle dies
and splits into a finite number of offspring. The law of the offspring (offspring number and positions relative to their parent) is given by a point process $L$. The offspring evolve independently and obey the same rule as their parent. The procedure goes on.
We will use $\mathbb{T}$ to denote the genealogical tree of the branching random walk, $\mathcal{N}(n)$ to denote the collection of particles in the
$n$th generation, $|x|$ to denote the generation of particle $x$ and
$\left\{ V(x), x\in \mathcal{N}(n) \right\}$
to denote the positions of the particles in the $n$th generation.
 We will use $\P$ to denote the law of the branching random walk above and
 use $\E$ to denote the expectation with respect to $\P$.
 If the initial particle is located at $x\in \R$ instead of the origin, we will use $\P_x$
to denote the law of the corresponding
branching random walk
and use $\E_x$ to denote the expectation with respect to $\P_x$
 For $n\ge 0$, we denote by $\mathcal{F}_n$
the $\sigma$-field generated by the branching random walk up to generation $n$.  We will always assume that
\begin{itemize}
	\item [{\bf(A1)}]
	\[
	\mathbb{E}\bigg(\sum_{x\in \mathcal{N}(1)} e^{-V(x)} \bigg)	=1,\quad \mathbb{E}\bigg(\sum_{x\in \mathcal{N}(1)} V(x) e^{-V(x)}\bigg)=0
	\]
and
	\[
	\sigma^2 := \mathbb{E}\bigg(\sum_{x\in \mathcal{N}(1)} \left(V(x)\right)^2 e^{-V(x)}\bigg)< \infty.
	\]
\end{itemize}
Under {\bf(A1)},
\[
W_n:=\sum_{x\in \mathcal{N}(n)} e^{-V(x)},\quad D_n:= \sum_{x\in \mathcal{N}(n)} V(x)e^{-V(x)},\quad n\geq 0,
\]
are martingales with respect to $\{\mathcal{F}_n: n\ge 0\}$.
They are called the  additive martingale and the  derivative martingale of the  branching random walk respectively. Suppose that
\begin{itemize}
	\item [{\bf(A2)}]
	\begin{align*}
		\mathbb{E}\left(W_1 \left(\log_+ W_1\right)^2\right) + \mathbb{E}\left(\widetilde{W}_1 \log_+ \widetilde{W}_1\right)< \infty,
	\end{align*}
\end{itemize}
where $\log_+ y :=\max\{0, \log y\}$ and
\[
\widetilde{W}_1:= \sum_{x\in \mathcal{N}(1)} (V(x))_+ e^{-V(x)}
\]
with $(V(x))_+:= \max\{V(x), 0\}$.
It was proved in Biggins and Kyprianou \cite{BK2} and Chen \cite{Chen} that, under \textbf{(A1)}, $D_n$ converges almost surely to a non-negative limit $D_\infty$ with $\P (D_\infty >0)>0$  if and only if \textbf{(A2)} holds.  A\"{i}d\'ekon and Shi \cite{AESZ} studied the relationship between $W_n$ and $D_n$, and showed that, under the assumptions {\bf(A1)} and {\bf(A2)},
\begin{align}\label{Seneta-Heyde-scaling}
	\lim_{n\to\infty} \sqrt{n}W_n = \sqrt{\frac{2}{\pi \sigma^2}}D_\infty,\quad \mbox{in probability}.
\end{align}
 Under {\bf(A1)}, {\bf(A2)} and
  the additional assumption
\begin{itemize}
	\item [{\bf(A3)}] The branching random walk is non-arithmetic, i.e., for any $\delta>0$,
	\[
	\P\left(L\left(\R \setminus \delta \mathbb{Z}\right) > 0\right)>0,
	\]
\end{itemize}
Buraczewski, Iksanov and Mallein \cite{BIM} proved that
\begin{align}\label{Constant_c}
	\lim_{y\to+\infty} \left( \mathbb{E} \left(D_\infty 1_{\{D_\infty \leq y\}}\right) - \log y\right) =c_0
\end{align}
for some real number $c_0$ if and only if
\begin{align}
	\mathbb{E} \left(W_1^+\left(\log_+ W_1^+\right)^3\right) + \mathbb{E}\left(\widetilde{W}_1 \left(\log_+ \widetilde{W}_1\right)^2\right) + \mathbb{E}\bigg(\sum_{x\in\mathcal{N}(1)}e^{-V(x)} \left(-V(x)\right)_+^3\bigg)
	<\infty \label{Assumption_5}
\end{align}
and
\begin{align}
	\mathbb{E} \left(W_{1}^{-}\left(\log_+ W_{1}^{-}\right)^{3} 1_{\left\{\widehat{W}_1>C_{0}\right\}}\right)<\infty \text { for some } C_{0}>0. \label{Assumption_6}
\end{align}
Here, $W_1^+, W_1^-$ and $\widehat{W}_1$ are defined respectively by
\begin{align*}
	W_1^+&:= \sum_{x\in\mathcal{N}(1)} e^{-V(x)}1_{\{V(x) \geq  0\}} ,\quad W_1^-:= \sum_{x\in\mathcal{N}(1)} e^{-V(x)}1_{\{V(x) < 0\}},\\
	\widehat{W}_1&:= \sum_{x\in\mathcal{N}(1)}\left(1+V(x)-\min_{y\in\mathcal{N}(1)}V(y)\right) e^{\min_{y\in\mathcal{N}(1)}V(y)-V(x)} 1_{\left\{V(x)<0\right\}}.
\end{align*}
The following  sufficient condition for \eqref{Constant_c} was given in  \cite[Remark 2.3(2)]{BIM}:
\begin{itemize}
	\item [{\bf(A4)}]
	\[
	\mathbb{E} \left(W_1\left(\log_+ W_1\right)^3\right) + \mathbb{E}\left(\widetilde{W}_1 \left(\log_+ \widetilde{W}_1\right)^2\right) <\infty.
	\]
\end{itemize}
\cite[Theorem 2.4]{BIM} says that, under conditions {\bf(A1)}--{\bf(A4)}, for any bounded continuous function $f:\R \to \R$, it holds that
\begin{align}\label{PROP-2-1}
	\mathbb{E}\left(f\left(\sqrt{n}\left(D_\infty-D_n+ \frac{\log n}{2} W_{n}\right)\right) \Big|\mathcal{F}_{n}\right) \stackrel{\mathbb{P}}{\rightarrow} \mathbb{E}\left(f(D_\infty X_1) \mid
 D_{\infty}\right),   \quad n \rightarrow \infty,
\end{align}
where $X_1$ is  a spectrally positive 1-stable random variable
independent of $D_{\infty}$ with generating triplet $\left((c_0+1-\gamma)\sqrt{2 /\left(\pi \sigma^{2}\right)},\sqrt{\pi /\left(2 \sigma^{2}\right)}, 1\right)$, $\gamma$ is the Euler-Mascheroni constant and $c_0$ is the constant in \eqref{Constant_c}. More precisely, the characteristic function of $X_1$ is given by
\begin{align}\label{Charact-Function}
	\mathbb{E}\left( e^{\mathrm{i} \lambda X_1} \right)&=\exp \left\{\mathrm{i}(c_0+1-\gamma)\sqrt{2/\left(\pi\sigma^2\right)} \lambda-\sqrt{\pi/\left(2\sigma^2\right)}|\lambda|(1+\mathrm{i} \operatorname{sgn}(\lambda)(2 / \pi) \log |\lambda|)\right\}\nonumber\\ &=
\exp\left\{-\psi_{\sqrt{\pi/(2\sigma^2)}, (c_0+1-\gamma)\sqrt{2/(\pi\sigma^2)}}\right\}, \quad \lambda \in \mathbb{R},
\end{align}
where $\psi_{\sigma,\mu}(\lambda):= \sigma |\lambda|(1+\mathrm{i} \operatorname{sgn}(\lambda)(2 / \pi) \log |\lambda|) - \mathrm{i}\mu \lambda.$
Combining \eqref{Seneta-Heyde-scaling} and \eqref{PROP-2-1}, we can easily get the following fact: $\forall \varepsilon>0$,
\begin{align}\label{step_31}
	 \lim_{n\to\infty} \mathbb{P}\left(|D_\infty - D_n| > n^{\varepsilon -1/2}\right)=0.
\end{align}
The goal of this paper is to prove the counterpart of \eqref{FDD-BBM} for branching random walks. We assume the additional  assumption
\begin{itemize}
	\item [{\bf(A5)}] There exists a constant $\alpha \in (0,1]$ such that
	\[
 	\mathbb{E}
 	 	\bigg(
 	\sum_{u\in\mathcal{N}(1)} e^{-(1+\alpha)V(u)}
  	\bigg) <\infty.
	\]
\end{itemize}
Let $\{(S_n)_{n\geq 0}, \mathbf{P}\}$
be the random walk defined in \eqref{step_1}
below and let $\mathbf{E}$ denote the corresponding expectation. Then
{\bf(A5)} says that
$\mathbf{E}(e^{-\alpha S_1})<\infty$.
This assumption  is only used in \eqref{step_18}.

\section{Main result}\label{s:mainresults}
We will always assume that {\bf(A1)}--{\bf(A5)} hold. Define $\lceil y \rceil:= \min\{k\in\mathbb{Z}: k\geq y\}$.

\begin{thrm}\label{thm1}
	Let $\left(X_{t}\right)_{t \geq 0}$ be a spectrally positive 1-stable L\'{e}vy process
	with characteristic function given in \eqref{Charact-Function}, independent of $D_{\infty}$. Then the conditional law of
$$
\bigg(\sqrt{n}\bigg( D_{\infty}-D_{\lceil an \rceil}+\frac{\log n}{2} W_{\lceil an \rceil}\bigg)\bigg)_{a \geq 1}
$$
given $\mathcal{F}_{n}$ converges weakly in probability (in the sense of finite-dimensional distributions) to the conditional law of $\left(X_{a^{-1/2}D_{\infty}}\right)_{a \geq 1} \text { given } D_{\infty}.$  In other words, for all $ m \geq 1, a_{1}, \ldots, a_{m} \in[1, \infty)$ and bounded continuous $f: \mathbb{R}^{m} \rightarrow \mathbb{R}$, we have
	\begin{align}
		&\mathbb{E}\bigg( f \bigg(	\sqrt{n}\bigg(D_{\infty}-D_{\lceil a_k n\rceil}+\frac{\log n}{2} W_{\lceil a_kn \rceil}	\bigg),	1 \leq k \leq m\bigg) \Big| \mathcal{F}_{n}\bigg)
		\nonumber\\
		&\underset{n \rightarrow \infty}{\longrightarrow} \mathbb{E}\left(f\left(X_{a_k^{-1/2}D_{\infty} }, 1 \leq k \leq m\right) \Big| D_{\infty}\right), \quad \text { in probability. }
	\end{align}
\end{thrm}

Recall that $\{(S_n)_{n\geq 0}, \mathbf{P} \}$ is the random walk defined in \eqref{step_1} below.
It follows from \cite[(2.8)]{AESZ} that there exists $\theta^*>0$ such that
\begin{align}
         \lim_{n\to\infty} \sqrt{n}\mathbf{P}\Big(
	\min_{j\leq n} S_j \geq 0 \Big)=\theta^*.
\end{align}
Set
\begin{align}\label{Delta-n}
	\delta_n:= \left(\theta^*\right)^{-1}\sqrt{n}
	\mathbf{P}\Big(\min_{j\leq n} S_j \geq 0 \Big).
\end{align}

\begin{prop}\label{prop5}
	There exists a $\delta_+> 0$ such that
	\[
	\lim_{n\to\infty}\mathbb{P}\bigg(\bigg|
	\sqrt{n}W_n -\sqrt{\frac{2}{\pi \sigma^2}}\delta_n D_\infty \bigg| \geq n^{-\delta_+} \bigg)=0.
	\]
	Consequently, for all $ m \geq 1, a_{1}, \ldots, a_{m} \in[1, \infty)$ and bounded continuous $f: \mathbb{R}^{m} \rightarrow \mathbb{R}$,
	\begin{align}
		&\mathbb{E}\bigg(f\bigg( \sqrt{n}\bigg( D_{\infty}-D_{\lceil a_k n\rceil}+\frac{\log n}{\sqrt{2\pi \sigma^2 \lceil a_k n\rceil}} \delta_{\lceil a_k n\rceil} D_\infty \bigg), 1 \leq k \leq m\bigg) \Big| \mathcal{F}_{n} \bigg)
		\nonumber\\
		&\underset{n \rightarrow \infty}{\longrightarrow} \mathbb{E}\left(f\left(X_{a_k^{-1/2}D_{\infty} }, 1 \leq k \leq m\right) \mid D_{\infty}\right), \quad \text { in probability. }
	\end{align}
\end{prop}

If we want to replace $\delta_n$ by $1$, we will need a slightly stronger condition:
\begin{itemize}
	\item [{\bf(A6)}]
	For some $\gamma_0>0$,
	\[
	\mathbb{E} \bigg( \sum_{u\in\mathcal{N}(1)} \left|V(u)\right|^{2+\gamma_0}e^{-V(u)} \bigg)<\infty.
	\]
\end{itemize}

The assumption {\bf(A6)} says that the random walk $S_n$ has finite $(2+\gamma_0)$-th moment, which
implies that  $\delta_n-1 = o(n^{-\varepsilon_0})$ with some $\varepsilon_0>0$ according to \cite{GX}.
We summarize the result of \cite[Theorem 2.7]{GX} as follows:
\begin{lemma}\label{lemma14}
	 If {\bf(A1)}--{\bf(A6)} hold, then there exists $\varepsilon_0>0$ such that
	\[
	\lim_{n\to\infty} n^{\varepsilon_0} \left|\delta_n -1\right|=0.
	\]
\end{lemma}
Combining Theorem \ref{thm1}, Proposition \ref{prop5} and Lemma \ref{lemma14}, we immediately get
the following theorem:

\begin{thrm}
	Assume that {\bf(A1)}--{\bf(A6)} hold. Then for all $ m \geq 1, a_{1}, \ldots, a_{m} \in[1, \infty)$ and bounded continuous $f: \mathbb{R}^{m} \rightarrow \mathbb{R}$, we have
	\begin{align*}
		&\mathbb{E}\bigg(f\bigg(\sqrt{n} \bigg( D_{\infty}-D_{\lceil a_k n\rceil}+\frac{\log n}{\sqrt{2\pi \sigma^2 \lceil a_k n\rceil}}D_\infty \bigg), 1 \leq k \leq m \bigg) \Big| \mathcal{F}_{n} \bigg)
		\nonumber\\
		&\underset{n \rightarrow \infty}{\longrightarrow} \mathbb{E}\left(f\left(X_{a_k^{-1/2}D_{\infty} }, 1 \leq k \leq m\right) \mid D_{\infty}\right), \quad \text { in probability. }
	\end{align*}
\end{thrm}

The main idea of this paper is a modification of that of \cite{MP}.  To get the fluctuation of $D_\infty-D_{\lceil an \rceil}$, we choose a level $\gamma_n$ and define a quantity
$D_m^{\lceil an\rceil, \gamma_n}$, for $m\ge \lceil an \rceil$,  which roughly takes care of the contributions to $D_m$ by the
paths that stay above the level $\gamma_n$ between generations $\lceil an \rceil$ and $m$. We first show that $D_m^{\lceil an\rceil, \gamma_n}$ converges to a limit $D_\infty^{\lceil an\rceil, \gamma_n}$ as $m\to\infty$ and get a rate of convergence for $D_{\lceil an\rceil}^{\lceil an\rceil, \gamma_n}$ as $n\to\infty$, see Lemma \ref{lemma2}. Then we analyze the contribution
of $D_\infty^{\lceil an\rceil, \gamma_n}$ to the limit behavior of $D_{\lceil an \rceil}$ in Proposition \ref{prop3} .  For contributions to $D_\infty$
by the collection $\mathcal{L}^{\lceil an\rceil, \gamma_n}$ of particles  $x$ with $|x|>\lceil an\rceil$, $V(x)<\gamma_n$ and $\min_{j\in [an, |x|-1]\cap \mathbb{Z}}V(x_j)\ge \gamma_n$,
we separate $\mathcal{L}^{\lceil an\rceil, \gamma_n}$ into two sets $\mathcal{L}_{good}^{\lceil an\rceil, \gamma_n}$ and $\mathcal{L}_{bad}^{\lceil an\rceil, \gamma_n}$ and look at their respective contributions  to the limit behavior of $D_\infty$,
see \eqref{Decomposition-2} in which $F_{good}^{\lceil an \rceil, \gamma_n}$ represents the contributions by $\mathcal{L}_{good}^{\lceil an\rceil, \gamma_n}$  and
$F_{bad}^{\lceil an \rceil, \gamma_n}$
represents the contributions by $\mathcal{L}_{bad}^{\lceil an\rceil, \gamma_n}$. We show in Proposition \ref{prop1} that  $F_{bad}^{\lceil an \rceil, \gamma_n}$
is asymptotically negligible. For the contributions by $\mathcal{L}_{good}^{\lceil an\rceil, \gamma_n}$,
we define a sequence of random variables $\widehat{N}_{good}^{\lceil an \rceil, \gamma_n}$
(see \eqref{e:rs5}). By using the branching property and the  tail behavior of $D_\infty$,
we show in Proposition \ref{prop4} that
$\sqrt{n}(F_{good}^{\lceil an \rceil,\gamma_n}-\widehat{N}_{good}^{\lceil an \rceil, \gamma_n})$
converges in distribution to $c^*X_{a^{-1/2}D_\infty}$
with $c^*$ being the positive constant defined in \eqref{Limit-of-R(y)/y} below,
 which leads to the main result.

Although the general approach of this paper is similar to that of  \cite{MP},  adapting it to the case of branching random walk is pretty challenging. In \cite{MP}, the continuity of the sample paths of Brownian motion makes things a lot easier. For instance, the counterpart of $\widehat{N}_{good}^{\lceil an \rceil, \gamma_n}$ in the branching Brownian motion case takes care the contributions by the particles that hit a certain level at some  time after $at$ due to the continuity of Brownian motion. The main difficulty in the case of branching random walks is that a branching random walk can jump across
the level and one needs to take care of the landing positions of the particles after crossing the level. This leads to many complications and many subtle modifications are needed to actually carry out the program.

\section{Preliminaries}\label{Sec2}

We will use $f(x)\lesssim g(x),\ x\in E,$  to denote that there exists a constant $C$  independent of $x\in E$ such that
\[
f(x)\leq Cg(x), \quad x\in E.
\]
We will use $f(x)\asymp g(x), x\in E$ to denote $f(x)\lesssim g(x), x\in E$ and $g(x)\lesssim f(x), x\in E$.

 \subsection{Spine decomposition}\label{ss:spine}
Define a random walk
$\{(S_n)_{n\geq 0}, \mathbf{P}\}$ such that
for any $n\in\mathbb{N}$ and measurable function  $g: \R^n \to [0,\infty)$,
\begin{align}\label{step_1}
	\mathbb{E} \bigg( \sum_{x\in\mathcal{N}(n)} g\left(V(x_1),...,V(x_n)\right) \bigg) = \mathbf{E}\left(e^{S_n} g\left(S_1,...,S_n\right)\right),
\end{align}
where $\mathbf{E}$ stands for expectation with respect to $\mathbf{P}$ and
for $x\in \mathcal{N}(n)$ and $j\leq n$,  $x_j$ denotes the ancestor of $x$ in the
$j$th generation.
\eqref{step_1} is also known as the  many-to-one formula.
See \cite[Theorem 1.1]{Shi} for more information about the random walk $\{S_n, n\geq 0\}$.
By taking $n=1$, $g(x)= xe^{-x}$ and $g(x)= x^2e^{-x}$ respectively in \eqref{step_1}, we get
that {\bf(A1)} and {\bf(A2)} imply that
$\mathbf{E} S_1 = 0, \sigma^2 = \mathbf{E} S_1^2 < \infty.$
For any  $y\in \R$, we use
$\mathbf{P}_y$ to denote the law of $\{y+S_n, n\geq 0\}$ and $\mathbf{E}_y$ to denote the expectation with respect to $\mathbf{P}_y$.
Note that, under  $\mathbf{P}_y$, $\{y+S_n, n\geq 0\}$ is a random walk starting from $y$.

We define a probability $\Q$ such that for all $n\geq 0$,
 $$\frac{\mathrm{d}\Q}{\mathrm{d}\mathbb{P}}\bigg|_{\mathcal{F}_n}:= W_n.$$
Denote by $\widehat{L}$ the  law of $L$ under $\Q$. Lyons \cite{Lyons} gave the following description of
the law of the branching random walk under $\Q$: there is a spine process denoted by  $\{w_n\}_{n\ge 0}$ with $w_0 =\emptyset$ and the initial position of the spine  is $V(w_0)=0$. At generation $n=1$, $w_0$ dies and splits into a random number of offspring with law $\widehat{L}$.  Choose one offspring $x$ from all the offspring of $w_0$ with probability proportional to $e^{-V(x)}$, and call it $w_1$. $w_1$ evolves independently as $w_0$ and the other unmarked offspring evolve  independently as in the original branching random walk. By Lyons \cite{Lyons}, for any $x\in \mathcal{N}(n)$, we have
\begin{align}\label{Many-to-one-1}
	\Q\left(w_n =x |\mathcal{F}_n\right) = \frac{e^{-V(x)}}{W_n}.
\end{align}
Moreover,  the position process $\left\{V(w_n)\right\}_{n\geq 0}$ along the spine
 under $\Q$ is equal in the law to $\left\{S_n\right\}_{n\geq 0}$ defined in \eqref{step_1}.
Also, for $y\in \R$, we will use  $\Q_y$ to denote the counterpart of $\Q$ in the case of branching random walk with
the initial particle located at $y$.

Let $\tau^+:=\inf\{k\geq 1: S_k \geq 0\}$. Define the renewal function $R(y)$ by
\begin{align}\label{Rwnewal-Function}
	R(y):=
		\mathbf{E}\Big(\sum_{j=1}^{\tau^+-1} 1_{\{S_j \geq -y\}}\Big), \quad y\in \R.
\end{align}
Using the facts that $\mathbf{E} S_1 = 0$ and $\mathbf{E}S_1^2 < \infty$,
one can easily get that  (see, for example, \cite[Section 2.2]{AESZ})
 $R(0)=1,$ $R(y)=0$ for $y<0$ and
\begin{align}\label{step_27}
	R(y) \asymp (1+y), \quad y\ge 0,
\end{align}
and the limit
\begin{align}\label{Limit-of-R(y)/y}
	c^*:= \lim_{y\to+\infty} \frac{R(y)}{y}
\end{align}
exists in $(0, +\infty)$.  According to \cite[(2.4)]{AESZ}, we also have
\begin{align}\label{step_52}
R(y)= \sum_{k=0}^\infty
\mathbf{P}(|H_k|\leq y),
\end{align}
where $H_k:=S_{\sigma_k^-}$ with $\sigma_0^-:=0$ and  $\sigma_k^-:= \inf\left\{i> \sigma_{k-1}^-: S_i <\min_{0\leq j\leq \sigma_{k-1}^-} S_j \right\}$.  For $y\geq 0$, define
$\tau_{-y}^{-,H}:= \inf\{k\geq 1: H_k < -y\}$ and $\tau_{-y}^-:= \inf\{n\geq 1: S_n <-y \}$.
Then we can rewrite \eqref{step_52} as
\begin{align}
	R(y)= \sum_{k=0}^\infty
	\mathbf{P}\left( H_k \geq -y\right)=  \sum_{k=0}^\infty  \mathbf{P}\left(\tau_{-y}^{-,H} > k\right) = \mathbf{E} \left(\tau_{-y}^{-,H}\right).
\end{align}
Note that
$\mathbf{E}(H_1)\in (0,\infty)$
(see \cite[Lemma A.4.(a)]{BIM}) and that
 $H_k - k \mathbf{E}(H_1)$ is a martingale.
Thus, combining the optional sampling theorem and the fact that
$H_{\tau_{-y}^{-, H}} = S_{\tau_{-y}^-}$,  we obtain that
\begin{align}\label{step_53}
(-\mathbf{E} H_1) \mathbf{E} \left(\tau_{-y}^{-,H}\right) = \mathbf{E} \left(-H_{\tau_{-y}^{-, H}} \right) \quad \Longleftrightarrow\quad R(y) \mathbf{E} (|H_1|) = -\mathbf{E} \left(S_{\tau_{-y}^-}\right)= y - \mathbf{E}_y \left(S_{\tau_{0}^-}\right).
\end{align}
By \cite[the first paragraph in the proof of Lemma 2.1]{AESZ},  we  have
$c^*= \mathbf{E} (|H_1|) ^{-1}$.
Note that, as a consequence of {\bf(A4)}, we have
$\mathbf{E} \left((-S_1)_+^3\right)< \infty$.
Thus, by \cite[Lemma A.4.(d)]{BIM},  there exists an $\alpha^*\in (0,\infty)$ such that
\begin{align}\label{Second-Exp-of-R(y)}
	\lim_{y\to+\infty} \left(R(y) -c^*y\right) = \alpha^*.
\end{align}
One can easily check that
\begin{align}\label{Expectation}
	R(y)=
		\mathbf{E}\left(R(S_1 +y)1_{\{S_1\geq -y\}}\right), \quad y\ge 0.
\end{align}
Hence the sequence of random variables
$$
D_n^{-y}:= \sum_{x\in\mathcal{N}(n)}
R\left(V(x) +y\right)e^{-V(x)}1_{\left\{\min_{j\leq n} V(x_j)\geq -y \right\}}
$$
is a non-negative $\P$-martingale with respect to $\{\mathcal{F}_n\}_{n\ge 0}$
with $\mathbb{E}\left( D_n^{-y}\right)= R(y)$ for all $n\geq 0$. Define a new probability measure $\Q^{-y}$ such that for all $n\geq 0$,
\begin{equation}\label{def-Q-y}
\frac{\mathrm{d} \Q^{-y}}{\mathrm{d} \mathbb{P}}\bigg|_{\mathcal{F}_n}:= \frac{D_n^{-y}}{R(y)}.
\end{equation}
Similar to the spine decomposition under $\Q$, we can also describe the spine decomposition for the branching random walk under $\Q^{-y}$ with a spine denoted by $\{w_n^{-y}\}_{n\ge 0}$ and with spatial displacement following the law of the random walk $\{S_n\}$ conditioned to stay in $[-y, +\infty)$:
there is a spine process denoted by  $\{w_n^{-y}\}_{n\ge 0}$ with $w_0^{-y} =\emptyset$ and the initial position of the spine  is $V(w_0^{-y})=0$. At generation $n=1$, $w_0^{-y}$ dies and splits into a random number of offspring
according to the law of $L$ under $\Q^{-y}$.
Choose one offspring $x$ from all the offspring of $w_0^{-y}$ with probability proportional to
$R(V(x)+y)e^{-V(x)}1_{\{V(x)\geq -y\}}$,
and call it $w_1^{-y}$. At generation $n=2$, given $V(w_1^{-y})$,
$w_1^{-y}$ gives birth to a point process according to $(L, \Q^{-(V(w_1^{-y})+y)})$ and again choose one offspring $x$ from all the offspring of $w_1^{-y}$ with probability proportional to
$R(V(x)+y)e^{-V(x)}1_{\{V(x)\geq -y\}}$
named $w_2^{-y}$. The other unmarked offspring evolve  independently as in the original branching random walk. The procedure goes on.
According to \cite[Fact 3.2]{AESZ} or \cite[Section 2.2]{Chen},  for $x\in\mathcal{N}(n)$,
\begin{align}\label{Many-to-one-2}
	\Q^{-y}\left(w_n^{-y} =x |\mathcal{F}_n\right) = \frac{R\left(V(x)+y\right)e^{-V(x)}1_{\left\{\min_{j\leq n} V(x_j)\geq -y \right\}} }{D_n^{-y}}
\end{align}
and the position process $\left(V(w_n^{-y}) \right)_{n\geq 1}$ along the spine is equal in law to
 $\{S_n\}_{n\geq 1}$
 conditioned to stay in $[-y,+\infty)$.

\subsection{Elementary properties for centered random walk}
\begin{lemma}\label{lemma1}
	(i) For all $a \geq 0$ and $n\geq 1$, it holds that
	\[
		\mathbf{P}_a
	\Big(\min_{j\leq n} S_j \geq 0 \Big) \lesssim \frac{(1+a)}{\sqrt{n}}.
	\]
	(ii) For all $a, u\geq  0$, $b>0$ and $n\geq 1$, it holds that
	\[
		\mathbf{P}_a
	\Big(\min_{j\leq n} S_j \geq 0, u\leq S_n \leq b+u\Big)\lesssim \frac{(b+1)(b+u+1)(a+1)}{\sqrt{n^3}}.
	\]
	(iii) For all $a,b\geq 0$, it holds that
	\[
	\sum_{n=0}^\infty
		\mathbf{P}_a
	\Big(\min_{j\leq n} S_j \geq 0, S_n \leq b\Big)\lesssim (1+b)\left(1+\left(a\land b\right)\right).
	\]
	Here $a\land b:= \min\{a,b \}.$
	\newline
	(iv)
	For any $\lambda > 0$, there exists a constant $ C_1(\lambda)>0$ such that
	\[
	\sum_{k=1}^\infty
		\mathbf{E}_a
	\Big(e^{-\lambda S_k}1_{\left\{\min_{j\leq k}S_j \geq 0 \right\}}\Big)\leq C_1(\lambda),\quad a\ge 0.
	\]
\end{lemma}
\textbf{Proof:} For (i), see \cite[(2.7)]{AE}; for (ii), see \cite[Lemma 2.2]{AESZ}; for (iii) and (iv), see \cite[Lemma B.2 (i) and (iii)]{AE}.
\hfill$\Box$

\begin{lemma}\label{lemma12}
	For all $a \geq 0$ and $n\geq 1$, it holds that
	\[
		\mathbf{E}
	\Big(S_n^2 1_{\left\{\min_{j\leq n} S_j \geq -a \right\}} \Big) \lesssim (1+a)\sqrt{n}.
	\]
\end{lemma}
\textbf{Proof: }
Note that under $\mathbf{P}$,
$\left\{S_n^2 -\sigma^2 n: n\geq 1 \right\}$ is a mean 0 martingale. Thus, by Lemma \ref{lemma1} (i),
\begin{align*}
&\mathbf{E}\Big(S_n^2 1_{\left\{\min_{j\leq n} S_j \geq -a \right\}} \Big) = \sigma^2 n\mathbf{P}\Big(\min_{j\leq n} S_j \geq -a \Big)+ \mathbf{E}\Big(\left(S_n^2-\sigma^2 n \right)1_{\left\{\min_{j\leq n} S_j \geq -a \right\}} \Big)\\
	& \lesssim (1+a)\sqrt{n}- \mathbf{E}\Big(\left(S_n^2-\sigma^2 n \right)1_{\left\{\min_{j\leq n} S_j < -a \right\}} \Big)\\
	& = (1+a)\sqrt{n} + \sum_{\ell=1}^n \mathbf{E}\Big(\left(\sigma^2 n-S_n^2\right)1_{\left\{\min_{j\leq \ell-1} S_j \geq -a \right\}}1_{\{S_\ell <-a\}} \Big)\\
	& = (1+a)\sqrt{n} + \sum_{\ell=1}^n \mathbf{E}\Big(\left(\sigma^2\ell -S_\ell^2\right)1_{\left\{\min_{j\leq \ell-1} S_j \geq -a \right\}}1_{\{S_\ell <-a\}} \Big)\\	& \leq (1+a)\sqrt{n} + \sigma^2 \sum_{\ell=1}^n \ell\, \mathbf{P}\Big(\min_{j\leq \ell-1} S_j \geq -a, S_\ell <-a \Big).
\end{align*}
Using Lemma \ref{lemma1} (i) again, we get
\begin{align*}
&\sum_{\ell=1}^n \ell\, \mathbf{P}\Big(\min_{j\leq \ell-1} S_j \geq -a, S_\ell <-a \Big)=\sum_{\ell=1}^n \ell\,  \mathbf{P}\Big(\min_{j\leq \ell-1} S_j \geq -a \Big) - \sum_{\ell=1}^n \ell\, \mathbf{P}\Big(\min_{j\leq \ell} S_j \geq -a\Big) \\
	& =  1-(n+1)\mathbf{P}\Big(\min_{j\leq n}S_j \geq -a \Big)+\sum_{\ell=1}^n \mathbf{P}\Big(\min_{j\leq \ell} S_j \geq -a\Big)\\
	& \lesssim 1+  (1+a)\sum_{\ell=1}^n \frac{1}{\sqrt{\ell}}\leq  1+  (1+a) \int_0^n \frac{1}{\sqrt{x}}\mathrm{d}x = 1+ 2 (1+a) \sqrt{n}\lesssim (1+a)\sqrt{n}.
\end{align*}
Combining the two displays above, we get the desired conclusion.
\hfill$\Box$

\begin{lemma}\label{lemma3}
If  $X$ and $\widetilde{X}$ are non-negative random variables such that
$$
\mathbb{E}\Big(X \left(\log_+ X\right)^3\Big)+\mathbb{E}\Big(\widetilde{X} \left(\log_+ \widetilde{X}\right)^2\Big)<\infty,
$$
	then
	\[
	\mathbb{E}\bigg(X \Big(\log_+ \left(\widetilde{X}+ X\right)\Big)^3\bigg) + \mathbb{E}\bigg(\widetilde{X} \Big(\log_+ \left(\widetilde{X}+ X\right)\Big)^2 \bigg)<\infty.
	\]
\end{lemma}
\textbf{Proof:}
By the trivial inequality $\log_+\left(x+y\right)\leq  \log_+\left(2x\right)+\log_+\left(2y\right)$,  we only need to show that
\[
\mathbb{E}\Big(X \Big(\log_+ \widetilde{X}\Big)^3\Big) + \mathbb{E}\Big(\widetilde{X} \Big(\log_+ X \Big)^2 \Big)<\infty.
\]
For this, it suffices to prove that for any $x, \widetilde{x}>0$,
\begin{align}
	x\left(\log_+ \widetilde{x}\right)^3 \leq  8x\left(\log_+ x\right)^3 + 2\widetilde{x}\left(\log_+ \widetilde{x}\right)^2,\label{step_29}\\
	\widetilde{x} \left(\log_+ x \right)^2 \leq   4 \widetilde{x} \left(\log_+ \widetilde{x} \right)^2 + 2x\left(\log_+ x\right).\label{step_33}
\end{align}
We will only prove \eqref{step_29}, the proof of \eqref{step_33} is similar.
Assume that $\widetilde{x}\geq 1$. If $\widetilde{x}\leq x^2$, then $x\left(\log \widetilde{x}\right)^3 \leq x\left(\log \left(x^2\right)\right)^3= 8x \left(\log x\right)^3$.
If $\widetilde{x} \geq x^2$, then by trivial inequality
$$
\log \widetilde{x}\leq 2\sqrt{\widetilde{x}}, \quad \widetilde{x}\geq 1,
$$
we have
$x\left(\log \widetilde{x}\right)^3  \leq \sqrt{\widetilde{x}}\left(\log \widetilde{x}\right)^3 \leq 2\widetilde{x}\left(\log \widetilde{x}\right)^2$.
The proof is complete.
\hfill$\Box$

\subsection{Moment estimates for the truncated martingales}

For $u\in \mathcal{N}(n),\Omega(u):= \left\{ v\in\mathcal{N}(n): v\neq u: v> u_{n-1} \right\}$
 denotes the set of  siblings of $u$. For $\kappa,y > 0$ and $m\in\mathbb{N}$, define
\begin{align*}
	B_{m,\kappa}& := \Big\{x\in \mathcal{N}(m): \forall 1\leq j \leq m, \\
	& \qquad \sum_{u\in \Omega(x_j)} \left(1+\left(V(u)-V(x_{j-1})\right)_+\right)e^{-\left(V(u)-V(x_{j-1})\right)} \leq \kappa e^{\left(V(x_{j-1})+y\right)/2} \Big\},\\
	D_{m,\kappa}^{-y}&:=	\sum_{x\in\mathcal{N}(m)}
	R\left(V(x) +y\right)e^{-V(x)}1_{\left\{\min_{j\leq m} V(x_j)\geq -y \right\}}1_{\{x\in B_{m,\kappa} \}}.
\end{align*}
\begin{lemma}\label{lemma4}
	There exists a decreasing function $h: \R_+ \to \R_+$ satisfying $\lim_{z\to+\infty} h(z)=0$ such that for all $y,\kappa > 0$ and $m\in\mathbb{N}$,
	\[
    0\le \mathbb{E}\left(D_m^{-y} - D_{m,\kappa}^{-y}\right)\lesssim h(\kappa).
	\]
\end{lemma}
\textbf{Proof:}
The first inequality is trivial, so we only prove the second.   For $j\geq 1$, set
\begin{align*}
	E_j(y,\kappa ):= \bigg\{\sum_{u\in \Omega(w_j^{-y})\cup\{ w_j^{-y}\}} \left(1+\left(V(u)-V(w_{j-1}^{-y})\right)_+\right)e^{-\left(V(u)-V(w_{j-1}^{-y})\right)} > \kappa e^{\left(V(w_{j-1}^{-y})+y\right)/2} \bigg\}.
\end{align*}
It follows from  \eqref{Many-to-one-2} that
\begin{align}\label{step_36}
	&\mathbb{E}\left(D_m^{-y} - D_{m,\kappa}^{-y}\right)
	= \mathbb{E}\bigg(	\sum_{x\in\mathcal{N}(m)}
	R\left(V(x) +y\right)e^{-V(x)}1_{\left\{\min_{j\leq m} V(x_j)\geq -y \right\}}1_{\{x\notin B_{m,\kappa} \}}\bigg)\nonumber\\
	& = \mathbb{E}_{\Q^{-y}}\bigg(\frac{R(y)}{D_n^{-y}}	\sum_{x\in\mathcal{N}(m)}
	R\left(V(x) +y\right)e^{-V(x)}1_{\left\{\min_{j\leq m} V(x_j)\geq -y \right\}}1_{\{x\notin B_{m,\kappa} \}}\bigg)\nonumber\\
	&= R(y)\mathbb{E}_{\Q^{-y}}\bigg(	\sum_{x\in\mathcal{N}(m)}
	\Q^{-y}\left(w_m^{-y} =x |\mathcal{F}_m\right) 1_{\{x\notin B_{m,\kappa} \}}\bigg)\nonumber\\
	&= R(y)\Q^{-y} \left( w_m^{-y}\notin B_{m,\kappa}\right)\leq R(y)\sum_{j=1}^\infty \Q^{-y}\left( E_j(y,\kappa)\right).
\end{align}
By the Markov property, for any $z\geq -y,$
\begin{align*}
	&\Q^{-y}\left( E_j(y,\kappa) \big|V( w_{j-1}^{-y}) = z \right)=
	\Q^{-y-z}\bigg(	\sum_{u\in\mathcal{N}(1)}	\left(1+\left(V(u)\right)_+\right)e^{-V(u)} > \kappa e^{\left(z+y\right)/2} \bigg)\nonumber\\
	 & = \mathbb{E}\bigg(\frac{	\sum_{u\in \mathcal{N}(1)}	R(V(u)+z+y)e^{-V(u)}1_{\{V(u)\geq -y -z\}}}{R(z+y)} 1_{\left\{	\sum_{u\in\mathcal{N}(1)}	\left(1+\left(V(u)\right)_+\right)e^{-V(u)} > \kappa e^{\left(z+y\right)/2} \right\}}\bigg)
	 \nonumber\\
	& \leq 	\E 	\bigg(\frac{\sum_{u\in\mathcal{N}(1)}	R(V(u)+z+y)e^{-V(u)}}{R(z+y)} 1_{\left\{ \sum_{u\in\mathcal{N}(1)} \left(1+\left(V(u)\right)_+\right)e^{-V(u)} > \kappa e^{\left(z+y\right)/2} \right\}}\bigg).
\end{align*}
Using \eqref{step_27}, we have
\[
\frac{ R(V(u)+z+y)}{R(z+y)}\lesssim \frac{\left(V(u)\right)_+ +z+y+1}{z+y+1} = 1+ \frac{\left(V(u)\right)_+}{z+y+1}.
\]
Thus,
\begin{align}\label{step_28}
	&\Q^{-y}\left( E_j(y,\kappa) \big| V(w_{j-1}^{-y}) = z \right)\nonumber \\
	&\lesssim \E \bigg( \sum_{u\in\mathcal{N}(1)}	\left(1+ \frac{\left(V(u)\right)_+}{z+y+1}\right) e^{-V(u)} 1_{\left\{	\sum_{u\in\mathcal{N}(1)}
	\left(1+\left(V(u)\right)_+\right)e^{-V(u)} > \kappa e^{\left(z+y\right)/2} \right\}}\bigg)\nonumber\\
	& =  \E \bigg( \bigg(W_1 + \frac{\widetilde{W}_1}{z+y+1}\bigg) 1_{\left\{ W_1 +\widetilde{W}_1 > \kappa e^{\left(z+y\right)/2} \right\}}\bigg).
\end{align}
Since the law of $(W_1, \widetilde{W}_1)$ is independent of $z$, we deduce from \eqref{step_28} that
\begin{align}\label{step_35}
	&\Q^{-y}\left( E_j(y,\kappa)  \right)\lesssim
	\left(\mathbb{E}_{\mathbb{Q}^{-y}}\otimes \mathbb{E}\right)
	\bigg( \bigg(W_1 + \frac{\widetilde{W}_1}{V(w_{j-1}^{-y})+y+1}\bigg) 1_{\big\{ W_1 +\widetilde{W}_1 > \kappa e^{\left(V(w_{j-1}^{-y})+y\right)/2} \big\}}\bigg)\nonumber\\
    & = \E \bigg( \mathbb{E}_{\mathbb{Q}^{-y}}\bigg( \bigg(z_1 + \frac{\widetilde{z}_1}{V(w_{j-1}^{-y})+y+1}\bigg) 1_{\big\{ V(w_{j-1}^{-y})+y < 2\log \left(\frac{z_1 +\widetilde{z}_1 }{\kappa}\right) \big\}}\bigg)\bigg|_{{z_1 = W_1, \widetilde{z}_1 = \widetilde{W}_1}}\bigg).
\end{align}
Here under $\mathbb{Q}^{-y}\otimes \P$,  $\big(W_1, \widetilde{W}_1\big)$ is independent of $V(w_{j-1}^{-y}).$ Next, note that the law of $V(w_{j}^{-y})$ under $\mathbb{Q}^{-y}$ is equal
to the law of the random walk $S_j$ conditioned to stay in  $[-y, +\infty)$.
Summing $j$ from $1$ to $\infty$, and using \eqref{Expectation} and  the fact that  $R(y)\lesssim 1+y$, we get
\begin{align}\label{step_34}
& R(y)\sum_{j=1}^\infty \mathbb{E}_{\mathbb{Q}^{-y}}\bigg( \bigg(z_1 + \frac{\widetilde{z}_1}{V(w_{j-1}^{-y})+y+1} \bigg) 1_{\big\{V( w_{j-1}^{-y})+y < 2\log \left(\frac{z_1 +\widetilde{z}_1 }{\kappa}\right) \big\}}\bigg)\nonumber\\
 & = \sum_{j=0}^\infty \mathbf{E}\bigg( R(S_{j}+y)1_{\{\min_{\ell \leq j } S_\ell \geq -y \}}\bigg(z_1 + \frac{\widetilde{z}_1}{S_{j}+y+1}\bigg)
  1_{\big\{ S_{j}+y < 2\log \left(\frac{z_1 +\widetilde{z}_1 }{\kappa}\right) \big\}}\bigg)\nonumber\\
 & \lesssim \sum_{j=0}^\infty \mathbf{E}\Big( \left(z_1\left(S_j+y+1\right) + \widetilde{z}_1\right) 1_{\{\min_{\ell \leq j } S_\ell \geq -y \}} 1_{\big\{ S_{j}+y < 2\log \left(\frac{z_1 +\widetilde{z}_1 }{\kappa}\right) \big\}}\Big)\nonumber\\
	& \leq \bigg(z_1\bigg(1+ 2\log_+ \bigg(\frac{z_1 +\widetilde{z}_1 }{\kappa}\bigg)  \bigg) + \widetilde{z}_1 \bigg)
	\sup_{y\in \mathbb{R}}\sum_{j=0}^\infty \mathbf{P}\bigg( \min_{\ell \leq j} S_\ell \geq -y, S_j +y <2\log_+ \bigg(\frac{z_1+\widetilde{z}_1 }{\kappa}\bigg)   \bigg)\nonumber\\ &=
\bigg(z_1\bigg(1+ 2\log_+ \bigg(\frac{z_1 +\widetilde{z}_1 }{\kappa}\bigg)  \bigg) + \widetilde{z}_1 \bigg) F\bigg(2\log_+\bigg(\frac{z_1 +\widetilde{z}_1}{\kappa}\bigg)\bigg),
\end{align}
where
$F(x):=\sup_{y\in\R}\sum^\infty_{j=0} \mathbf{P}_y(\min_{l\leq j}S_l\geq 0, S_j<x)$.
Taking $\P$-expectation in \eqref{step_34},
and combining the result with \eqref{step_36} and \eqref{step_35}, we get
\begin{align}
	&\mathbb{E}\left(D_m^{-y} - D_{m,\kappa}^{-y} \right)
	\leq R(y)\sum_{j=1}^\infty \Q^{-y}\left(E_j(y,\kappa) \right)\nonumber \\ &\lesssim
	\E  \bigg(\bigg(W_1\bigg(1+ 2\log_+ \bigg(\frac{W_1 +\widetilde{W}_1 }{\kappa}\bigg)  \bigg) + \widetilde{W}_1 \bigg) F\bigg(2\log_+ \bigg(\frac{W_1 +\widetilde{W}_1 }{\kappa}\bigg)\bigg)\bigg)=:h(\kappa).
\end{align}
It follows from Lemma \ref{lemma1} (iii) that $F(x)\lesssim (1+x)^2$ for all $x\geq 0$. By Lemma \ref{lemma3}, we know that $h(\kappa)$ is finite for all $\kappa>0$ and $h$ is a decreasing function. Since $F(0)=0$, we have for $\kappa >1$,
\begin{align*}
	& h(\kappa) \leq \mathbb{E}\left(\left(W_1\left( 1+ 2\log_+\left(W_1 +\widetilde{W}_1\right)  \right)+\widetilde{W}_1\right)F\left(\log_+(W_1 +\widetilde{W}_1)\right)1_{\left\{W_1 +\widetilde{W}_1 > \kappa \right\}}\right) \\
	& \lesssim \mathbb{E}\left(\left(W_1\left( 1+ \log_+(W_1 +\widetilde{W}_1) \right) + \widetilde{W}_1\right)\left( 1+ \log_+(W_1 +\widetilde{W}_1) \right)^2 1_{\left\{W_1 +\widetilde{W}_1 > \kappa \right\}} \right).
\end{align*}
The right hand side above tends to 0 as $\kappa \to\infty$. The proof is complete.
\hfill$\Box$

\begin{lemma}\label{lemma5}
	For all $y\geq 0, \kappa \geq 1$ and $m\geq 0$,
	\[
	\mathbb{E}\left( \left( D_{m,\kappa}^{-y} \right)^2	\right)\lesssim \kappa e^{y}.
	\]
\end{lemma}
\textbf{Proof:} Using \eqref{Many-to-one-2} and the fact that $D_{m,\kappa}^{-y}\leq D_m^{-y}$, we get
\begin{align}\label{step_39}
	&\mathbb{E}\left(\left( D_{m,\kappa}^{-y} \right)^2\right)  = \mathbb{E}\bigg( D_{m,\kappa}^{-y}	\sum_{x\in\mathcal{N}(m)}
	R\left(V(x) +y\right)e^{-V(x)}1_{\left\{\min_{j\leq m} V(x_j)\geq -y \right\}}1_{\{x\in B_{m,\kappa} \}} \bigg)\nonumber\\
	& = R(y)\mathbb{E}_{\Q^{-y}}\bigg( D_{m,\kappa}^{-y}	\sum_{x\in\mathcal{N}(m)}
	\Q^{-y}\left(w_m^{-y} =x |\mathcal{F}_m\right) 1_{\{x\in B_{m,\kappa} \}} \bigg)\nonumber\\
	& \leq R(y) \mathbb{E}_{\Q^{-y}}\left( D_{m}^{-y} 1_{\{w_m^{-y}\in B_{m,\kappa} \}} \right).
\end{align}
By the spine decomposition, we have
\begin{equation}\label{spine-decom-1}
D_{m}^{-y}= R\left(V\left(w_m^{-y}\right) + y\right)e^{-V\left(w_m^{-y} \right)}+ \sum_{\ell=1}^{m} \sum_{z\in\Omega\left(w_\ell^{-y}\right)}R\left(V(z)+y\right)e^{-V(z)}.
\end{equation}
Since
\[
R(x+y)\leq R(x_+ +y)\lesssim (1+x_+ +y)\leq (1+x_+)(1+y), \quad y\geq 0, \,x\in\R,
\]
we get, for any $1\leq \ell \leq m$,
$$
    R\left(V(z)+y\right)\lesssim \left(1+y+V\left(w_{\ell-1}^{-y}\right)\right)\Big(1+\left( V(z)-V\left(w_{\ell-1}^{-y}\right)\right)_+\Big), \quad z\in \Omega\left(w_\ell^{-y}\right).
$$
Thus, for $w_m^{-y}\in B_{m,\kappa}$ and $1\leq \ell \leq m$,
\begin{align}\label{step_37}
		&\sum_{z\in\Omega\left(w_\ell^{-y}\right)}R\left(V(z)+y\right)e^{-V(z)}\nonumber \\&\lesssim \left(1+y+V\left(w_{\ell-1}^{-y}\right)\right)e^{-V\left(w_{\ell-1}^{-y}\right)} \sum_{z\in\Omega\left(w_\ell^{-y}\right)}\Big(1+\left( V(z)-V\left(w_{\ell-1}^{-y}\right)\right)_+\Big)e^{-\left(V(z)-V\left(w_{\ell-1}^{-y}\right)\right)} \nonumber\\ & \leq \kappa e^{y/2}  \left(1+y+V\left(w_{\ell-1}^{-y}\right)\right)e^{-V\left(w_{\ell-1}^{-y}\right)/2} .
\end{align}
Combining this with \eqref{spine-decom-1} and  the fact that $V(w_m^{-y})+y\geq 0$, we get
\begin{align}\label{step_38}
	&R(y)\mathbb{E}_{\Q^{-y}}\left( D_{m}^{-y} 1_{\{w_m^{-y}\in B_{m,\kappa} \}} \right)\nonumber\\
& \lesssim R(y)\mathbb{E}_{\Q^{-y}}\bigg(R\left(V\left(w_m^{-y}\right) + y\right)e^{-V\left(w_m^{-y} \right)}+ \sum_{\ell =1}^{m}   \kappa e^{y/2}  \left(1+y+V\left(w_{\ell-1}^{-y}\right)\right)e^{-V\left(w_{\ell -1}^{-y}\right)/2}\bigg)\nonumber
	\\& \lesssim  \kappa e^{y}R(y)\sum_{\ell =0}^m \mathbb{E}_{\Q^{-y}}\left( \left(1+y+V\left(w_\ell ^{-y}\right)\right)e^{-\left(y+V\left(w_\ell^{-y}\right)\right)/2}\right).
\end{align}
Since $\left( (1+y)^2 e^{-y/4}\right) \lesssim 1$ on $[0, \infty)$,
we have for all $\ell \geq 0$ and $y\geq 0$,
\begin{align*}
		&R(y)\mathbb{E}_{\Q^{-y}}\left( \left(1+y+V\left(w_\ell^{-y}\right)\right)e^{-\left(y+V\left(w_\ell^{-y}\right)\right)/2}\right)	=
	        \mathbf{E}_y \left(R(S_\ell)\left(1+S_\ell \right)e^{-S_\ell /2}1_{\left\{\min_{j\leq \ell} S_j \geq 0\right\}}\right)\\	& \lesssim \mathbf{E}_y \left(\left(1+S_\ell \right)^2e^{-S_\ell /2}1_{\left\{\min_{j\leq \ell } S_j \geq 0\right\}}\right)\lesssim \mathbf{E}_y \left(e^{-S_\ell /4}1_{\left\{\min_{j\leq \ell} S_j \geq 0\right\}}\right),
\end{align*}
where in the equality we used  \eqref{def-Q-y}.
Applying Lemma \ref{lemma1} (iv) with $\lambda =\frac{1}{4}$, we get that for all $y\geq 0, m\geq 0$,
\begin{align}\label{step_40}
		&R(y)\sum_{\ell =0}^m \mathbb{E}_{\Q^{-y}}\left( \left(1+y+V\left(w_\ell^{-y}\right)\right)e^{-\left(y+V\left(w_\ell^{-y}\right)\right)/2}\right)\nonumber\\ &\lesssim \sum_{\ell=0}^\infty
				\mathbf{E}_y
		\left(e^{-S_\ell /4}1_{\left\{\min_{j\leq \ell} S_j \geq 0\right\}}\right)\leq e^{-y/4} +  C_1(1/4)\lesssim 1.
\end{align}
Combining \eqref{step_39}, \eqref{step_38} and \eqref{step_40}, we get the desired conclusion.
\hfill$\Box$

\subsection{Moment estimate for weighted number of particles hitting $-y$}\label{subsec3.4}

 For $y\geq 0$, $\kappa>0$ and $n,m \in\mathbb{N}:=\{1, 2, \cdots\}$  with $n\le m$, define
\begin{align}\label{step_11}
		    & N_{[n,m]}^y:= \sum_{\ell =n}^m 	\sum_{x\in\mathcal{N}(\ell )}	e^{-V(x)} 1_{\left\{V(x)<-y,\, \min_{j\leq \ell -1} V(x_j)\geq -y \right\}},\nonumber\\
	&\mathcal{A}_\kappa^{y} := \bigg\{x\in \mathbb{T}: \forall 1\leq j \leq |x|, \sum_{u\in \Omega(x_j)} \left(1+\left(V(u)-V(x_{j-1})\right)_+\right)e^{-\left(V(u)-V(x_{j-1})\right)} \leq \kappa e^{\left(V(x_{j-1})+y\right)/2} \bigg\},\nonumber \\ &
		N_{[n,m],\kappa}^y:= \sum_{\ell =n}^m	\sum_{x\in\mathcal{N}(\ell)}	e^{-V(x)} 1_{\left\{V(x)<-y,\, \min_{j\leq \ell -1} V(x_j)\geq -y \right\}}1_{\left\{ x\in \mathcal{A}_{\kappa}^y\right\}}.
\end{align}
We will use the notation $N_{[1,\infty)}^y:= \lim_{m\to\infty} N_{[1,m]}^y$ and  $N_{[1,\infty),\kappa}^y:= \lim_{m\to\infty} N_{[1,m],\kappa}^y$.

\begin{lemma}\label{lemma7}
	(i) For any $y\geq 0$,
	\[
	\mathbb{E}\left(N_{[1,\infty)}^y\right) =1.
	\]
	(ii)
	There exists a decreasing function $g:\R_+ \to \R_+$ with $\lim_{z\to+\infty}g(z)=0$ such that
	\[
	\mathbb{E}\left(N_{[1,\infty)}^y - N_{[1, \infty), \kappa}^y\right)\lesssim \frac{g(\kappa)}{\log \kappa},
	\quad y >0, \kappa > 1.
	\]
\end{lemma}
\textbf{Proof: } (i)  By the definition of $\Q$, we have
\begin{align*}
 \mathbb{E}\left(N_{[1,\infty)}^y\right) &= \sum_{k=1}^\infty \mathbb{E}_\Q \Big( \sum_{x\in\mathcal{N}(k)}
	\frac{e^{-V(x)}}{W_k} 1_{\left\{V(x)<-y,\, \min_{j\leq k-1} V(x_j)\geq -y \right\}}\Big)\\
 &=\sum_{k=1}^\infty \mathbb{E}_\Q \Big( \sum_{x\in\mathcal{N}(k)}
	 \Q \Big( w_k=x, V(w_k)<-y, \min_{j\leq k-1} V(w_j)\geq -y |\mathcal{F}_k\Big) \Big)\\
&=\sum_{k=1}^\infty \Q \Big(V(w_k)<-y, \min_{j\leq k-1} V(w_j)\geq -y\Big)=1,
\end{align*}
where in the second equality we used \eqref{Many-to-one-1}.

(ii)
For any $m \in \mathbb{N}$, by the definition of $\Q$ and \eqref{Many-to-one-1}, we have
\begin{align}\label{step_13}
		\mathbb{E}\left(N_{[1,m]}^y - N_{[1, m], \kappa}^y\right)&= \sum_{\ell =1}^m \mathbb{E} \bigg( \sum_{x\in\mathcal{N}(\ell)}e^{-V(x)}1_{\left\{V(x)<-y, \min_{j\leq \ell -1} V(x_j)\geq -y \right\}}1_{\left\{ x\notin \mathcal{A}_{\kappa}^y\right\}}\bigg)\nonumber\\ & = \sum_{\ell=1}^m \Q \Big(  V(w_\ell )<-y, \min_{j\leq \ell-1} V(w_j)\geq -y , w_\ell \notin \mathcal{A}_\kappa^y \Big)\nonumber\\ & \leq \sum_{\ell =1}^\infty \Q \Big(  V(w_\ell)<-y, \min_{j\leq \ell-1} V(w_j)\geq -y , w_\ell \notin \mathcal{A}_\kappa^y \Big).
\end{align}
Since
\[
1_{\left\{w_\ell \notin \mathcal{A}_\kappa^y \right\}} \leq \sum_{q=1}^{\ell} 1_{\left\{ \sum_{u\in \Omega(w_q)} \left(1+\left(V(u)-V(w_{q-1})\right)_+\right)e^{-\left(V(u)-V(w_{q-1})\right)} > \kappa e^{\left(V(w_{q-1})+y\right)/2}\right\}}=: \sum_{q=1}^\ell 1_{G_q},
\]
we have
\begin{align*}
		&\sum_{\ell =1}^\infty \Q \Big(  V(w_\ell)<-y, \min_{j\leq \ell-1} V(w_j)\geq -y , w_\ell \notin \mathcal{A}_\kappa^y \Big)\\ & \leq \sum_{\ell =1}^\infty \sum_{q=1}^\ell  \Q \Big(  V(w_\ell )<-y, \min_{j\leq \ell-1} V(w_j)\geq -y , G_q \Big)\\ & = \sum_{q=1}^\infty \sum_{\ell=q}^\infty \Q \Big(  V(w_\ell )<-y, \min_{j\leq\ell -1} V(w_j)\geq -y , G_q\Big) = \sum_{q=1}^\infty \mathbb{Q}\Big(G_q, \min_{j\leq q -1} V(w_j)\geq -y \Big)	\\ &\leq  \sum_{q=1}^\infty \mathbb{E}_{\mathbb{Q}}\bigg(1_{\{\min_{j\leq q -1} V(w_j)\geq -y\}}\mathbb{Q}\bigg( \sum_{u\in\mathcal{N}(1)}	\left(1+\left(V(u)\right)_+\right)e^{-V(u)} > \kappa e^{\left(z+y\right)/2}\bigg)\bigg|_{z= V(w_{q -1})} \bigg).
\end{align*}
Recalling the definition of $F$ in \eqref{step_34} and using an argument similar to that of the proof of Lemma \ref{lemma4}, we get
\begin{align*}
	&\sum_{q=1}^\infty \mathbb{E}_{\mathbb{Q}}\bigg(1_{\{\min_{j\leq q-1} V(w_j)\geq -y\}}\mathbb{Q}\bigg( \sum_{u\in\mathcal{N}(1)} \left(1+\left(V(u)\right)_+\right)e^{-V(u)} > \kappa e^{\left(z+y\right)/2}\bigg)\bigg|_{z= V(w_{q-1})} \bigg)\\ & =  \mathbb{E}_{\mathbb{Q}}\bigg(\sum_{q=1}^\infty \mathbb{Q}\Big(\min_{j\leq q-1} V(w_j)\geq -y, V(w_{q-1})  +y< 2\log_+ \left(\frac{m_1}{\kappa}\right) \Big)\bigg|_{m_1 = W_1 +\widetilde{W}_1} \bigg)\\& \leq  \mathbb{E}_{\mathbb{Q}}\bigg(\sup_{y\in \mathbb{R}}\sum_{q=1}^\infty \mathbb{Q}\Big(\min_{j\leq q -1} V(w_j)\geq -y, V(w_{q-1})  +y< 2\log_+ \left(\frac{m_1}{\kappa}\right) \Big)\bigg|_{m_1 = W_1 +\widetilde{W}_1} \bigg)\\ &= \mathbb{E}_{\mathbb{Q}}\bigg(F\bigg(2\log_+\bigg(\frac{W_1 +\widetilde{W}_1}{\kappa}\bigg)\bigg) \bigg)=\mathbb{E}\bigg(W_1F\bigg(2\log_+\bigg(\frac{W_1 +\widetilde{W}_1}{\kappa}\bigg)\bigg) \bigg).
\end{align*}
 It follows from Lemma \ref{lemma1}(iii) that $F(x)\lesssim (1+x)^2$ for $x\ge 0$. Note that $F(0)=0$, we have
\begin{align*}
	&	\mathbb{E}\left(N_{[1,\infty)}^y - N_{[1, \infty), \kappa}^y\right) \leq \mathbb{E}
		    \bigg(W_1F\bigg(2\log_+\bigg(\frac{W_1 +\widetilde{W}_1}{\kappa} \bigg)\bigg) 1_{\left\{W_1 +\widetilde{W}_1>\kappa \right\}}\bigg)
	\\ &  \lesssim \mathbb{E} \bigg(W_1\bigg(1+2\log_+\bigg(\frac{W_1 +\widetilde{W}_1}{\kappa}\bigg) \bigg)^2 1_{\left\{W_1 +\widetilde{W}_1 >\kappa \right\}}\frac{\log_+ \left(W_1+\widetilde{W}_1\right)}{\log \kappa}\bigg) =: \frac{g(\kappa)}{\log\kappa}.
\end{align*}
The proof is complete.
\hfill$\Box$

\begin{lemma}\label{lemma8}
	Let $\alpha$ be the constant in {\bf(A5)}. Then
	\[
	\mathbb{E}\Big(\left(N_{[1,\infty) ,\kappa}^y\right)^{1+\alpha}\Big)\lesssim \kappa^\alpha e^{\alpha y},
	\quad y\geq 0, \kappa \geq 1.
	\]
\end{lemma}
\textbf{Proof: } By Lemma \ref{lemma7} (i), we have
\begin{align}\label{step_12}
	\mathbb{E}\big(N_{[1,\infty) ,\kappa}^y\big) \leq \mathbb{E}\big(N_{[1,\infty)}^y\big)=1.
\end{align}
Define $\tau_{-y}^- := \inf\left\{k\geq 1: V(w_k) < -y\right\}$.
For $m\geq 1$,  it holds that
\begin{align}\label{step_14}
	\mathbb{E}\Big(\left(N_{[1,m] ,\kappa}^y \right)^{1+\alpha}\Big)& =\mathbb{E}\bigg( \sum_{k=1}^m
	\sum_{x\in\mathcal{N}(k)}
	e^{-V(x)}1_{\left\{V(x)<-y, \min_{j\leq k-1} V(x_j)\geq -y \right\}}1_{\left\{ x\in \mathcal{A}_{\kappa}^y\right\}}\left(N_{[1,m] ,\kappa}^{y}\right)^{\alpha}\bigg) \nonumber \\
	& = \mathbb{E}_{\Q}\bigg( 1_{\left\{\tau_{-y}^- \leq m \right\}}1_{\left\{ w_{\tau_{-y}^-}\in \mathcal{A}_{\kappa}^y\right\}}\left(N_{[1,m] ,\kappa}^y\right)^\alpha\bigg).
\end{align}
 Using the trivial inequality
$$
(x+y)^\alpha \leq x^\alpha + y^\alpha, \quad x,y\geq0, \alpha\in (0, 1],
$$
we get that
\begin{align}\label{I+II}
	&\mathbb{E}\left(\left(N_{[1,m] ,\kappa}^y \right)^{1+\alpha}\right) \leq  \mathbb{E}_{\Q}\bigg( 1_{\left\{\tau_{-y}^- \leq m \right\}}1_{\big\{ w_{\tau_{-y}^-}\in \mathcal{A}_{\kappa}^y\big\}}\bigg(N_{[1,m] ,\kappa}^y-e^{-V\big(w_{\tau_{-y}^-}\big)} \bigg)^\alpha\bigg)\nonumber \\
	&\qquad\qquad\qquad\qquad \quad + \mathbb{E}_{\Q}\bigg( 1_{\big\{\tau_{-y}^- \leq m \big\}}1_{\big\{ w_{\tau_{-y}^-}\in \mathcal{A}_{\kappa}^y\big\}}e^{-\alpha V\big(w_{\tau_{-y}^-}\big)}\bigg)=:I+II.
\end{align}
By the spine decomposition, we have
\begin{align*}
	    N_{[1,m] ,\kappa}^y =e^{-V\big(w_{\tau_{-y}^- \land m}\big)}+ \sum_{k=1}^{\tau_{-y}^- \land m} \sum_{u \in\Omega(w_k)} e^{-V(u)}N_{[1, m-k],\kappa}^{y+V(u)}.
\end{align*}
By the branching property, \eqref{step_12} and  using the trivial inequality $\mathbb{E}(|X|^\alpha)\leq \mathbb{E}(|X|)^\alpha$  (since $\alpha \in (0,1]$),  we have
\begin{align}\label{step_15}
     I^{1/\alpha}
     & \leq \mathbb{E}_{\Q}\bigg( 1_{\left\{\tau_{-y} ^- \leq m \right\}}1_{\big\{ w_{\tau_{-y}^-}\in \mathcal{A}_{\kappa}^y\big\}}\Big(N_{[1,m] ,\kappa}^y-e^{-V(w_{\tau_{-y}^-})} \Big)\bigg)  \nonumber\\
		 &=\mathbb{E}_{\Q}\bigg( 1_{\left\{\tau_{-y}^- \leq m \right\}}1_{\big\{ w_{\tau_{-y}^-}\in \mathcal{A}_{\kappa}^y\big\}}\sum_{\ell =1}^{\tau_{-y}^-} \sum_{u \in\Omega(w_\ell)}e^{-V(u)}\mathbb{E}_{\Q}\left( N_{[1, m-\ell], \kappa}^{z}\right)_{z=V(u)+y} \bigg)\nonumber\\
		& \leq \mathbb{E}_{\Q}\bigg( 1_{\big\{ w_{\tau_{-y}^-}\in \mathcal{A}_{\kappa}^y\big\}}\sum_{\ell =1}^{\tau_{-y}^-} \sum_{u \in\Omega(w_\ell )}e^{-V(u)} \bigg)\nonumber\\
		& \leq \mathbb{E}_{\Q}\bigg( 1_{\big\{ w_{\tau_{-y}^-}\in \mathcal{A}_{\kappa}^y\big\}}\sum_{\ell =1}^{\tau_{-y}^-} e^{-V(w_{\ell -1})}\sum_{u \in\Omega(w_\ell)}\left(1+\left(V(u)-V(w_{\ell-1}) \right)_+\right)e^{-\left(V(u)-V(w_{\ell-1})\right)} \bigg)\nonumber\\
		& \leq  \kappa \mathbb{E}_{\Q}\bigg( 1_{\big\{ w_{\tau_{-y}^-}\in \mathcal{A}_{\kappa}^y\big\}}\sum_{\ell=1}^{\tau_{-y}^-} e^{-\left(V(w_{\ell-1})-y\right)/2} \bigg)\leq \kappa e^{y} \mathbb{E}_{\Q}\bigg(\sum_{\ell =1}^{\tau_{-y}^-} e^{-\left(V(w_{\ell -1})+y\right)/2} \bigg).
\end{align}
Using Lemma \ref{lemma1} (iv) with $\lambda=1/2$,
we get that
\begin{align}\label{step_16}
		&\mathbb{E}_{\Q}\bigg(\sum_{\ell=1}^{\tau_{-y}^-} e^{-\left(V(w_{\ell -1})+y\right)/2} \bigg) \leq 1+ \mathbb{E}_{\Q}\bigg(\sum_{\ell=1}^{\tau_{-y}^- -1} e^{-\left(V(w_{\ell})+y\right)/2} \bigg)\nonumber\\
				& \leq  1+ \sum_{\ell =1}^\infty \mathbf{E}_y\Big(e^{-S_\ell/2}1_{\left\{\min_{j\leq \ell} S_j \geq 0\right\}}\Big)
		\leq 1+ C_1(1/2)\lesssim 1.
\end{align}
Then we have
\begin{align}\label{I}
I&\lesssim\kappa^\alpha e^{\alpha y}.
\end{align}
Finally, using Lemma \ref{lemma1} (iv) with $\lambda=1$ and {\bf(A5)}, we get
\begin{align}\label{step_18}
II&
\leq \mathbb{E}_{\Q}\bigg( e^{-\alpha V\big(w_{\tau_{-y}^-}\big)}\bigg)=
\mathbf{E}\bigg(e^{-\alpha S_{\tau_{-y}^-}}\bigg)\nonumber\\
	 & = e^{\alpha y} \sum_{\ell=1}^\infty
		 \mathbf{E}_y
	 \Big(e^{-\alpha S_\ell} 1_{\left\{\min_{k\leq \ell -1}S_k \geq 0\right\}} 1_{\left\{S_\ell <0\right\}}\Big)\nonumber\\
	 & \leq e^{\alpha y} \sum_{\ell=1}^\infty
		 \mathbf{E}_y
	 \Big(e^{-\alpha S_{\ell-1}} 1_{\left\{\min_{k\leq \ell -1}S_k \geq 0\right\}} \Big)
		 \mathbf{E}_y
	 \Big(e^{-\alpha (S_\ell -S_{\ell-1})}\Big)\nonumber\\
	 &\leq C_1(1)e^{\alpha y}
		 \mathbf{E}
	 \big(e^{-\alpha S_1}\big) = C_1(1)e^{\alpha y} \mathbb{E}\bigg(\sum_{u\in\mathcal{N}(1)} e^{-(1+\alpha)V(u)}\bigg)\lesssim e^{\alpha y}.
\end{align}
Hence, combining \eqref{I+II}, \eqref{I} and \eqref{step_18}, we get that
$$\mathbb{E}\left(\left(N_{[1,m] ,\kappa}^y \right)^{1+\alpha}\right)\lesssim \kappa^\alpha e^{\alpha y}.$$
This completes the proof of the Lemma.
\hfill$\Box$

\medskip

For a sequence $(\beta_n)_{n\geq 1}$ of positive numbers, define
\begin{align*}
	& \widehat{N}_{[\ell ,m]}^{y,n}:= \sum_{q=\ell }^m \sum_{x\in\mathcal{N}(q)}	\left(V(x)+\beta_n +y\right)e^{-V(x)} 1_{\left\{-y -\beta_n/2 \leq V(x)<-y,\, \min_{j\leq q-1} V(x_j)\geq -y \right\}},\nonumber\\
	&
		\widehat{N}_{[\ell,m],\kappa}^{y,n}:= \sum_{q=\ell}^m \sum_{x\in\mathcal{N}(q)} \left(V(x)+\beta_n +y\right)e^{-V(x)} 1_{\left\{-y -\beta_n/2 \leq V(x)<-y,\, \min_{j\leq q-1} V(x_j)\geq -y \right\}}1_{\left\{ x\in \mathcal{A}_{\kappa}^y\right\}},
\end{align*}
and let   $\widehat{N}_{[\ell,\infty)}^{y,n}:=\lim_{m\to\infty} \widehat{N}_{[\ell,m]}^{y,n}$ and $\widehat{N}_{[\ell,\infty),\kappa}^{y,n}:=\lim_{m\to\infty}\widehat{N}_{[\ell,m],\kappa}^{y,n}$.
    Then we have the following result similar to that in Lemma \ref{lemma7} and Lemma \ref{lemma8} for $\widehat{N}_{[1,\infty)}^{y,n}$ and $\widehat{N}_{[1,\infty),\kappa}^{y,n}$:

\begin{lemma}\label{lemma10}
	(i)
	Let $g$ be the function in Lemma \ref{lemma7} (ii). Then
	\[
	\mathbb{E}\Big(\widehat{N}_{[1,\infty)}^{y,n} - \widehat{N}_{[1, \infty), \kappa}^{y,n}\Big)\lesssim \beta_n \frac{g(\kappa)}{\log \kappa}, \quad y >0, \kappa > 1, n\geq 1.
	\]
	(ii)
	Let $\alpha$ be the constant in {\bf(A5)}. Then
	\[
	 \mathbb{E}\Big(\left(\widehat{N}_{[1,\infty) ,\kappa}^{y,n}\right)^{1+\alpha}\Big)\lesssim \beta_n^{1+\alpha} \kappa^\alpha e^{\alpha y}, \quad y\geq 0, \kappa \geq 1, n\geq 1.
	\]
\end{lemma}
\textbf{Proof: }(i) By direct calculation and Lemma \ref{lemma7} (ii), we have
\begin{align*}
&\quad	\mathbb{E}\left(\widehat{N}_{[1,\infty)}^{y,n} - \widehat{N}_{[1, \infty), \kappa}^{y,n}\right)\\
  &\leq  \mathbb{E}\bigg(\sum_{q=1}^\infty \sum_{x\in\mathcal{N}(q)}	\left(V(x)+\beta_n +y\right)_+	e^{-V(x)} 1_{\left\{V(x)<-y, \,\min_{j\leq q-1} V(x_j)\geq -y \right\}}1_{\left\{ x\notin \mathcal{A}_{\kappa}^y\right\}}\bigg)\\
		& \leq \beta_n \mathbb{E}\bigg(\sum_{q=1}^\infty \sum_{x\in\mathcal{N}(q)}	e^{-V(x)} 1_{\left\{V(x)<-y,\, \min_{j\leq q-1} V(x_j)\geq -y \right\}}1_{\left\{ x\notin \mathcal{A}_{\kappa}^y\right\}}\bigg)\\
	& = \beta_n\mathbb{E}\left(N_{[1,\infty)}^y - N_{[1, \infty), \kappa}^y\right)\lesssim \beta_n\frac{g(\kappa)}{\log \kappa}.
\end{align*}

(ii)  Combining Lemma \ref{lemma8} and  the inequality
$0\leq \widehat{N}_{[1,\infty) ,\kappa}^{y,n} \leq \beta_n N_{[1,\infty) ,\kappa}^y$, we immediately get the desired conclusion.
\hfill$\Box$

\section{Proof of Theorem \ref{thm1}}\label{main}

For $n\geq 1$, set
\begin{equation}\label{e:rs1}
\gamma_n:= \frac{1}{2}\log n + \beta_n,
\end{equation}
here $\left(\beta_n\right)_{n\geq1}$ is a sequence of positive numbers with
\begin{align}\label{step_4}
	\beta_n -\frac{1}{2}\log n \to +\infty,\quad \frac{\beta_n}{n^{1/16}} \to 0,\quad \mbox{as}\  n\to\infty.
\end{align}

We first give the main idea of the proof of Theorem \ref{thm1}.
Recall that $c^*$ and $\alpha^*$ are the constants in \eqref{Limit-of-R(y)/y} and \eqref{Second-Exp-of-R(y)}.
Let $\ell (y):= c^* y +\alpha^*$.
Using \eqref{Seneta-Heyde-scaling}, we get
\begin{align}\label{Decomposition-2a}
	c^* D_\infty & = \lim_{m\to\infty}  \sum_{x\in \mathcal{N}(m)}
	\ell \left(V(x) - \gamma_n \right)e^{-V(x)}.
\end{align}
For $x\in\mathbb{T}$ and $m\ge 0$, let $\mathbb{T}_x$ be the subtree of $\mathbb{T}$ with root at $x$ and $\mathcal{N}(x, m)$ be the collections of particles in the $m$-th generation of $\mathbb{T}_x$. Define
\begin{equation}\label{e:D_n(x)}
	D_m(x)
	:=\sum_{u\in \mathcal{N}(x, m)}(V(u)-V(x))e^{-(V(u)-V(x))},
\end{equation}
and let $D_\infty(x)$ be the limit of $D_m(x)$.

For $n \ge 1$ and $m\ge \lceil an \rceil$, we define a quantity
$D_m^{\lceil an\rceil, \gamma_n}$ which roughly takes care of the contributions to $D_m$ by the
paths that stay above the level $\gamma_n$ between generations $\lceil an \rceil$ and $m$. We  show  in Lemma \ref{lemma2} (i) that $D_m^{\lceil an\rceil, \gamma_n}$ converges to a limit $D_\infty^{\lceil an\rceil, \gamma_n}$ as $m\to\infty$.

For any $n\in \mathbb{N}$ and $a\geq 1$, we define
\begin{align*}
	&\mathcal{L}^{\lceil an\rceil, \gamma_n}:= \Big\{x\in \mathbb{T}: 	|x|> \lceil an \rceil,
	V(x) < \gamma_n \ \mbox{and}\ 	 \min_{j\in [an,|x|-1]\cap\mathbb{Z}} V(x_j)\geq \gamma_n
	\Big\}.
\end{align*}
We will use $\mathcal{F}_{\mathcal{L}^{\lceil an\rceil, \gamma_n}}$ to denote the smallest $\sigma$-field containing all the information about the particles $x\in \mathcal{L}^{\lceil an\rceil, \gamma_n}$ (about $x_j$ and $V(x_j)$ for all $x\in \mathcal{L}^{\lceil an\rceil, \gamma_n}$ and $j\le |x|$).
 It is obvious that $\mathcal{F}_{\lceil an\rceil}\subset \mathcal{F}_{\mathcal{L}^{\lceil an\rceil, \gamma_n}}$.

To consider contributions to $D_\infty$
by particles in  $\mathcal{L}^{\lceil an\rceil, \gamma_n}$,
we separate $\mathcal{L}^{\lceil an\rceil, \gamma_n}$ into two sets $ \mathcal{L}_{good}^{\lceil an\rceil, \gamma_n}$ and $ \mathcal{L}_{bad}^{\lceil an\rceil, \gamma_n}$:
 \begin{align*}
  	&\mathcal{L}_{good}^{\lceil an\rceil, \gamma_n} := \Big\{ x \in\mathcal{L}^{\lceil an\rceil, \gamma_n}: \min_{n\leq j\leq \lceil an\rceil} V(x_j)\geq \gamma_n , \ V(x)\geq  \gamma_n -\frac{\beta_n}{2}\Big\}, \\
	&\mathcal{L}_{bad}^{\lceil an\rceil, \gamma_n} := \Big\{ x \in\mathcal{L}^{\lceil an\rceil, \gamma_n}: \min_{n\leq j\leq \lceil an\rceil} V(x_j) < \gamma_n \Big\}\\
	&\quad \quad\quad\quad\bigcup  \Big\{ x \in\mathcal{L}^{\lceil an\rceil, \gamma_n}:
	\min_{n\leq j\leq |x|-1} V(x_j)\geq  \gamma_n , \ V(x) <\gamma_n -\frac{\beta_n}{2}\Big\}.
\end{align*}
By \eqref{Decomposition-2a}, the branching property  and  Lemma \ref{lemma2} (i) below, we get
\begin{align}\label{Decomposition-2}
	c^* D_\infty & 	=  D_\infty^{\lceil an\rceil, \gamma_n} + \sum_{x\in \mathcal{L}^{\lceil an\rceil, \gamma_n}} c^* e^{-V(x)} D_\infty(x)\nonumber \\
	& =  D_\infty^{\lceil an\rceil, \gamma_n} + \sum_{x\in \mathcal{L}_{good}^{\lceil an\rceil, \gamma_n}} c^* e^{-V(x)} D_\infty(x)+
	\sum_{x\in \mathcal{L}_{bad}^{\lceil an\rceil, \gamma_n}} c^* e^{-V(x)} D_\infty(x)\nonumber\\
	& =: D_\infty^{\lceil an\rceil, \gamma_n} + F_{good}^{\lceil an \rceil, \gamma_n} +
	F_{bad}^{\lceil an \rceil, \gamma_n}.
\end{align}
Using this, we get that
\begin{align}\label{Decomposition-3}
	D_\infty - D_{\lceil an \rceil} + \frac{\log n}{2}W_{\lceil an \rceil}& = \frac{1}{c^*}\left(D_\infty^{\lceil an\rceil, \gamma_n} -\left(c^* D_{\lceil an \rceil} -\left(c^* \gamma_n -\alpha^*\right)W_{\lceil an \rceil}\right) \right)\nonumber \\
	&  + \frac{1}{c^*}\left(F_{good}^{\lceil an \rceil, \gamma_n}-\left(c^*\beta_n-\alpha^*\right) W_{\lceil a n\rceil }\right) +
	 \frac{1}{c^*}F_{bad}^{\lceil an \rceil, \gamma_n}.
\end{align}

In Section \ref{subsec4.1}, we prove that $D_\infty^{\lceil an\rceil, \gamma_n} \approx c^* D_{\lceil an \rceil} -\left(c^* \gamma_n -\alpha^*\right)W_{\lceil an \rceil}$, see Proposition \ref{prop3} below. In Section \ref{subsec4.2}, we analyze the weighted number of particles $\widehat{N}_{good}^{\lceil an \rceil, \gamma_n }$ defined in \eqref{Widehat-N} below, and prove that $\widehat{N}_{good}^{\lceil an \rceil, \gamma_n }\approx \left(c^* \beta_n -\alpha^*\right)W_{\lceil an \rceil} $, see  Corollary \ref{cor2} below.
In Section \ref{Proof-Prop1}, we prove that
$F_{bad}^{\lceil an \rceil, \gamma_n}$ is negligible, see Proposition \ref{prop1}.

Then, by \eqref{Decomposition-3},  we see that
 \begin{align}
     c^* \Big(D_\infty -D_{\lceil a n\rceil}   + \frac{\log n}{2} W_{\lceil an \rceil}\Big) \approx F_{good}^{\lceil an \rceil, \gamma_n} -\widehat{N}_{good}^{\lceil an \rceil, \gamma_n}.
\end{align}
In Subsection \ref{subsec4.3},  we  prove the convergence of $\sqrt{n}(F_{good}^{\lceil an \rceil, \gamma_n} -\widehat{N}_{good}^{\lceil an \rceil, \gamma_n})$ to
$c^*X_{a^{-1/2} D_\infty}$ in distribution,
see Proposition \ref{prop4} below. Using these, we can easily get the conclusion of Theorem \ref{thm1}.

\subsection{Modifications of the martingales with level $\gamma_n$}\label{subsec4.1}

For $a\geq 1$, $n\in\mathbb{N}$ and $m\geq \lceil an \rceil$, define
\begin{align*}
	\widetilde{D}_m^{\lceil an\rceil, \gamma_n}&:=
	\sum_{x\in\mathcal{N}(m)}
	R\left(V(x)-\gamma_n\right)e^{-V(x)}1_{\left\{ \min_{j\in [an,m]\cap \mathbb{Z}} V(x_j)\geq  \gamma_n \right\}}, \\
	D_m^{\lceil an\rceil, \gamma_n}&:=
	\sum_{x\in\mathcal{N}(m)}
	\ell \left(V(x)-\gamma_n\right)e^{-V(x)}1_{\left\{ \min_{j\in [an,m]\cap \mathbb{Z}} V(x_j)\geq  \gamma_n \right\}}.
\end{align*}

\begin{lemma}\label{lemma2}
	Let $a\ge 1$. Then
	
	(i) $D_\infty^{\lceil an\rceil, \gamma_n}:=\lim_{m\to\infty} D_m^{\lceil an\rceil, \gamma_n}$ exists $\P$-almost surely.
	
	(ii) Moreover, under $\P$,
	$$ \lim_{n\to\infty} \sqrt{n} \left|D_\infty^{\lceil an\rceil, \gamma_n} - D_{\lceil an \rceil}^{\lceil an\rceil, \gamma_n} \right| = 0,\quad \mbox{in probability}.$$
\end{lemma}
\textbf{Proof:} (i) For $m\geq \lceil an\rceil$, by the branching property,
\begin{align*}
	\mathbb{E}\left(\widetilde{D}_{m+1}^{\lceil an\rceil, \gamma_n} \big| \mathcal{F}_m\right) &= 	\sum_{x\in\mathcal{N}(m)} 1_{\left\{ \min_{j\in [an,m]\cap \mathbb{Z}} V(x_j)\geq \gamma_n \right\}}e^{-V(x)}\\
	&\quad\times  \mathbb{E}\bigg(	\sum_{x'\in\mathcal{N}(1)} 	R\left(V(x') + z-\gamma_n\right)e^{-V(x')} 1_{\left\{ V(x')+z\geq \gamma_n \right\}}	\bigg)\bigg|_{z= V(x)}\\
	& =
	\sum_{x\in\mathcal{N}(m)}
	1_{\left\{ \min_{j\in [an,m]\cap \mathbb{Z}} V(x_j)\geq \gamma_n \right\}}e^{-V(x)}  R\left(  V(x)-\gamma_n\right)= \widetilde{D}_m^{\lceil an\rceil, \gamma_n}.
\end{align*}
Thus $\left(\widetilde{D}_{m}^{\lceil an\rceil, \gamma_n}\right)_{m\geq \lceil an\rceil}$ is a non-negative martingale and hence $D_\infty^{\lceil an\rceil, \gamma_n}:=\lim_{m\to\infty} \widetilde{D}_{m}^{\lceil an\rceil, \gamma_n}$ exists $\P$-almost surely.
It follows from \eqref{Second-Exp-of-R(y)} that $\sup_{y\geq 0} \left|R(y)-\ell(y)\right| <\infty$. Therefore,
\begin{align*}
	\left|D_m^{\lceil an\rceil, \gamma_n} - \widetilde{D}_m^{\lceil an\rceil, \gamma_n}\right| \leq \sup_{y\geq 0} \left|R(y)-\ell(y)\right| 	\sum_{x\in\mathcal{N}(m)} e^{-V(x)} =\sup_{y\geq 0} \left|R(y)-\ell(y)\right|  W_m,  \quad \P\text{--a.s.}
\end{align*}
Since $\lim_{m\to\infty} W_m =0$, we get the desired conclusion.

(ii) For $\kappa \geq 1$ and $m\geq \lceil an \rceil+1$, define
\begin{align*}
	B_{m,\kappa}^{\lceil an \rceil, \gamma_n}&:= \bigg\{x\in\mathcal{N}(m):\  \forall \lceil an \rceil+1 \leq j \leq m,\\
	& \qquad \sum_{\Omega(x_j)} \left(1+\left(V(u)-V(x_{j-1})\right)_+\right)e^{-\left(V(u)-V(x_{j-1})\right)} \leq \kappa e^{\left(V(x_{j-1})-\gamma_n \right)/2}
	\bigg\},\\	\widetilde{D}_{m,\kappa}^{\lceil an\rceil, \gamma_n}&:=
	\sum_{x\in\mathcal{N}(m)}
	R\left(V(x)-\gamma_n\right)e^{-V(x)}1_{\left\{ \min_{j\in [an,m]\cap \mathbb{Z}} V(x_j)\geq \gamma_n \right\}} 1_{\left\{x\in B_{m,\kappa}^{\lceil an \rceil, \gamma_n} \right\}}.
\end{align*}
Then by the branching property and Lemma \ref{lemma4},
\begin{align}\label{step_2}
	& \mathbb{E}\left( \widetilde{D}_m^{\lceil an\rceil, \gamma_n} -  \widetilde{D}_{m,\kappa}^{\lceil an\rceil, \gamma_n} \Big| \mathcal{F}_{\lceil an \rceil}\right)\nonumber\\
	 & =  \sum_{u\in \mathcal{N}(\lceil an \rceil)}
	 e^{-V(u)} 1_{\{V(u)\geq \gamma_n \}}
	 \mathbb{E} \left(D_{m-  \lceil an \rceil }^{-y} - D_{m- \lceil an \rceil ,\kappa}^{-y} \right)\bigg|_{y= V(u)-\gamma_n} \lesssim W_{\lceil an\rceil} h(\kappa),
\end{align}
where $h$ is the function in Lemma \ref{lemma4}. By Lemma \ref{lemma5}, we have
\begin{align}\label{step_3}
	\textup{Var}\left( \widetilde{D}_{m,\kappa}^{\lceil an\rceil, \gamma_n} \Big| \mathcal{F}_{\lceil an \rceil}\right)& \leq 	\sum_{u\in \mathcal{N}(\lceil an \rceil)}
	e^{-2V(u)} \mathbb{E}\Big(\left( D_{m- \lceil an \rceil ,\kappa}^{-y}\right)^2 \Big)\Big|_{y= V(u)-\gamma_n}\nonumber \\ &\lesssim \kappa 	 \sum_{u\in \mathcal{N}(\lceil an \rceil)}
	 e^{-2V(u)} e^{V(u)-\gamma_n}= \kappa e^{-\gamma_n}W_{\lceil an \rceil}.
\end{align}
Combining \eqref{step_2}--\eqref{step_3} with $\mathbb{E}\left(\widetilde{D}_m^{\lceil an\rceil, \gamma_n}\Big| \mathcal{F}_{\lceil an \rceil} \right) = \widetilde{D}_{\lceil an \rceil}^{\lceil an\rceil, \gamma_n}$, we get that, for any $\varepsilon >0$,
\begin{align*}
	& \P\Big(\left|\widetilde{D}_m^{\lceil an\rceil, \gamma_n} -\widetilde{D}_{\lceil an\rceil}^{\lceil an\rceil, \gamma_n} \right|\geq \frac{3\varepsilon}{\sqrt{\lceil an \rceil}}\Big| \mathcal{F}_{\lceil an \rceil}\Big)\\ &
		\leq \P\Big(\widetilde{D}_m^{\lceil an\rceil, \gamma_n} -\widetilde{D}_{m,\kappa}^{\lceil an\rceil, \gamma_n} \geq \frac{\varepsilon}{\sqrt{\lceil an \rceil}}\Big| \mathcal{F}_{\lceil an \rceil}\Big)\\
	&\quad + \P\Big(\left|\widetilde{D}_{m,\kappa}^{\lceil an\rceil, \gamma_n} -\mathbb{E}\left(\widetilde{D}_{m,\kappa}^{\lceil an\rceil, \gamma_n}\Big|\mathcal{F}_{\lceil an \rceil} \right) \right|\geq \frac{\varepsilon}{\sqrt{\lceil an \rceil}}\Big| \mathcal{F}_{\lceil an \rceil}\Big)\\
	& \quad + 1_{\left\{\mathbb{E}\left(\widetilde{D}_m^{\lceil an\rceil, \gamma_n} -\widetilde{D}_{m,\kappa}^{\lceil an\rceil, \gamma_n}  \Big|\mathcal{F}_{\lceil an \rceil} \right) \geq \varepsilon /\sqrt{\lceil an \rceil }\right\}}\\
	& \leq 2\frac{\lceil an \rceil }{\varepsilon} \mathbb{E}\left(\widetilde{D}_m^{\lceil an\rceil, \gamma_n} -\widetilde{D}_{m,\kappa}^{\lceil an\rceil, \gamma_n}  \Big|\mathcal{F}_{\lceil an \rceil}\right) + \frac{\lceil an \rceil^2}{\varepsilon^2}\textup{Var}\left( \widetilde{D}_{m,\kappa}^{\lceil an\rceil, \gamma_n} \Big| \mathcal{F}_{\lceil an \rceil}\right)\\
	& \lesssim \lceil an \rceil W_{\lceil an \rceil}\left(\frac{2 h(\kappa)}{\varepsilon}+ \frac{\lceil an \rceil}{\varepsilon^2} \kappa e^{-2^{-1}\log n -\beta_n}\right).
\end{align*}
Letting $m\to\infty$, we get that for all $\kappa \geq 1$ and $n\geq 1$,
\begin{align}\label{step_41}
	\P\Big(\left|{D}_\infty^{\lceil an\rceil, \gamma_n} -\widetilde{D}_{\lceil an\rceil}^{\lceil an\rceil, \gamma_n} \right|\geq \frac{3\varepsilon}{\sqrt{\lceil an \rceil}}\Big| \mathcal{F}_{\lceil an \rceil}\Big)\lesssim \lceil an \rceil W_{\lceil an \rceil}\left( h(\kappa) +  \kappa e^{-\left(\beta_n -2^{-1}\log n\right)}\right).
\end{align}
Using \eqref{Seneta-Heyde-scaling},  \eqref{step_4} and \eqref{step_41},
first letting $n\to\infty$ and then $\kappa \to\infty$, we get
\begin{align}\label{step_7}
	\lim_{n\to\infty} \sqrt{n} \left|D_\infty^{\lceil an\rceil, \gamma_n} - \widetilde{D}_{\lceil an \rceil}^{\lceil an\rceil, \gamma_n} \right| = 0,\quad \mbox{in probability}.
\end{align}
On the other hand, by \eqref{Second-Exp-of-R(y)} we know that there exists $C>0$ such that $ \left|R(y)-\ell(y)\right|\leq C$ for all $y\geq 0$ and that, for any $\eta >0$, there exists $K=K(\eta)>0$ such that $\left|R(y)-\ell(y)\right|<\eta$ for  $y > K$.
Since $\min_{x \in\mathcal{N}(m)} V(x) \to +\infty$ as $m\to\infty$,
for any $\delta >0$, there exists $L>0$ such that
\begin{align}\label{Global-Minimum}
	 \mathbb{P}\big(  \min_{x\in \mathbb{T}} V(x) < -L \big) \leq \delta.
\end{align}
Therefore,
\begin{align}\label{step_42}
	\left| \widetilde{D}_{\lceil an \rceil}^{\lceil an\rceil, \gamma_n} -D_{\lceil an \rceil}^{\lceil an\rceil, \gamma_n}\right|\leq \eta W_{\lceil an \rceil}+ C
	\sum_{x\in \mathcal{N}(\lceil an \rceil)}
	e^{-V(x)}1_{\{\gamma_n \leq V(x)< \gamma_n +K \}}.
\end{align}
Combining \eqref{Global-Minimum}, \eqref{step_42} and Markov's inequality, we get that, for any $\varepsilon>0$,
\begin{align}
	&\mathbb{P}\Big(\left|\widetilde{D}_{\lceil an \rceil}^{\lceil an\rceil, \gamma_n} -D_{\lceil an \rceil}^{\lceil an\rceil, \gamma_n}\right|>\frac{\varepsilon}{\sqrt{n}} \Big)\leq \delta + \P\left( \eta \sqrt{n}W_{\lceil an \rceil} > \varepsilon/2\right) \nonumber\\
	&\qquad\quad\quad\quad +\P\Big(C \sqrt{n}	\sum_{x\in \mathcal{N}(\lceil an \rceil)}	e^{-V(x)}1_{\{\gamma_n \leq V(x)< \gamma_n +K \}}1_{\left\{\min_{j\leq \lceil an\rceil} V(x_j)\geq -L \right\}} > \varepsilon/2 \Big)\nonumber\\
	& \leq \delta + \P\left( \eta \sqrt{n}W_{\lceil an \rceil} > \varepsilon/2\right) \nonumber\\
	&\qquad\quad\quad\quad+ \frac{2C \sqrt{n}}{\varepsilon}\mathbb{E}\Big(	\sum_{x\in \mathcal{N}(\lceil an \rceil)}	e^{-V(x)}1_{\{\gamma_n \leq V(x)< \gamma_n +K \}}1_{\left\{\min_{j\leq \lceil an\rceil} V(x_j)\geq -L \right\}}\Big) \nonumber\\
	& = \delta + \P\left( \eta \sqrt{n}W_{\lceil an \rceil} > \varepsilon/2\right)+ \frac{2C \sqrt{n}}{\varepsilon}
	\mathbf{P}
	\big(\gamma_n \leq S_{\lceil an \rceil}< \gamma_n +K ,\min_{j\leq \lceil an\rceil} S_j\geq -L \big).\label{e:rs2}
\end{align}
Since $\gamma_n = o(n^{1/16})$, by Lemma \ref{lemma1}(ii) we have
\begin{align*}
	&\limsup_{n\to\infty} \frac{2C \sqrt{n}}{\varepsilon}
	\mathbf{P}
	\big(\gamma_n \leq S_{\lceil an \rceil}< \gamma_n +K ,\min_{j\leq \lceil an\rceil} S_j>-L \big) \\ & \lesssim \lim_{n\to\infty} \frac{2C \sqrt{n}}{\varepsilon} \frac{(L+1)(K+1)(\gamma_n + L+ K +1)}{\sqrt{\lceil an \rceil^3}}=0.
\end{align*}
Thus, using \eqref{Seneta-Heyde-scaling} and letting $n\to\infty$ in \eqref{e:rs2}, we get
\begin{align}\label{step_8}
	\limsup_{n\to\infty} \mathbb{P}\Big(\left|\widetilde{D}_{\lceil an \rceil}^{\lceil an\rceil, \gamma_n} -D_{\lceil an \rceil}^{\lceil an\rceil, \gamma_n}\right|>\frac{\varepsilon}{\sqrt{n}} \Big) \leq \delta +
	\P\Big(\eta \sqrt{\frac{2}{\pi\sigma^2}}a^{-1/2} D_\infty> \varepsilon/2 \Big).
\end{align}
Letting $\delta ,\eta \to 0$ in \eqref{step_8}, and then combining the resulting fact with  \eqref{step_7}, we get the desired conclusion.
\hfill$\Box$

\bigskip

For $m\geq n$, define
\begin{equation}\label{e:rs3}
\widetilde{W}_{m}^{n, \gamma_n}:= \sum_{x\in\mathcal{N}(m)}
e^{-V(x)}1_{\left\{ \min_{j\in [n,m]\cap \mathbb{Z}} V(x_j)\geq  \gamma_n \right\}}.
\end{equation}

\begin{lemma}\label{lemma6}
	(i) For all $a\geq 1, b\in\R $, as $n\to\infty$,  $\left(\gamma_n +b\right) \sqrt{n} \big(W_{\lceil an \rceil} - \widetilde{W}_{\lceil an \rceil}^{\lceil an \rceil, \gamma_n}\big) \to 0$ in probability.
	
	(ii) For all $a\geq 1, b\in\R$, as $n\to\infty$,  $\left(\gamma_n +b\right) \sqrt{n} \big(W_{\lceil an \rceil} - \widetilde{W}_{\lceil an \rceil}^{n, \gamma_n}\big) \to 0$ in probability.
	Moreover, if there exists a sequence of random variables $\left\{J_n \right\}$ such that for all $n\geq 1$,
	\[
	\mathbb{E}\left( \left|J_n\right| \Big| \mathcal{F}_{\lceil an \rceil}\right) \lesssim \left(\gamma_n +b\right) \sqrt{n} \left(W_{\lceil an \rceil} - \widetilde{W}_{\lceil an \rceil}^{n, \gamma_n}\right),
	\]
	 then $J_n \to 0$ in probability.
\end{lemma}
\textbf{Proof: } (i) Without loss of generality, we assume $b\geq 0$. Fix $\delta>0$ and  let $L$ be the constant in  \eqref{Global-Minimum}.
Using Markov's inequality and Lemma \ref{lemma1} (ii), we know that for any $\varepsilon>0$,
\begin{align*}
	& \P \Big(W_{\lceil an \rceil} - \widetilde{W}_{\lceil an \rceil}^{\lceil an \rceil, \gamma_n}\geq \frac{\varepsilon}{\left(\gamma_n+b\right) \sqrt{n}} \Big)
	\\ & \leq \delta + \P \Big(W_{\lceil an \rceil} - \widetilde{W}_{\lceil an \rceil}^{\lceil an \rceil, \gamma_n}\geq \frac{\varepsilon}{\left(\gamma_n +b\right) \sqrt{n}},\, \min_{j\leq \lceil an\rceil }	\min_{x\in\mathcal{N}(\lceil an \rceil)}	V(x_j) \geq -L\Big)\nonumber
	\\ & \leq \delta + \frac{\left(\gamma_n +b\right) \sqrt{n}}{\varepsilon}\mathbb{E}\bigg(
	\sum_{x\in\mathcal{N}(\lceil an \rceil )}
	e^{-V(x)}1_{\left\{ V(x)<  \gamma_n \right\}} 1_{\{\min_{j\leq \lceil an\rceil } V(x_j)\geq -L \}}\bigg)\nonumber\\
	& = \delta + \frac{\left(\gamma_n +b\right) \sqrt{n}}{\varepsilon}
		\mathbf{P}_L
	\Big(S_{\lceil an \rceil} < \gamma_n+L , \min_{j\leq \lceil an\rceil } S_j \geq 0\Big)\lesssim \delta + \frac{\left(\gamma_n +b\right) \sqrt{n}}{\varepsilon}\cdot\frac{(\gamma_n +1)^2 (L+1)}{\sqrt{\lceil an \rceil^3}}.
\end{align*}
Since $\gamma_n = o(n^{1/16})$,  by taking $n\to\infty$ first and then $\delta \to 0$ in the display above, we get the conclusion of (i).

(ii)
Similarly we assume that $b\geq 0$. It suffices to prove the second result since the first one holds by taking $J_n= \left(\gamma_n +b\right) \sqrt{n} \left(W_{\lceil an \rceil} - \widetilde{W}_{\lceil an \rceil}^{n, \gamma_n}\right).$  Using the same argument as in (i) gives
\begin{align}\label{step_5}
	& \P\left(\left|J_n\right|> \varepsilon\right) \leq \delta + \frac{1}{\varepsilon}\mathbb{E}\left(\mathbb{E}\left(\left|J_n\right| \Big| \mathcal{F}_{\lceil an \rceil}\right) 1_{\left\{\min_{j\leq \lceil an\rceil }	\min_{x\in\mathcal{N}(\lceil an \rceil)}
	V(x_j) \geq -L \right\}} \right) \nonumber\\
	& \leq \delta + \frac{\left(\gamma_n +b\right) \sqrt{n}}{\varepsilon}\mathbb{E}\left(\left(W_{\lceil an \rceil} - \widetilde{W}_{\lceil an \rceil}^{n, \gamma_n}\right) 1_{\left\{\min_{j\leq \lceil an\rceil }
	\min_{x\in\mathcal{N}(\lceil an \rceil)}	V(x_j) \geq -L \right\}} \right)\nonumber
	\\ &
	\leq
	 \delta + \frac{\left(\gamma_n+b\right) \sqrt{n}}{\varepsilon}\mathbb{E}\bigg(
	\sum_{x\in\mathcal{N}(\lceil an \rceil)}
	e^{-V(x)}1_{\left\{\min_{j\in [n, \lceil an \rceil]\cap \mathbb{Z}} V(x_j)<  \gamma_n \right\}} 1_{\{\min_{j\leq \lceil an\rceil } V(x_j)\geq -L \}}\bigg)\nonumber\\
	& = \delta + \frac{\left(\gamma_n+b\right) \sqrt{n}}{\varepsilon}
	\mathbf{P}_L
	\Big( \min_{n\leq j\leq \lceil an \rceil} S_j < \gamma_n+L , \min_{j\leq \lceil an\rceil } S_j \geq 0\Big).
\end{align}
Let $f_k(S_k):=\mathbf{P}_{S_k}\left(\min_{j\leq \lceil an \rceil-k } S_j \geq 0\right)$.
By Lemma \ref{lemma1} (i), we know that $f_k(S_k)\lesssim (1+S_k) (\lceil an\rceil -k)^{-1/2}$. By Lemma \ref{lemma1} (ii), we have
\begin{align*}
	& \mathbf{P}_L \Big( \min_{n\leq j\leq \lceil an \rceil} S_j < \gamma_n+L , \min_{j\leq \lceil an\rceil } S_j \geq 0\Big) \leq \sum_{k=n}^{\lceil an \rceil} \mathbf{P}_L \Big(  S_k < \gamma_n+L , \min_{j\leq \lceil an\rceil } S_j \geq 0\Big)\\
	& \leq \sum_{k=n}^{\lceil an \rceil - [\sqrt{n}]} \mathbf{E}_L \Big( f_k(S_k) 1_{\left\{S_k < \gamma_n+L\right\} }1_{\left\{ \min_{j\leq k } S_j \geq 0\right\}} \Big)+ \sum_{k=\lceil an \rceil - [\sqrt{n}]+1}^{\lceil an \rceil}\mathbf{P}_L \Big(  S_k < \gamma_n+L , \min_{j\leq k } S_j \geq 0\Big)\\
	& \lesssim  (1+\gamma_n +L)\sum_{k=n}^{\lceil an \rceil - [\sqrt{n}]}\frac{1}{\sqrt{\lceil an \rceil-k}} \mathbf{P}_L \Big(S_k < \gamma_n+L, \min_{j\leq k } S_j \geq 0 \Big) + \sqrt{n}\frac{(1+\gamma_n+L)^2(1+L)}{\sqrt{\left(\lceil an \rceil - [\sqrt{n}]+1\right)^3}}\\
	& \lesssim \frac{(1+\gamma_n +L)}{\sqrt{[\sqrt{n}]}}\cdot\lceil an \rceil \cdot  \frac{(1+\gamma_n +L)^2 (1+L)}{\sqrt{n^3}}+ \sqrt{n}\frac{(1+\gamma_n+L)^2(1+L)}{\sqrt{\left(\lceil an \rceil - [\sqrt{n}]+1\right)^3}}=:b_n.
\end{align*}
Using $\gamma_n\to\infty$ as $n\to\infty$, we can easily get
$n^{3/4}b_n/\left(\gamma_n^3\right) \lesssim 1.$ Combining this with
  \eqref{step_5} and the fact that $\gamma_n = o(n^{1/16})$, we get that as $n\to\infty$,
\[
\P\left(\left|J_n\right|> \varepsilon\right) \lesssim \delta + \frac{\left(\gamma_n+b\right) \sqrt{n}}{\varepsilon}b_n\lesssim \delta + \frac{\gamma_n^4}{n^{1/4}}.
\]
Letting $n\to\infty$ first and then $\delta \to 0$, we get the desired result.
\hfill$\Box$

\bigskip

Recall that $c^*$ and $\alpha^*$ are the constants in \eqref{Limit-of-R(y)/y} and \eqref{Second-Exp-of-R(y)} and $\ell (y)= c^* y +\alpha^*$.

\begin{prop}\label{prop3}
	For any $a\geq 1$, under $\mathbb{P}$,
	\[
	 \lim_{n\to\infty} \sqrt{n}\left|D_\infty^{\lceil an\rceil, \gamma_n} -\left(c^* D_{\lceil an \rceil} -\left(c^* \gamma_n -\alpha^*\right)W_{\lceil an \rceil}\right) \right| = 0 \qquad \mbox{in probability}.
	\]
\end{prop}
\textbf{Proof: } Define $\widetilde{D}_m^{\gamma_n}:=
\sum_{x\in\mathcal{N}(m)} \ell\left(V(x) -\gamma_n \right)e^{-V(x)} = c^* D_m -\left(c^* \gamma_n -\alpha^*\right)W_m$.
Note that
\begin{align}\label{step_6}
	\sqrt{n}\left|D_\infty^{\lceil an\rceil, \gamma_n} -\widetilde{D}_{\lceil an \rceil}^{\gamma_n} \right|  \leq \sqrt{n} \left|{D}_\infty^{\lceil an\rceil, \gamma_n} - {D}_{\lceil an \rceil}^{\lceil an\rceil, \gamma_n} \right| + \sqrt{n}\left|{D}_{\lceil an\rceil}^{\lceil an\rceil, \gamma_n} -\widetilde{D}_{\lceil an \rceil}^{\gamma_n} \right|.
\end{align}
By Lemma \ref{lemma2} (ii), $\lim_{n\to\infty}\sqrt{n} \left|{D}_\infty^{\lceil an\rceil, \gamma_n} - {D}_{\lceil an \rceil}^{\lceil an\rceil, \gamma_n} \right|=0$ in probability with respect to $\P$.
On  the set $\left\{\min_{|x| =\lceil an \rceil} V(x)>0 \right\}$, we have
\begin{align}
	\sqrt{n}\Big|D_{\lceil an\rceil}^{\lceil an\rceil, \gamma_n} -\widetilde{D}_{\lceil an \rceil}^{\gamma_n} \Big|& = \sqrt{n}\Big|	\sum_{x\in\mathcal{N}(\lceil an \rceil)}
	\ell \left(V(x)-\gamma_n\right)e^{-V(x)}1_{\left\{  V(x) <  \gamma_n \right\}}\Big| \nonumber\\ & \leq \left(c^* \gamma_n +\alpha^*\right)\sqrt{n}\left(W_{\lceil an \rceil} - \widetilde{W}_{\lceil an \rceil}^{\lceil an \rceil, \gamma_n}\right).
\end{align}
Combining Lemma \ref{lemma6} (i) with the fact that $\lim_{m\to\infty}\P\left(\min_{x\in \mathcal{N}(m)} V(x)>0\right)=1$, we immediately get that
$\lim_{n\to\infty}\sqrt{n}\left|{D}_{\lceil an\rceil}^{\lceil an\rceil, \gamma_n} -\widetilde{D}_{\lceil an \rceil}^{\gamma_n} \right|=0$ in probability with respect to $\P$. Now  the desired conclusion follows.
\hfill$\Box$

\subsection{Approximation of the martingales via weighted number of particles} \label{subsec4.2}

Recall that $\gamma_n$ is defined in \eqref{e:rs1}. For any $n\in \mathbb{N}$ and $a\geq 1$, we define
\begin{align*}
	N_{good}^{\lceil an \rceil, \gamma_n }:= \sum_{k=\lceil an \rceil+1}^\infty
	\sum_{x\in\mathcal{N}(k)}
	e^{-V(x)}1_{\left\{\min_{n\leq j\leq k-1} V(x_j)\geq \gamma_n, V(x)< \gamma_n \right\}}.
\end{align*}

\begin{lemma}\label{lemma9}
	For every $a\geq 1$, under $\P$,
	$$\lim_{n\to\infty} \sqrt{n}\beta_n \left|N_{good}^{\lceil an \rceil, \gamma_n }-\mathbb{E}\left( N_{good}^{\lceil an \rceil, \gamma_n }\Big|\mathcal{F}_{\lceil an \rceil} \right) \right|=0,\quad \mbox{in probability}.$$
\end{lemma}
\textbf{Proof: } Define
\begin{align}
	&\mathcal{A}_\kappa^{\lceil an \rceil, \gamma_n} := \Big\{x\in \mathcal{L}_{good}^{\lceil an \rceil ,\gamma_n}: \forall \lceil an \rceil+1 \leq j \leq |x|,\nonumber\\
	&\quad \sum_{u\in \Omega(x_j)} \left(1+\left(V(u)-V(x_{j-1})\right)_+\right)e^{-\left(V(u)-V(x_{j-1})\right)} \leq \kappa e^{\left(V(x_{j-1})-\gamma_n \right)/2} \Big\},
	\label{e:rs4}\\
		& N_{good,\kappa }^{\lceil an \rceil, \gamma_n }:= \sum_{\ell =\lceil an \rceil+1}^\infty \sum_{x\in\mathcal{N}(\ell )} e^{-V(x)}1_{\left\{\min_{n\leq j\leq \ell-1} V(x_j)\geq \gamma_n, V(x)< \gamma_n \right\}}1_{\big\{x \in \mathcal{A}_\kappa^{\lceil an \rceil, \gamma_n} \big\}}.\nonumber
 \end{align}
Then for any $\varepsilon>0$, we have
\begin{align}\label{Ineq-1}
	& \mathbb{P}\left(\sqrt{n}\beta_n\left|N_{good}^{\lceil an \rceil, \gamma_n }-\mathbb{E}\left( N_{good}^{\lceil an \rceil, \gamma_n }\Big|\mathcal{F}_{\lceil an \rceil} \right) \right| > 3\varepsilon \Big| \mathcal{F}_{\lceil a n \rceil}\right)\nonumber \\
	& \leq \mathbb{P}\left(\sqrt{n}\beta_n\left( N_{good}^{\lceil an \rceil, \gamma_n }-N_{good,\kappa }^{\lceil an \rceil, \gamma_n } \right) > \varepsilon \Big| \mathcal{F}_{\lceil a n \rceil}\right) \nonumber \\
	&\quad	+ \mathbb{P}\left(\sqrt{n}\beta_n\left|N_{good,\kappa }^{\lceil an \rceil, \gamma_n }-\mathbb{E}\left( N_{good,\kappa }^{\lceil an \rceil, \gamma_n }\Big|\mathcal{F}_{\lceil an \rceil} \right) \right| > \varepsilon \Big| \mathcal{F}_{\lceil a n \rceil}\right)+ 1_{\left\{\sqrt{n}\beta_n\mathbb{E}\left( N_{good}^{\lceil an \rceil, \gamma_n } -N_{good,\kappa }^{\lceil an \rceil, \gamma_n } \Big|\mathcal{F}_{\lceil an \rceil} \right)> \varepsilon \right\}}\nonumber \\
	& \leq \frac{2\sqrt{n}\beta_n}{\varepsilon }\mathbb{E}\left( N_{good}^{\lceil an \rceil, \gamma_n } -N_{good,\kappa }^{\lceil an \rceil, \gamma_n } \Big|\mathcal{F}_{\lceil an \rceil} \right)\nonumber \\ &\quad +\left(\frac{\sqrt{n}\beta_n}{\varepsilon }\right)^{1+\alpha} \mathbb{E}\left(\left|N_{good,\kappa }^{\lceil an \rceil, \gamma_n }-\mathbb{E}\left( N_{good,\kappa }^{\lceil an \rceil, \gamma_n }\Big|\mathcal{F}_{\lceil an \rceil} \right) \right|^{1+\alpha}\Big| \mathcal{F}_{\lceil an \rceil} \right).
\end{align}
Using Lemma \ref{lemma7} (ii) and the Markov property, we get
\begin{align}\label{Ineq-2}
	&\mathbb{E}\left( N_{good}^{\lceil an \rceil, \gamma_n } -N_{good,\kappa }^{\lceil an \rceil, \gamma_n } \Big|\mathcal{F}_{\lceil an \rceil} \right) \nonumber \\&
		=  \mathbb{E} \bigg(\sum_{\ell =\lceil an \rceil+1}^\infty \sum_{x\in\mathcal{N}(\ell )} e^{-V(x)}1_{\left\{\min_{n\leq j\leq \ell -1} V(x_j)\geq \gamma_n, V(x)< \gamma_n \right\}}1_{\big\{x \notin \mathcal{A}_\kappa^{\lceil an \rceil, \gamma_n} \big\}}\bigg| \mathcal{F}_{\lceil an \rceil}\bigg) \nonumber \\
		& \leq \mathbb{E} \bigg(\sum_{\ell =\lceil an \rceil+1}^\infty \sum_{x\in\mathcal{N}(\ell)} e^{-V(x)}1_{\left\{\min_{\lceil an \rceil \leq j\leq \ell-1} V(x_j)\geq \gamma_n, V(x)< \gamma_n \right\}}1_{\big\{x \notin \mathcal{A}_\kappa^{\lceil an \rceil, \gamma_n} \big\}}\bigg| \mathcal{F}_{\lceil an \rceil}\bigg) \nonumber
	\\& =	\sum_{u\in\mathcal{N}(\lceil an \rceil)}
	e^{-V(u)}1_{\left\{V(u)\geq  \gamma_n \right\}} \mathbb{E}\left(N_{[1,\infty)}^{z}-  N_{[1,\infty),\kappa}^{z}\right)\Big|_{z=V(u)-\gamma_n }\lesssim \frac{g(\kappa)}{\log \kappa }W_{\lceil an \rceil},
\end{align}
where $g$ is the function in Lemma \ref{lemma7} (ii).
It follows from \cite[Theorem 2]{vonB-E} that, for any random variables $X_1, \dots, X_m$ with finite $(1+\alpha)$-moment satisfying
$\mathbb{E}(X_{j+1}|X_1+\cdots X_j)=0$
for all $j=1, \dots, m-1$, we have $\mathbb{E}\left( \left|\sum_{i} X_i \right|^{1+\alpha} \right)\leq 2\sum_{i}\mathbb{E}\left(|X_i|^{1+\alpha}\right)$.
Note that,  by the branching property, $N_{good,\kappa }^{\lceil an \rceil, \gamma_n }= \sum_{u\in \mathcal{N} (\lceil an \rceil)} H(u)$ with $\{H(u): u\in \mathcal{N} (\lceil an \rceil)\}$ being independent random variables conditioned on $\mathcal{F}_{\lceil an \rceil}$.
More precisely, for $u\in \mathcal{N} (\lceil an \rceil)$,
\[
H(u)= \sum_{\ell =\lceil an \rceil+1}^\infty \sum_{x\in\mathcal{N}(\ell ), x> u} e^{-V(x)}1_{\left\{\min_{n\leq j\leq \ell-1} V(x_j)\geq \gamma_n, V(x)< \gamma_n \right\}}1_{\big\{x \in \mathcal{A}_\kappa^{\lceil an \rceil, \gamma_n} \big\}}.
\]
Using these two observations, the branching property, the trivial inequality $\mathbb{E}\left(|X- \mathbb{E} X|^{1+\alpha}\right)\leq 2^{1+\alpha}\mathbb{E}(|X|^{1+\alpha})$  and  Lemma \ref{lemma8}, we get
\begin{align}\label{Ineq-3}
	&\mathbb{E}\left(\left|N_{good,\kappa }^{\lceil an \rceil, \gamma_n }-\mathbb{E}\left( N_{good,\kappa }^{\lceil an \rceil, \gamma_n }\Big|\mathcal{F}_{\lceil an \rceil} \right) \right|^{1+\alpha}\Big| \mathcal{F}_{\lceil an \rceil} \right)	\nonumber \\
	&= \mathbb{E}\bigg(\bigg|	 \sum_{u\in\mathcal{N}(\lceil an \rceil)}	 \left(H(u)- \mathbb{E}\left(H(u)\big| \mathcal{F}_{\lceil an \rceil}\right)\right) \bigg|^{1+\alpha}\Big| \mathcal{F}_{\lceil an \rceil} \bigg) \nonumber \\
    &\lesssim      \sum_{u\in\mathcal{N}(\lceil an \rceil)}
     \mathbb{E}\left(\left|H(u)- \mathbb{E}\left(H(u)\big| \mathcal{F}_{\lceil an \rceil}\right) \right|^{1+\alpha}\Big| \mathcal{F}_{\lceil an \rceil} \right) \lesssim
     \sum_{u\in\mathcal{N}(\lceil an \rceil)}
     \mathbb{E}\left(\left|H(u)\right|^{1+\alpha}\Big| \mathcal{F}_{\lceil an \rceil} \right) \nonumber  \\
	 & \lesssim   	\sum_{u\in\mathcal{N}(\lceil an \rceil)}
	e^{-(1+\alpha)V(u)}1_{\left\{V(u)\geq \gamma_n \right\}}\mathbb{E}\Big(\left(N_{[1,\infty),\kappa}^{z} \right)^{1+\alpha}\Big)\Big|_{z=V(u)-\gamma_n } \nonumber \\
	& \lesssim  \kappa^\alpha 	\sum_{u\in\mathcal{N}(\lceil an \rceil)}
	e^{-(1+\alpha)V(u)} e^{\alpha \left(V(u)-\gamma_n \right)}= \kappa^\alpha e^{-\alpha\gamma_n }W_{\lceil an\rceil}.
\end{align}
Combining \eqref{Ineq-1}, \eqref{Ineq-2} and \eqref{Ineq-3}, taking $\kappa = e^{\beta_n /2}$ and using the definition of $\gamma_n$, we get
\begin{align}\label{step_43}
	&\mathbb{P}\left(\sqrt{n}\beta_n\left|N_{good}^{\lceil an \rceil, \gamma_n }-\mathbb{E}\left( N_{good}^{\lceil an \rceil, \gamma_n }\Big|\mathcal{F}_{\lceil an \rceil} \right) \right| > 3\varepsilon\Big| \mathcal{F}_{\lceil an \rceil}\right)\nonumber \\
	& \lesssim \sqrt{n}W_{\lceil an \rceil}\bigg(\frac{g\left(e^{\beta_n /2}\right)}{\varepsilon} +\frac{ \beta_n^{1+\alpha} e^{-\alpha\beta_n /2}}{\varepsilon^{1+\alpha}}\bigg).
\end{align}
Letting $n\to\infty$, and combining \eqref{Seneta-Heyde-scaling}, \eqref{step_4} and the fact that $\lim_{z\to\infty}g(z)=0$, we  get the desired conclusion.
\hfill$\Box$

\bigskip

Define
\begin{align}\label{Widetilde-N}
\widetilde{N}_{good}^{\lceil an \rceil, \gamma_n}:= \sum_{k=\lceil an \rceil+1}^\infty \sum_{x\in\mathcal{N}(k)}
e^{-V(x)}1_{\left\{\min_{n\leq j\leq k-1} V(x_j)\geq \gamma_n, \gamma_n -\beta_n /2 \leq  V(x)< \gamma_n \right\}}.
\end{align}

\begin{prop}\label{prop2}
For any $a\geq 1$, under $\P$, it holds that
	$$\lim_{n\to\infty} \sqrt{n}\beta_n \left|\widetilde{N}_{good}^{\lceil an \rceil, \gamma_n}- W_{\lceil an \rceil} \right|=0,\quad \mbox{in probability}.$$
\end{prop}
\textbf{Proof:} By the branching property and Lemma \ref{lemma7} (i), we have
\begin{align}\label{step_21}
	\mathbb{E}\left( N_{good}^{\lceil an \rceil, \gamma_n }\Big|\mathcal{F}_{\lceil an \rceil} \right) & = \mathbb{E}\bigg( \sum_{k=\lceil an \rceil+1}^\infty
	\sum_{x\in\mathcal{N}(k)}
	e^{-V(x)}1_{\left\{\min_{n\leq j\leq k-1} V(x_j)\geq \gamma_n, V(x)< \gamma_n \right\}} \Big| \mathcal{F}_{\lceil an \rceil}\bigg)\nonumber\\
	& =  \sum_{u\in\mathcal{N}(\lceil an \rceil)}
	 e^{-V(u)}1_{\left\{\min_{n\leq j\leq \lceil an \rceil} V(u_j)\geq \gamma_n \right\}}\mathbb{E}\big(N_{[1,\infty)}^{z}\big)\big|_{z=V(u)-\gamma_n}\nonumber\\
	& =  \sum_{u\in\mathcal{N}(\lceil an \rceil)}
	e^{-V(u)}1_{\left\{\min_{n\leq j\leq \lceil an \rceil} V(u_j)\geq \gamma_n \right\}} = \widetilde{W}_{\lceil an \rceil}^{n,\gamma_n }.
\end{align}
Then we have
\begin{align*}
	\sqrt{n}\beta_n \left|\widetilde{N}_{good}^{\lceil an \rceil, \gamma_n}- W_{\lceil an \rceil} \right| & \leq \sqrt{n}\beta_n \left|N_{good}^{\lceil an \rceil, \gamma_n}- \widetilde{N}_{good}^{\lceil an \rceil, \gamma_n} \right|+ \sqrt{n}\beta_n \left|N_{good}^{\lceil an \rceil, \gamma_n }-\mathbb{E}\left( N_{good}^{\lceil an \rceil, \gamma_n }\Big|\mathcal{F}_{\lceil an \rceil} \right) \right|\\
	& \qquad + \sqrt{n}\beta_n \left|\widetilde{W}_{\lceil an \rceil}^{n,\gamma_n }- W_{\lceil an \rceil} \right|.
\end{align*}
It follows from Lemma \ref{lemma6} (ii) and Lemma \ref{lemma9} that the second and third terms
on the right hand side of the above inequality  converge to $0$ in probability as $n\to\infty$. To prove the desired result, we only need to prove
\begin{align}\label{step_22}
	\lim_{n\to\infty} \sqrt{n}\beta_n \left|N_{good}^{\lceil an \rceil, \gamma_n}- \widetilde{N}_{good}^{\lceil an \rceil, \gamma_n} \right|=0,\quad \mbox{in probability}.
\end{align}

Recall that $(S_n, \mathbf{P}_y)$
is the random walk defined in \eqref{step_1} with $S_0=y$. By the Markov property and \eqref{step_1}, we have
\begin{align}\label{step_49}
		&\mathbb{E}\left( \sqrt{n}\beta_n \left(N_{good}^{\lceil an \rceil, \gamma_n}- \widetilde{N}_{good}^{\lceil an \rceil, \gamma_n} \right) \Big| \mathcal{F}_{\lceil an \rceil }\right)
		 \nonumber \\& = \sqrt{n}\beta_n \mathbb{E}\bigg(\sum_{k=\lceil an \rceil+1}^\infty
		   \sum_{x\in\mathcal{N}(k)}
		e^{-V(x)}1_{\left\{\min_{n\leq j\leq k-1} V(x_j)\geq \gamma_n, V(x)< \gamma_n-\beta_n/2 \right\}} \Big| \mathcal{F}_{\lceil an \rceil}\bigg)\nonumber \\
	& \leq \sqrt{n}\beta_n 	   \sum_{u\in\mathcal{N}(\lceil an \rceil)}
	e^{-V(u)} 1_{\left\{V(u)\geq \gamma_n \right\}}\sum_{k=1}^\infty
		\mathbf{P}_y
	\Big(\min_{j\leq k-1} S_j \geq 0, S_k < -\beta_n /2\Big)\Big|_{y = V(u)-\gamma_n}.
\end{align}
For $k\ge 2$, $U_k:= S_k- S_{k-1}$ is independent to $S_1,...,S_{k-1}$ and has the same law as $S_1-S_0$. Thus for all $y\geq 0$,
\begin{align*}
&\sum_{k=1}^\infty \mathbf{P}_y \big(\min_{j\leq k-1} S_j \geq 0, S_k < -\beta_n /2\big) =  \mathbf{E}_y \Big(\sum_{k=1}^\infty\mathbf{P}_y \Big(\min_{j\leq k-1} S_j \geq 0, S_{k-1} < -u_k-\beta_n /2\Big)\Big|_{u_k = U_k} \Big)\\
	& =   \mathbf{E}_y \Big(\sum_{k=1}^\infty \mathbf{P}_y \Big(\min_{j\leq k-1} S_j \geq 0, S_{k-1} < -u-\beta_n /2\Big)\Big|_{u = S_1-y} \Big).
\end{align*}
By  Lemma \ref{lemma1} (iii), we see that
\[
\sum_{k=1}^\infty
 \mathbf{P}_y
\big(\min_{j\leq k-1} S_j \geq 0, S_{k-1} < -u-\beta_n /2\big) \lesssim (1-u-\beta_n/2)^21_{\{-u-\beta_n/2 >0\}}.
\]
Therefore, using the fact that
 $(S_1-y, \mathbf{P}_y)\stackrel{\mathrm{d}}{=} (S_1, \mathbf{P})$,
we obtain
\begin{align}
\sum_{k=1}^\infty \mathbf{P}_y \big(\min_{j\leq k-1} S_j \geq 0, S_k < -\beta_n /2\big) 	\lesssim  \mathbf{E}\big( \left(1-S_1-\beta_n/2 \right)^21_{\left\{S_1<-\beta_n/2 \right\}}\big).
\end{align}
Note that for $u<-\beta_n /2$ and $\beta_n /2 >1$,
\[
\left(1-u-\beta_n/2\right)^2 \leq 2\big(u^2 + \left(\beta_n/2 -1\right)^2\big)\leq 4u^2.
\]
Hence, for all $y\geq 0$ and $n\geq 1$,
\begin{align}\label{step_17}
	\sum_{k=1}^\infty \mathbf{P}_y \big(\min_{j\leq k-1} S_j \geq 0, S_k < -\beta_n /2\big) \lesssim 4 \mathbf{E}\left(S_1^2 1_{\left\{S_1 < -\beta_n/2\right\}}\right)\lesssim \frac{1}{\beta_n}\mathbf{E}\big(\left(-S_1\right)^31_{\left\{S_1 < -\beta_n/2 \right\}}\big).
\end{align}
Combining \eqref{step_49} and \eqref{step_17}, we obtain
\begin{align*}
	&\mathbb{E}
	\left( \sqrt{n}\beta_n \left(N_{good}^{\lceil an \rceil, \gamma_n}- \widetilde{N}_{good}^{\lceil an \rceil, \gamma_n} \right) \Big| \mathcal{F}_{\lceil an \rceil }\right) \lesssim \sqrt{n}W_{\lceil an \rceil}
	 \mathbf{E}
	\left(\left(-S_1\right)^3 1_{\left\{S_1 < -\beta_n/2 \right\}}\right),
\end{align*}
which implies \eqref{step_22}. The proof is complete.
\hfill$\Box$

\begin{cor}\label{cor1}
	Let $m\geq 1, 1\leq a_1 < ... < a_m < a_{m+1}=\infty$ and $\left(z^n\right)_{n\in \mathbb{N}}\in\left(\mathbb{C}^m\right)^{\mathbb{N}}$, here $z^n=(z^n_1, z^n_2,\cdots, z^n_m)$.
	Assume that for all $1\leq k\leq m$ and $n\geq 0$, $\mathrm{Re}\left(z_k^n\right)\leq 0$ and  $z^n$ converges to some  $z=(z_1, \cdots, z_m)\in\mathbb{C}^m$.
Then under $\P$,
	\begin{align*}
		&\lim_{n\to\infty} \mathbb{E} \bigg(\exp\bigg\{\sum_{k=1}^m z_k^n \sqrt{n}\left(\widetilde{N}_{good}^{\lceil a_k n\rceil, \gamma_n }- \widetilde{N}_{good}^{\lceil a_{k+1} n\rceil, \gamma_n }\right) \bigg\}\bigg| \mathcal{F}_n\bigg) \\
		& =\exp\bigg\{\sqrt{\frac{2}{\pi \sigma^2}}D_\infty \sum_{k=1}^m z_k \bigg(\frac{1}{\sqrt{a_k}}- \frac{1}{\sqrt{a_{k+1}}}\bigg) \bigg\},\quad \qquad \qquad\quad\quad\quad\mbox{in probability}.
	\end{align*}
\end{cor}
\textbf{Proof:}
Using \eqref{Seneta-Heyde-scaling} and Proposition \ref{prop2}, we have that, for all $1\leq a<b \leq \infty$,
$$ \lim_{n\to\infty} \bigg|\sqrt{n}\left(\widetilde{N}_{good}^{\lceil an \rceil, \gamma_n} -\widetilde{N}_{good}^{\lceil bn \rceil, \gamma_n} \right)- \sqrt{\frac{2}{\pi \sigma^2}}D_n \left(\frac{1}{\sqrt{a}} - \frac{1}{\sqrt{b}}\right) \bigg|=0,\quad \mbox{in probability}.
$$
Noticing that
$\widetilde{N}_{good}^{\lceil an \rceil, \gamma_n} \geq \widetilde{N}_{good}^{\lceil bn \rceil, \gamma_n}$
 and $\mathrm{Re}\left(z_k^n\right)\leq 0$, using \cite[Remark A.3]{MP}, it suffices to prove that
\begin{align*}
	&\lim_{n\to\infty} \mathbb{E} \bigg(\exp\bigg\{\sum_{k=1}^m z_k^n \sqrt{\frac{2}{\pi \sigma^2}}D_n \left(\frac{1}{\sqrt{a_k}} - \frac{1}{\sqrt{a_{k+1}}}\right) \bigg\}\bigg| \mathcal{F}_n\bigg) \\ & =\exp\bigg\{\sqrt{\frac{2}{\pi \sigma^2}}D_\infty \sum_{k=1}^m z_k \bigg(\frac{1}{\sqrt{a_k}}- \frac{1}{\sqrt{a_{k+1}}}\bigg) \bigg\},\quad \qquad \qquad\quad\quad\quad\mbox{in probability}.
\end{align*}
Since $D_n \in \mathcal{F}_n $ and $\lim_{n\to\infty}D_n = D_\infty$, the equality above is trivial.
\hfill$\Box$

\bigskip

Define
\begin{equation}\label{e:rs5}
 \widehat{N}_{good}^{\lceil an \rceil, \gamma_n}:= \sum_{k=\lceil an \rceil +1}^\infty
 \sum_{x\in \mathcal{N}(k)}
 c^* e^{-V(x)}\Big(V(x)-\frac{\log n}{2}\Big) 1_{\left\{\min_{n\leq j\leq k-1} V(x_j)\geq \gamma_n \right\}}1_{\left\{ \gamma_n -\beta_n/2\leq V(x) <\gamma_n \right\}}.
\end{equation}

\begin{lemma}\label{lemma11}
For any $a\geq 1$, under $\P$, it holds that
	$$\lim_{n\to\infty} \sqrt{n}\left| \widehat{N}_{good}^{\lceil an \rceil, \gamma_n} -\left(c^* \beta_n -\alpha^*\right)\widetilde{N}_{good}^{\lceil an \rceil, \gamma_n}\right|=0,\quad \mbox{in probability}.$$
\end{lemma}
\textbf{Proof:}  Note that
\begin{align*}
&\sqrt{n}\left| \widehat{N}_{good}^{\lceil an \rceil, \gamma_n} -\left(c^* \beta_n -\alpha^*\right)\widetilde{N}_{good}^{\lceil an \rceil, \gamma_n}\right|\\
&\leq \sqrt{n}\left|\widehat{N}_{good}^{\lceil an \rceil, \gamma_n} -\mathbb{E}\left( \widehat{N}_{good}^{\lceil an \rceil, \gamma_n }\Big|\mathcal{F}_{\lceil an \rceil} \right) \right| \\
& \quad + \sqrt{n}\left|\mathbb{E}\left( \widehat{N}_{good}^{\lceil an \rceil, \gamma_n }\Big|\mathcal{F}_{\lceil an \rceil} \right) -\left(c^* \beta_n -\alpha^*\right)\mathbb{E}\left( N_{good}^{\lceil an \rceil, \gamma_n }\Big|\mathcal{F}_{\lceil an \rceil} \right)\right|\\
&\quad +\left(c^* \beta_n -\alpha^*\right)\sqrt{n}\left|\mathbb{E}\left( N_{good}^{\lceil an \rceil, \gamma_n }\Big|\mathcal{F}_{\lceil an \rceil} \right) - N_{good}^{\lceil an \rceil, \gamma_n }\right|\\
&\quad + \left(c^* \beta_n -\alpha^*\right)\sqrt{n}\left|N_{good}^{\lceil an \rceil, \gamma_n } - \widetilde{N}_{good}^{\lceil an \rceil, \gamma_n }\right|\\
&=:I_n+II_n+III_n+IV_n.
\end{align*}
It follows from Lemma \ref{lemma9} and \eqref{step_22} that $III_n$ and $IV_n$ tend to $0$ in probability as $n\to\infty$.

We first show that $\lim_{n\to\infty}I_n=0$ in probability.  Define
\begin{align}\label{Widehat-N}
	&\widehat{N}_{good, \kappa}^{\lceil an \rceil, \gamma_n}
		:= \sum_{\ell =\lceil an \rceil +1}^\infty 	\sum_{x\in\mathcal{N}(\ell )} c^* e^{-V(x)}\Big(V(x)-\frac{\log n}{2}\Big) \nonumber \\
	&\quad\quad\times  1_{\left\{\min_{n\leq j\leq \ell -1} V(x_j)\geq \gamma_n \right\}}1_{\left\{ \gamma_n -\beta_n/2\leq V(x) <\gamma_n \right\}}1_{\big\{x \in \mathcal{A}_\kappa^{\lceil an \rceil, \gamma_n} \big\}},
\end{align}
where $\mathcal{A}_\kappa^{\lceil an \rceil, \gamma_n}$ is defined in \eqref{e:rs4}.
Using an argument similar to that in first part of proof of Lemma \ref{lemma9}, we get that, for any $\varepsilon>0$,
\begin{align*}
	&\mathbb{P}\left(\sqrt{n}\left|\widehat{N}_{good}^{\lceil an \rceil, \gamma_n }-\mathbb{E}\left( \widehat{N}_{good}^{\lceil an \rceil, \gamma_n }\Big|\mathcal{F}_{\lceil an \rceil} \right) \right| > 3\varepsilon\Big| \mathcal{F}_{\lceil an \rceil}\right)\\
	& \leq \frac{2\sqrt{n}}{\varepsilon }\mathbb{E}\left( \widehat{N}_{good}^{\lceil an \rceil, \gamma_n } -\widehat{N}_{good,\kappa }^{\lceil an \rceil, \gamma_n } \Big|\mathcal{F}_{\lceil an \rceil} \right) +\left(\frac{\sqrt{n}}{\varepsilon }\right)^{1+\alpha} \mathbb{E}\left(\left| \widehat{N}_{good,\kappa }^{\lceil an \rceil, \gamma_n }- \mathbb{E}\left(\widehat{N}_{good,\kappa }^{\lceil an \rceil, \gamma_n } \Big| \mathcal{F}_{\lceil an \rceil} \right) \right|^{1+\alpha}\Big| \mathcal{F}_{\lceil an \rceil} \right).
\end{align*}
By an argument similar to that leading to \eqref{step_43}, taking $\kappa =e^{\beta_n/2}$ and using Lemma \ref{lemma10}, we can get
\begin{align*}
	&\mathbb{P}\left(\sqrt{n}\left|\widehat{N}_{good}^{\lceil an \rceil, \gamma_n }-\mathbb{E}\left( \widehat{N}_{good}^{\lceil an \rceil, \gamma_n }\Big|\mathcal{F}_{\lceil an \rceil} \right) \right| > 3\varepsilon\Big| \mathcal{F}_{\lceil an \rceil}\right)\\
	& \lesssim \sqrt{n}W_{\lceil an \rceil}\bigg(\frac{g\left(e^{\beta_n /2}\right)}{\varepsilon} +\frac{ \beta_n^{1+\alpha} e^{-\alpha\beta_n /2}}{\varepsilon^{1+\alpha}}\bigg).
\end{align*}
Combining this with \eqref{Seneta-Heyde-scaling}, \eqref{step_4} and the fact $\lim_{z\to\infty}g(z)=0$, we immediately
get $\lim_{n\to\infty}I_n=0$ in probability.

Therefore, it remains to prove that $\lim_{n\to\infty}II_n=0$ in probability.
By the branching property, we have
	\begin{align}\label{step_51}
	&\mathbb{E}\left( \widehat{N}_{good}^{\lceil an \rceil, \gamma_n }\Big|\mathcal{F}_{\lceil an \rceil} \right)
	= c^* 	\sum_{u\in\mathcal{N}(\lceil an \rceil)}
	e^{-V(u)}1_{\left\{\min_{n\leq j \leq \lceil an \rceil} V(u_j)\geq \gamma_n \right\}}
	\nonumber\\
	& \quad \quad \times \sum_{k=1}^\infty\mathbb{E}\bigg(	\sum_{x\in\mathcal{N}(k)}
	e^{-V(x)}\left(V(x)+\beta_n +y\right)1_{\left\{\min_{j\leq k-1}V(x_j)\geq -y \right\}}1_{\left\{-y-\beta_n/2 \leq V(x)<-y \right\}} \bigg)\bigg|_{y= V(u)-\gamma_n}
	\nonumber\\
	& = c^* 	\sum_{u\in\mathcal{N}(\lceil an \rceil)}
	e^{-V(u)}1_{\left\{\min_{n\leq j \leq \lceil an \rceil} V(u_j)\geq \gamma_n \right\}} \nonumber\\
	&\quad\quad\times \sum_{k=1}^\infty
		 \mathbf{E}_{V(u)-\gamma_n}
	\left(\left(S_k+\beta_n \right)1_{\left\{\min_{j\leq k-1}S_j\geq 0 \right\}}1_{\left\{-\beta_n/2 \leq S_k<0\right\}} \right).
\end{align}
For $y\geq 0, n\geq 1$, define
\[
\Lambda_n(y)
:= \sum_{k=1}^\infty
 \mathbf{E}_{y}
\left(\left(S_k+\beta_n \right)1_{\left\{\min_{j\leq k-1}S_j\geq 0 \right\}}1_{\left\{S_k<-\beta_n/2 \right\}} \right).
\]
By \eqref{step_17},  we have
\[
\Lambda_n(y)
\leq \beta_n \sum_{k=1}^\infty\mathbf{P}_y \left(\min_{j\leq k-1} S_j \geq 0, S_k < -\beta_n /2\right)\lesssim \mathbf{E}\left(\left(-S_1\right)^31_{\left\{S_1 < -\beta_n/2 \right\}}\right).
\]
On the other hand,   we have
\begin{align*}
-\Lambda_n(y) &\leq  \sum_{k=1}^\infty\mathbf{E}_{y}\left(\left(-S_k \right)1_{\left\{\min_{j\leq k-1}S_j\geq 0 \right\}}1_{\left\{S_k<-\beta_n/2 \right\}} \right)\nonumber\\
	 & = \sum_{k=1}^\infty\mathbf{E}_{y}\left(\left(-S_{k}\right)1_{\left\{\min_{j\leq k-1}S_j\geq 0 \right\}}1_{\left\{S_{k-1}<-U_k-\beta_n/2 \right\}} \right)\nonumber\\
	 & \leq \sum_{k=1}^\infty\mathbf{E}_{y}\left(\left(-U_k\right)1_{\left\{\min_{j\leq k-1}S_j\geq 0 \right\}}1_{\left\{S_{k-1}<-U_k-\beta_n/2 \right\}} \right),
\end{align*}
where $U_k= S_k-S_{k-1}$, and in the last inequality we used the fact that $-S_k \leq -U_k$ on the set $\{S_{k-1}\geq 0\}$.
By Lemma \ref{lemma1} (iii), we have
\begin{align*}
	&\sum_{k=1}^\infty\mathbf{E}_{y}\big(\left(-U_k\right) 1_{\left\{\min_{j\leq k-1}S_j\geq 0 \right\}}
	1_{\left\{S_{k-1}<-U_k-\beta_n/2 \right\}} \big)\\
	& = \mathbf{E}\Big(\left(-U_1\right)1_{\{-U_1 > \beta_n /2\}}\sum_{k=1}^\infty \mathbf{E}_y \Big(\min_{j\leq k-1} S_j \geq 0, S_{k-1}< -u_1-\beta_n /2 \Big)\Big|_{u_1=U_1} \Big)\\
	& \lesssim \mathbf{E}\big(\left(-U_1\right)\left(1- U_1-\beta_n /2\right)^21_{\{-U_1 > \beta_n /2\}}\big)\\
	&\lesssim \mathbf{E}\big(\left(-U_1\right)^31_{\{-U_1 > \beta_n /2\}}\big)= \mathbf{E}\big(\left(-S_1\right)^31_{\left\{S_1 < -\beta_n/2 \right\}}\big).
\end{align*}
Therefore,  for all $y\geq 0, n\geq 1$,
\[
 \left|\Lambda_n(y) \right|
 \lesssim \mathbf{E}\big(\left(-S_1\right)^31_{\left\{S_1 < -\beta_n/2 \right\}}\big).
\]
Combining this inequality and the definition of  $\Lambda_n$, we get
\begin{align}\label{step_24}
	&c^* 	\sum_{u\in\mathcal{N}(\lceil an \rceil)}
	e^{-V(u)}1_{\left\{\min_{n\leq j \leq \lceil an \rceil} V(u_j)\geq \gamma_n \right\}}\bigg| \sum_{k=1}^\infty
 \mathbf{E}_{V(u)-\gamma_n}
\left(\left(S_k+\beta_n \right)1_{\left\{\min_{j\leq k-1}S_j\geq 0 \right\}}1_{\left\{S_k<-\beta_n/2 \right\}} \right)\bigg|\nonumber\\
	&=c^*
	 \sum_{u\in\mathcal{N}(\lceil an \rceil)}
	 e^{-V(u)}1_{\left\{\min_{n\leq j \leq \lceil an \rceil} V(u_j)\geq \gamma_n \right\}}\left| 	\Lambda_n(V(u)-\gamma_n)\right|\nonumber \\
	&\lesssim \mathbf{E}
	\big(\left(-S_1\right)^31_{\left\{S_1 < -\beta_n/2 \right\}}\big) W_{\lceil an \rceil}.
\end{align}
By \eqref{step_53},  we see that
$R(y)=c^*y -c^*\mathbf{E}_y\big(S_{\tau_0^-}\big)$ for $y\geq 0$. Thus
	\begin{align}\label{step_50}
	&c^* 	\sum_{u\in\mathcal{N}(\lceil an \rceil)}
	e^{-V(u)}1_{\left\{\min_{n\leq j \leq \lceil an \rceil} V(u_j)\geq \gamma_n \right\}}\sum_{k=1}^\infty
		\mathbf{E}_{V(u)-\gamma_n}
	\big(\left(S_k+\beta_n \right)1_{\left\{\min_{j\leq k-1}S_j\geq 0 \right\}}1_{\left\{ S_k<0\right\}} \big)
	\nonumber\\
	& = 	\sum_{u\in\mathcal{N}(\lceil an \rceil)}
	e^{-V(u)}1_{\left\{\min_{n\leq j \leq \lceil an \rceil} V(u_j)\geq \gamma_n \right\}}c^*
	 	  \mathbf{E}_{V(u)-\gamma_n}
	 \big(S_{\tau_0^-}+\beta_n \big)
	\nonumber\\
	& = 	\sum_{u\in\mathcal{N}(\lceil an \rceil)}
	e^{-V(u)}1_{\left\{\min_{n\leq j \leq \lceil an \rceil} V(u_j)\geq \gamma_n \right\}} \left(c^*\beta_n +c^*\left(V(u)-\gamma_n\right)-R\left(V(u)-\gamma_n\right)\right).
	\end{align}
For any $\varepsilon>0,$ let $K$ be large enough such that $\left|R(y)- (c^*y +\alpha^*)\right|<\varepsilon$  for $y>K$. Therefore, when $V(u)-\gamma_n > K$, we have
$$\left|\left(c^*\beta_n +c^*\left(V(u)-\gamma_n\right)-R\left(V(u)-\gamma_n\right)\right)-\left(c^*\beta_n -\alpha^*\right)\right| < \varepsilon.$$
Recall that by $ \eqref{step_21}$,
\[
\left(c^* \beta_n -\alpha^*\right)\mathbb{E}\left( N_{good}^{\lceil an \rceil, \gamma_n }\Big|\mathcal{F}_{\lceil an \rceil} \right)=\left(c^* \beta_n -\alpha^*\right) \sum_{u\in\mathcal{N}(\lceil an \rceil)}
e^{-V(u)}1_{\left\{\min_{n\leq j\leq \lceil an \rceil} V(u_j)\geq \gamma_n \right\}}.
\]
 Combining this with \eqref{step_51}, \eqref{step_24} and \eqref{step_50},  we get
\begin{align}\label{step_25}
&II_n=
\sqrt{n}\left|\mathbb{E}\left( \widehat{N}_{good}^{\lceil an \rceil, \gamma_n }\Big|\mathcal{F}_{\lceil an \rceil} \right) -\left(c^* \beta_n -\alpha^*\right)\mathbb{E}\left( N_{good}^{\lceil an \rceil, \gamma_n }\Big|\mathcal{F}_{\lceil an \rceil} \right)\right|\nonumber\\
	& \lesssim
		\mathbf{E}
	\left(\left(-S_1\right)^31_{\left\{S_1 < -\beta_n/2 \right\}}\right) \sqrt{n}W_{\lceil an \rceil}+ \varepsilon  \sqrt{n}	\sum_{u\in\mathcal{N}(\lceil an \rceil)}
	e^{-V(u)}1_{\left\{\min_{n\leq j \leq \lceil an \rceil} V(u_j)\geq \gamma_n \right\}}1_{\left\{V(u)> \gamma_n +K\right\}}\nonumber\\
	& \ + \sup_{y\in[0,K]} \left(\left|R(y)-(c^*y+\alpha^*)\right|\right)  \sqrt{n}
	\sum_{u\in\mathcal{N}(\lceil an \rceil)}
	e^{-V(u)}1_{\left\{\min_{n\leq j \leq \lceil an \rceil} V(u_j)\geq \gamma_n \right\}}1_{\left\{V(u)\leq \gamma_n +K\right\}}.
\end{align}
Let $L$ and $\delta$ be the  constants in \eqref{Global-Minimum}. Then for  any $\theta>0$,
\begin{align}\label{step_26}
	&\mathbb{P}\bigg(\sqrt{n} 	\sum_{u\in\mathcal{N}(\lceil an \rceil)}
	e^{-V(u)}1_{\left\{\min_{n\leq j \leq \lceil an \rceil} V(u_j)\geq \gamma_n \right\}}1_{\left\{V(u)\leq \gamma_n +K\right\}} > \theta\bigg)\nonumber \\
	& \leq \delta + \mathbb{P}\bigg(\sqrt{n}
	\sum_{u\in\mathcal{N}(\lceil an \rceil)}
	e^{-V(u)}1_{\left\{\min_{ j \leq \lceil an \rceil} V(u_j)\geq -L \right\}}1_{\left\{V(u)\leq \gamma_n +K\right\}} > \theta\bigg)\nonumber\\
	& \leq \delta + \frac{\sqrt{n}}{\theta}\mathbb{E}\bigg(
	\sum_{u\in\mathcal{N}(\lceil an \rceil)}
	e^{-V(u)}1_{\left\{\min_{ j \leq \lceil an \rceil} V(u_j)\geq -L \right\}}1_{\left\{V(u)\leq \gamma_n +K\right\}} \bigg)\nonumber
	\\& =  \delta + \frac{\sqrt{n}}{\theta}
		\mathbf{P}
	\Big( \min_{ j \leq \lceil an \rceil} S_j\geq -L , S_{\lceil an \rceil}\leq \gamma_n +K \Big)\nonumber \\
	&\lesssim \delta + \frac{\sqrt{n}}{\theta} \frac{(\gamma_n + L + K)^2 (L+1)}{\sqrt{\lceil an \rceil^3}},
\end{align}
where in the last line we used  Lemma \ref{lemma1} (ii).
Letting $n\to\infty$ first and then $\delta \to 0$ in \eqref{step_26}, we get that the third term on the  right-side hand of \eqref{step_25} converges to $0$ in probability as $n\to\infty$. Note that the second term on the right-side hand is bounded by $\varepsilon \sqrt{n}W_{\lceil an \rceil}$,
so letting $n\to\infty$ first and then $\varepsilon \to 0$ in \eqref{step_25},
by \eqref{Seneta-Heyde-scaling},  we get that $\lim_{n\to\infty}II_n=0$ in probability.
The proof is now complete.
\hfill$\Box$

\begin{cor}\label{cor2}
	For any $a\geq 1$, under $\P$,
	$$\lim_{n\to\infty} \sqrt{n}\Big| \widehat{N}_{good}^{\lceil an \rceil, \gamma_n} -\left(c^* \beta_n -\alpha^*\right)W_{\lceil an \rceil}\Big|=0,\quad \mbox{in probability}.$$
\end{cor}
\textbf{Proof:} This is a direct consequence of Proposition \ref{prop2}  and Lemma \ref{lemma11}.
\hfill$\Box$

\subsection{Limit behaviour  of  $F_{bad}^{\lceil an \rceil, \gamma_n}$.}\label{Proof-Prop1}

\begin{prop}\label{prop1}
	For any $a\geq 1$, under $\mathbb{P}$,
	\[
	 \lim_{n\to\infty} \sqrt{n} F_{bad}^{\lceil an \rceil, \gamma_n} = 0,\quad \mbox{in probability}.
	\]
\end{prop}
\textbf{Proof:}
Set $Y_n:= \sqrt{n}\left( F_{bad}^{\lceil an \rceil, \gamma_n} \land 1\right)>0$.
We only need to prove that $\lim_{n\to\infty}Y_n=0$ in probability. We claim that \eqref{Constant_c} implies that
\begin{align}\label{step_9}
	\mathbb{E}\left(D_\infty \land y\right)\lesssim\left(1+ \log_+ y\right), \quad y\ge 0.
\end{align}
Indeed, \eqref{Constant_c} implies that $\mathbb{E}\left(D_\infty 1_{\left\{D_\infty \leq y\right\}} \right)- \log_+ y \lesssim 1,y\geq 0$.
\eqref{Constant_c} also implies (see \cite[Theorem 2.2]{BIM})
that $\lim_{y\to+\infty} y\P\left(D_\infty > y\right)=1$, which means $y\P\left(D_\infty > y\right)\lesssim 1, y\geq 0$. Therefore, for all $y\geq 0$,
\[
\mathbb{E}\left(D_\infty \land y\right) = \mathbb{E}\left(D_\infty 1_{\left\{D_\infty \leq y\right\}} \right)+ y \mathbb{P}\left(D_\infty >y\right)\lesssim 1+ \log_+ y.
\]
Since  $V(x)<\gamma_n$ for all $x\in \mathcal{L}^{\lceil an\rceil, \gamma_n}$, we have
\begin{align*}
	F_{bad}^{\lceil an \rceil, \gamma_n}\land 1 &
	 \leq
	  \sum_{x\in \mathcal{L}_{bad}^{\lceil an\rceil, \gamma_n}}
	 c^* e^{-V(x)} \left(D_\infty(x)\land \left(\left(c^*\right)^{-1} e^{V(x)}\right)\right)\\
	 & \leq \sum_{x\in \mathcal{L}_{bad}^{\lceil an\rceil, \gamma_n}}
	c^* e^{-V(x)} \left(D_\infty(x)\land \left(\left(c^*\right)^{-1} e^{\gamma_n}\right)\right).
\end{align*}
It follows from \eqref{step_9} that
\begin{align}\label{step_10}
	&\mathbb{E}\left(F_{bad}^{\lceil an \rceil, \gamma_n}\land 1 \Big| \mathcal{F}_{\mathcal{L}^{\lceil an\rceil, \gamma_n}} \right) \leq \sum_{x\in \mathcal{L}_{bad}^{\lceil an\rceil, \gamma_n}} c^* e^{-V(x)} \mathbb{E}\left(D_\infty(x)\land \left(\left(c^*\right)^{-1} e^{\gamma_n}\right)\Big| \mathcal{F}_{\mathcal{L}^{\lceil an\rceil, \gamma_n}} \right)\nonumber \\ & \lesssim  \gamma_n \sum_{x\in \mathcal{L}_{bad}^{\lceil an\rceil, \gamma_n}} e^{-V(x)}.
\end{align}
By the branching property and the definition of $\mathcal{L}_{bad}^{\lceil an\rceil, \gamma_n}$, we have
\begin{align*}
	&\mathbb{E}\bigg(\sum_{x\in \mathcal{L}_{bad}^{\lceil an\rceil, \gamma_n}} e^{-V(x)}\Big| \mathcal{F}_{\lceil an \rceil}\bigg) =
	\sum_{x\in \mathcal{N}(\lceil an \rceil)}
	e^{-V(x)}1_{\left\{\min_{n\leq j\leq \lceil an\rceil} V(x_j) < \gamma_n \right\}} \\
	& \quad\quad +	\mathbb{E}\bigg(\sum_{k=\lceil an \rceil +1}^\infty
	\sum_{x\in \mathcal{N}(k)}
	e^{-V(x)}1_{\left\{\min_{n\leq j \leq k-1} V(x_j)\geq \gamma_n \right\}}1_{\left\{V(x)<\gamma_n -\beta_n/2\right\}}\Big| \mathcal{F}_{\lceil an \rceil}\bigg)\\
	& \leq \left(W_{\lceil an \rceil} - \widetilde{W}_{\lceil an \rceil}^{n, \gamma_n}\right)+\sum_{u\in \mathcal{N}(\lceil an \rceil)}	e^{-V(u)}1_{\{V(u)\geq \gamma_n\}}\\
	&\quad\quad\times
\mathbb{E}\bigg(\sum_{k= 1}^\infty 	\sum_{x\in \mathcal{N}(k)}	e^{-V(x)}1_{\left\{\min_{0\leq j \leq k-1} V(x_j)\geq y \right\}}1_{\left\{V(x)<y -\beta_n/2\right\}} \bigg)\bigg|_{y = V(u)-\gamma_n}\\
	&= \left(W_{\lceil an \rceil} - \widetilde{W}_{\lceil an \rceil}^{n, \gamma_n}\right)+
	\sum_{u\in \mathcal{N}(\lceil an \rceil)}
	e^{-V(u)}1_{\{V(u)\geq \gamma_n\}}\sum_{k= 1}^\infty
		\mathbf{P}
	\Big(\min_{0\leq j \leq k-1} S_j\geq y , S_k<y -\beta_n/2 \Big)\bigg|_{y = V(u)-\gamma_n},
\end{align*}
where $\widetilde{W}_{\lceil an \rceil}^{n, \gamma_n}$ is defined in \eqref{e:rs3}.
By \eqref{step_17} and \eqref{step_10}, we have
\begin{align*}
\mathbb{E}\left(Y_n \Big| \mathcal{F}_{\lceil an \rceil} \right)
\lesssim \gamma_n \sqrt{n}\left(W_{\lceil an \rceil} - \widetilde{W}_{\lceil an \rceil}^{n, \gamma_n}\right) +\frac{\gamma_n}{\beta_n}\sqrt{n}W_{\lceil an \rceil}
\mathbf{E}
\left(\left(-S_1\right)^3 1_{\left\{S_1 < -\beta_n/2\right\}}\right).
\end{align*}
By Lemma \ref{lemma6}, the first term on the right hand side above converges to $0$ in probability.
Note that  \eqref{step_4} implies that $\gamma_n = O(\beta_n)$. By \eqref{Seneta-Heyde-scaling},
 the second term on the right hand side above also converges to $0$ in probability.
 Then we have shown
$\mathbb{E}\left(Y_n \big| \mathcal{F}_{\lceil an \rceil} \right)$ converges to $0$ in probability.
Finally, applying  Jensen's inequality,  we have
\[
\mathbb{E}\left(e^{-Y_n} \right)= \mathbb{E} \left(\mathbb{E}\left(e^{-Y_n} \big| \mathcal{F}_{\lceil an \rceil} \right)\right)\geq \mathbb{E} \left(\exp\left\{- \mathbb{E}\left(Y_n \big| \mathcal{F}_{\lceil an \rceil} \right) \right\}\right).
\]
By the bounded convergence theorem, we get that
\[
1\geq \limsup_{n\to\infty}\mathbb{E}\left(e^{-Y_n} \right) \geq  \liminf_{n\to\infty}\mathbb{E}\left(e^{-Y_n} \right) \geq \lim_{n\to \infty}\mathbb{E} \left(\exp\left\{- \mathbb{E}\left(Y_n \big| \mathcal{F}_{\lceil an \rceil} \right) \right\}\right) = 1.
\]
Therefore, for any $\varepsilon>0$,
\[
\limsup_{n\to\infty} \mathbb{P}\left(Y_n >\varepsilon\right)=\limsup_{n\to\infty} \mathbb{P}\left(1- e^{-Y_n} >1- e^{-\varepsilon}\right)\leq \frac{1}{1-e^{-\varepsilon}}\mathbb{E}\left(1- e^{-Y_n} \right)=0,
\]
which implies that $\lim_{n\to\infty}Y_n=0$ in probability.
\hfill$\Box$

\subsection{Convergence in distribution for $ \sqrt{n}\left(F_{good}^{\lceil an \rceil,\gamma_n} - \widehat{N}_{good}^{\lceil an \rceil, \gamma_n}\right) $}\label{subsec4.3}

Recall that $F_{good}^{\lceil an \rceil,\gamma_n} $ and $\widehat{N}_{good}^{\lceil an \rceil, \gamma_n}$ are defined in \eqref{Decomposition-2} and \eqref{Widetilde-N} respectively.

\begin{prop}\label{prop4}
	Let $\left(X_t\right)_{t\geq 0}$ be the 1-stable L\'{e}vy process with characteristic function given by \eqref{Charact-Function}. For any $m\geq 1, 1\leq a_1 < ...<a_m$ and $\lambda\in \R^m$, under $\P$,
	\begin{align*}
	&\lim_{m\to\infty} \mathbb{E}\bigg(\exp\bigg\{i \sum_{k=1}^m \lambda_k \sqrt{n}\left(F_{good}^{\lceil a_kn \rceil,\gamma_n} - \widehat{N}_{good}^{\lceil a_kn \rceil, \gamma_n}\right) \bigg\}\bigg| \mathcal{F}_n\bigg)\nonumber\\
	& = \mathbb{E}\bigg(\exp\bigg\{i\sum_{k=1}^m c^*\lambda_k X_{a_k^{-1/2}D_\infty}\bigg\}\bigg| D_\infty \bigg),\quad \qquad \qquad  \mbox{in probability}.
	\end{align*}
\end{prop}
\textbf{Proof: } Let $\lambda_k^* := \lambda_1+...+\lambda_k$ and $a_{m+1}=\infty$,  by the definitions of $F_{good}^{\lceil a_kn \rceil,\gamma_n}$ and  $ \widehat{N}_{good}^{\lceil a_kn \rceil, \gamma_n}$ in \eqref{Decomposition-2} and \eqref{e:rs5},  we have
\begin{align*}
	&\sum_{k=1}^m \lambda_k \sqrt{n}\Big(F_{good}^{\lceil a_kn \rceil,\gamma_n} - \widehat{N}_{good}^{\lceil a_kn \rceil, \gamma_n}\Big) = \sum_{k=1}^m \lambda_k \sqrt{n}\bigg(\sum_{\ell = \lceil a_k n\rceil +1}^\infty 	\sum_{x\in\mathcal{N}(\ell)}	e^{-V(x)}\\
	&\quad\quad\quad\quad\times 1_{\left\{\min_{n\leq j\leq \ell-1} V(x_j)\geq \gamma_n, \gamma_n -\beta_n/2 \leq V(x)< \gamma_n \right\}} c^* \left(D_\infty (x) -V(x)+\frac{\log n}{2} \right)  \bigg) \\
	& = \sum_{k=1}^m  \sum_{\ell = \lceil a_k n\rceil +1}^{\lceil a_{k+1} n\rceil} 	\sum_{x\in\mathcal{N}(\ell)}
	1_{\left\{\min_{n\leq j\leq \ell-1} V(x_j)\geq \gamma_n, \gamma_n -\beta_n/2 \leq V(x)< \gamma_n \right\}} \\&\qquad\qquad \times c^*\lambda_k^*\sqrt{n}e^{-V(x)}\Big( D_\infty (x) - V(x)+\frac{\log n}{2}\Big)  .
\end{align*}
Define $\Psi_{D_\infty}(\lambda):= \mathbb{E}(e^{\mathrm{i}\lambda D_\infty})$, we obtain that
\begin{align*}
	&\mathbb{E}\bigg(\exp\bigg\{\mathrm{i} \sum_{k=1}^m \lambda_k \sqrt{n}\left(F_{good}^{\lceil a_kn \rceil,\gamma_n} -\widehat{N}_{good}^{\lceil a_kn \rceil, \gamma_n}\right) \bigg\}\bigg| \mathcal{F}_{\mathcal{L}^{n,\gamma_n}}\bigg)=\prod_{k=1}^{m}\prod_{\ell = \lceil a_kn \rceil +1}^{\lceil a_{k+1} n\rceil } \\
	 &\quad\quad\times \prod_{x\in \mathcal{N}(\ell), \min_{n\leq j\leq \ell-1} V(x_j)\geq \gamma_n, \gamma_n -\beta_n/2 \leq V(x)< \gamma_n }  \left(e^{-\mathrm{i}c^*\lambda_k^*\sqrt{n}e^{-V(x)}\left(V(x)- \log n/2\right) }\Psi_{D_\infty} \left(c^* \lambda_k^* \sqrt{n}e^{-V(x)}\right)\right),
\end{align*}
where $\mathcal{F}_{\mathcal{L}^{n,\gamma_n}}$ is defined at the beginning of
Section \ref{main} with $a=1$.
 By \cite[(1.12)]{MP} and \cite[Lemma 2.3]{MP}, there exist continuous functions $\eta, q: \mathbb{R}\to \mathbb{C}$ with $\eta(0)=q(0)=0$ such that
\begin{align*}
	& e^{-\mathrm{i}\lambda_k^*c^*\sqrt{n}e^{-V(x)}\left(V(x)-\log n/2\right) }\Psi_{D_\infty} \left(c^* \lambda_k^* \sqrt{n}e^{-V(x)}\right)\\ &= e^{-\mathrm{i}\lambda_k^*c^*\sqrt{n}e^{-V(x)}\left(V(x)-\log n/2\right) }
	 \psi_{\pi/2,c_0+1-\gamma}
	\left(c^* \lambda_k^* \sqrt{n}e^{-V(x)}\right) e^{c^* \lambda_k^* \sqrt{n}e^{-V(x)}\eta\left(c^* \lambda_k^* \sqrt{n}e^{-V(x)}\right) }\\& =: \exp\Big\{-\sqrt{n}e^{-V(x)}
	\psi_{\pi/2, c_0+1-\gamma}
	\left(c^*\lambda_k^*\right) + c^* \lambda_k^* \sqrt{n}e^{-V(x)}q\left(c^* \lambda_k^* \sqrt{n}e^{-V(x)} \right)\Big\},
\end{align*}
where $c_0$  is the constant in \eqref{Constant_c}. Define
\begin{align*}
	R_n&:=-\sum_{k=1}^m \sqrt{n}\left(\widetilde{N}_{good}^{\lceil a_kn \rceil, \gamma_n}-\widetilde{N}_{good}^{\lceil a_{k+1}n \rceil, \gamma_n} \right)
	\psi_{\pi/2, c_0+1-\gamma}(c^*\lambda_k^*),\\
	Z_n&:= \sum_{k=1}^m \sum_{\ell=\lceil a_k n \rceil +1}^{\lceil a_{k+1}n \rceil}\sum_{x\in \mathcal{N}(\ell)} 1_{\left\{\min_{n\leq j\leq \ell-1} V(x_j)\geq \gamma_n, \gamma_n -\beta_n/2 \leq V(x)< \gamma_n \right\}} c^*\lambda_k^* \sqrt{n}e^{-V(x)}q \left(c^* \lambda_k^* \sqrt{n}e^{-V(x)}\right).
\end{align*}
Since $\mathcal{F}_n \subset \mathcal{F}_{\mathcal{L}^{n,\gamma_n}}$,  by the argument above, we have
\begin{align}
	&\mathbb{E}\bigg(\exp\bigg\{\mathrm{i} \sum_{k=1}^m \lambda_k \sqrt{n}\left(F_{good}^{\lceil a_kn \rceil,\gamma_n} - \widehat{N}_{good}^{\lceil a_kn \rceil, \gamma_n}\right) \bigg\}\bigg| \mathcal{F}_{n}\bigg)=\mathbb{E}\left(\exp\left\{R_n + Z_n \right\}\Big| \mathcal{F}_n\right).
\end{align}

We claim that
\begin{align}\label{step_46}
	\Big|\mathbb{E}\left(\exp\left\{R_n+ Z_n \right\}\Big| \mathcal{F}_n\right) - \mathbb{E}\left(\exp\left\{R_n \right\}\Big| \mathcal{F}_n\right) \Big|\to 0\quad\mbox{in probability}.
\end{align}
Indeed, for any $\varepsilon\in (0,1)$
and any complex number $z=a+\mathrm{i}b$ with $|z|\leq \varepsilon$, we have $|a|, |b|\leq \varepsilon$, which implies that
 $$
 \big|e^{z}-1 \big|\leq \big|e^{a+\mathrm{i}b}- e^{\mathrm{i} b} \big|+\big|e^{\mathrm{i} b}-1\big|
 \lesssim \varepsilon.
 $$
 Thus, we get that, for any $\varepsilon>0$,
 \begin{align}\label{step_45}
	&\Big|\mathbb{E}\left(\exp\left\{R_n+ Z_n \right\}\Big| \mathcal{F}_n\right) - \mathbb{E}\left(\exp\left\{R_n \right\}\Big| \mathcal{F}_n\right) \Big| \lesssim \mathbb{P}\left(\left|Z_n \right|> \varepsilon \Big| \mathcal{F}_n\right) + \varepsilon.
\end{align}
Recalling that $\widetilde{N}_{good}^{\lceil a_1 n \rceil, \gamma_n}$ is defined in \eqref{Widetilde-N} with $a=a_1$, we see that, for any fixed $\lambda_1,\cdots, \lambda_m$, we have
\begin{align*}
	 \left|Z_n \right|&\lesssim \max_{1\leq k_0 \leq m}\max_{y\in [\beta_n/2, \beta_n)} \left|q\left(c^* \lambda_{k_0}^* e^{-y}\right)\right| \sum_{k=1}^m \sum_{\ell=\lceil a_k n \rceil +1}^{\lceil a_{k+1}n \rceil}	\sum_{x\in \mathcal{N}(\ell)} \\
	 &\quad\quad\quad 1_{\left\{\min_{n\leq j\leq \ell-1} V(x_j)\geq \gamma_n, \gamma_n -\beta_n/2 \leq V(x)< \gamma_n \right\}}	\sqrt{n}e^{-V(x)}\\
&\lesssim
\max_{1\leq k_0 \leq m}\max_{y\in [\beta_n/2, \beta_n)} \left|q\left(c^* \lambda_{k_0}^* e^{-y}\right)\right| \sqrt{n}\widetilde{N}_{good}^{\lceil a_1 n \rceil, \gamma_n}.
\end{align*}
Applying Corollary \ref{cor1} with $m=1$ and
$z^n=\mathrm{i}a$ with $a\in \R$, we see that $\left(\sqrt{n}\widetilde{N}_{good}^{\lceil a_1 n \rceil, \gamma_n}, \mathbb{P}(\cdot| \mathcal{F}_n)\right)$ converges in distribution.
Since $\lim_{z\to 0}q(z)=0$, we have $\lim_{n\to\infty}\mathbb{P}\left(\left|Z_n \right|> \varepsilon \Big| \mathcal{F}_n\right)=0$.
Letting $n\to\infty$ first and then $\varepsilon \to 0$ in \eqref{step_45}, we get \eqref{step_46}.

Combining  \eqref{step_46} and Corollary \ref{cor1}, we get
\begin{align}\label{step_19}
	&	\lim_{n\to\infty}  \mathbb{E} \bigg(\exp\bigg\{\mathrm{i} \sum_{k=1}^m \lambda_k \sqrt{n}\left(F_{good}^{\lceil a_kn \rceil,\gamma_n} - \widehat{N}_{good}^{\lceil a_kn \rceil, \gamma_n}\right) \bigg\}\bigg| \mathcal{F}_{n}\bigg) \nonumber\\
	&  =\exp\bigg\{\sqrt{\frac{2}{\pi \sigma^2}}D_\infty \sum_{k=1}^m  \psi_{\pi/2, c_0+1-\gamma}
	(c^*\lambda_k^*) \bigg(\frac{1}{\sqrt{a_k}}- \frac{1}{\sqrt{a_{k+1}}}\bigg) \bigg\},\quad \mbox{in probability}.
\end{align}
A standard computation yields
\begin{align*}
		&\mathbb{E}\bigg(\exp\bigg\{\mathrm{i}\sum_{k=1}^m c^*\lambda_k X_{a_k^{-1/2}D_\infty}\bigg\}\bigg| D_\infty \bigg)=\prod_{k=1}^m \mathbb{E}\bigg(\exp\bigg\{\mathrm{i}c^*\lambda_k^*\Big( X_{a_{k+1}^{-1/2}D_\infty}-X_{a_k^{-1/2}D_\infty}\Big)\bigg\}\bigg| D_\infty \bigg) \\ & =	\exp\bigg\{\sum_{k=1}^m D_\infty\bigg(\frac{1}{\sqrt{a_k}}-\frac{1}{\sqrt{a_{k+1}}}\bigg)	\psi_{\sqrt{\pi/2\sigma^2}, (c_0+1-\gamma)\sqrt{2/\pi \sigma^2}}	(c^*\lambda_k^*) \bigg\}.
\end{align*}
This implies the desired result since
$$
\psi_{\sqrt{\pi/2\sigma^2}, (c_0+1-\gamma)\sqrt{2/\pi \sigma^2}}(\lambda) =\sqrt{\frac{2}{\pi \sigma^2}} \psi_{\pi/2, c_0+1-\gamma}(\lambda), \quad \lambda \in \R.
$$
\hfill$\Box$

\textbf{Proof of Theorem \ref{thm1}}: By \eqref{Decomposition-2}, we get that
 \begin{align*} 	
\sqrt{n}\left(D_\infty - D_{\lceil an \rceil} + \frac{\log n}{2}W_{\lceil an \rceil}\right)
  =& \frac{\sqrt{n}}{c^*}\left(D_\infty^{\lceil an\rceil, \gamma_n} -\left(c^* D_{\lceil an \rceil} -\left(c^* \gamma_n -\alpha^*\right)W_{\lceil an \rceil}\right) \right)\\
 &+ \frac{\sqrt{n}}{c^*}\left(\widehat{N}_{good}^{\lceil an \rceil, \gamma_n}-\left(c^*\beta_n-\alpha^*\right) W_{\lceil a_k n\rceil }\right) + \frac{\sqrt{n}}{c^*} F_{bad}^{\lceil an \rceil, \gamma_n}\\
  &+\frac{\sqrt{n}}{c^*}\left(F_{good}^{\lceil an \rceil, \gamma_n} -\widehat{N}_{good}^{\lceil an \rceil, \gamma_n}\right). \end{align*}
It follows from Proposition \ref{prop3},  Corollary \ref{cor2} and Proposition \ref{prop1} that the first three terms
on the right hand side of the display above converge to $0$ in probability as $n\to\infty$.
Now the conclusion of Theorem \ref{thm1} follows from Proposition \ref{prop4}.
 \hfill$\Box$

\section{Proof of Proposition \ref{prop5}}

Recall that $\alpha\in (0, 1]$ is the constant in {\bf(A5)}. Fix an $r\in (0,1)$ such that
\begin{align}\label{step_32}
	\frac{1+\alpha}{2}< \frac{r}{2}(1-\alpha)+ \alpha \left((1+r)\land \left(\frac{3}{2}\right) \right),
\end{align}
which is equivalent to
\[
r> \max\left\{ \frac{1-\alpha}{1+\alpha}, \ \frac{1-2\alpha}{1-\alpha} \right\}.
\]
By the definition \eqref{e:rs3} of $\widetilde{W}_{m}^{n, \gamma_n}$, we have
\begin{align}\label{e:rs6}
	\widetilde{W}_{n}^{ \lceil n^{r}\rceil ,0}=
		\sum_{x\in\mathcal{N}(n)} e^{-V(x)}1_{\left\{ \min_{j\in [n^{r},n]\cap \mathbb{Z}} V(x_j)\geq  0 \right\}}.
\end{align}
Note that $\min_{x\in \mathcal{N}(n)} V(x)\to +\infty, \mathbb{P}$-a.s., so
\begin{align}\label{step_47}
	\lim_{n\to\infty} \mathbb{P}\left(
	\widetilde{W}_{n}^{ \lceil n^{r}\rceil ,0}
	\neq W_n \right) = 0.
\end{align}

\textbf{Proof of Proposition \ref{prop5}:}
Recall the definition \eqref{Delta-n} of $\delta_n$. Let $r\in (0, 1)$ satisfy \eqref{step_32}.
Let $\beta_+\in (0, \frac12((1-r)\wedge r))$ be sufficient small such that
\begin{align}\label{Beta}
	& \Big(\beta_+  +\frac{1}{2}\Big)(1+\alpha)
	< \min\Big\{ (1+r)\land \Big(\frac{3}{2}\Big), \ \frac{r}{2}(1-\alpha) + \alpha (1+r)\land \Big(\frac{3}{2}\Big)  \Big\}.
\end{align}
If we can show that for any $\beta \in (0,\beta_+)$,
\begin{equation}\label{e:rs7}
\lim_{n\to\infty}\mathbb{P}\bigg(n^\beta \bigg|\sqrt{n}W_n -\sqrt{\frac{2}{\pi \sigma^2}}\delta_n D_\infty \bigg|\geq 5 \bigg)=0,
\end{equation}
then the desired conclusion follows immediately.  Note that
\begin{align*}
&\mathbb{P}\bigg(n^\beta \bigg|\sqrt{n}W_n -\sqrt{\frac{2}{\pi \sigma^2}}\delta_n D_\infty \bigg|\geq 5 \bigg)\\
	& \leq \mathbb{P}\Big(
	\widetilde{W}_{n}^{ \lceil n^{r}\rceil ,0}
	 \neq W_n \Big)  + \mathbb{P}\Big(n^{\beta +\frac{1}{2}}\Big|
	\widetilde{W}_{n}^{ \lceil n^{r}\rceil ,0}
	 - \mathbb{E}\big(
	 \widetilde{W}_{n}^{ \lceil n^{r}\rceil ,0}
	  \big| \mathcal{F}_{\lceil n^r \rceil}\big) \Big|> 3 \Big)\\
	& \quad +  \mathbb{P}\bigg(\big(\sup_{k}\delta_k \big) \sqrt{\frac{2}{\pi \sigma^2}}|D_\infty - D_{\lceil n^r \rceil}| >  n^{-\beta} \bigg)+ \mathbb{P}\bigg(n^\beta \bigg|\sqrt{n}\mathbb{E}\left(
	\widetilde{W}_{n}^{ \lceil n^{r}\rceil ,0}
	 \Big| \mathcal{F}_{\lceil n^r \rceil}\right) - \sqrt{\frac{2}{\pi \sigma^2}}\delta_n D_{\lceil n^r \rceil}  \bigg|>1  \bigg).
\end{align*}
Using  \eqref{step_47}, \eqref{step_31}, we know that the first and  third  on the right hand side above tend to 0 as $n\to\infty$.
So it only remains to show that
\begin{align}\label{forth-term}
\lim_{n\to\infty}\mathbb{P}\bigg(n^\beta \bigg|\sqrt{n}\mathbb{E}\left(
\widetilde{W}_{n}^{ \lceil n^{r}\rceil ,0}
\Big| \mathcal{F}_{\lceil n^r \rceil}\right) - \sqrt{\frac{2}{\pi \sigma^2}}\delta_n D_{\lceil n^r \rceil}  \bigg|>1  \bigg)=0,
\end{align}
and
\begin{align}\label{step_48}
	\lim_{n\to \infty} \mathbb{P}\left(\left|
	\widetilde{W}_{n}^{ \lceil n^{r}\rceil ,0}
	- \mathbb{E}\left(
	\widetilde{W}_{n}^{ \lceil n^{r}\rceil ,0}
	 \Big| \mathcal{F}_{\lceil n^r \rceil}\right) \right|> 3 n^{-\beta-\frac{1}{2}}  \right)=0.
\end{align}
We will prove \eqref{forth-term}  and \eqref{step_48} in the following Lemma \ref{lemma13} and  Lemma \ref{proof-step48}, respectively.  Then the proof of Proposition \ref{prop5} will be complete.
\hfill$\Box$

\bigskip

\begin{lemma} \label{lemma13}
Let $r\in (0,1)$ be a fixed constant satisfying  \eqref{step_32}. For any  $\beta \in \left((0,  \frac{1}{2}\left(r\land(1-r)\right)\right)$, it holds that
	\[
	\lim_{n\to \infty}n^\beta \bigg|\sqrt{n}\mathbb{E}\left(
	\widetilde{W}_{n}^{ \lceil n^{r}\rceil ,0}
	 \Big| \mathcal{F}_{\lceil n^r \rceil}\right) - \sqrt{\frac{2}{\pi \sigma^2}}\delta_n D_{\lceil n^r \rceil}  \bigg|=0,\quad \mbox{ in probability}.
	\]
\end{lemma}
\textbf{Proof: } By the branching property, we have
\begin{align*}
	&\sqrt{n}\mathbb{E}\left(
	\widetilde{W}_{n}^{ \lceil n^{r}\rceil ,0}
	\Big| \mathcal{F}_{\lceil n^r \rceil}\right) - \sqrt{\frac{2}{\pi \sigma^2}}\delta_n D_{\lceil n^r \rceil} \\
	&= 	\sum_{x\in\mathcal{N}(\lceil n^r \rceil)}	e^{-V(x)}1_{\left\{V(x)\geq 0 \right\}}\bigg(\sqrt{n}
	 \mathbf{P}_{V(x)}
	\Big(\min_{j\leq n-\lceil n^r \rceil} S_j \geq 0\Big)-\sqrt{\frac{2}{\pi \sigma^2}}\delta_n V(x)\bigg).
\end{align*}
It follows from \cite[Lemma 2.2]{Hu} that,  for any $y\geq 0$ and positive integer $n$ with $n-\lceil n^r \rceil>1$,
\[
\Big|\mathbf{P}_y\Big(\min_{j\leq n-\lceil n^r \rceil} S_j \geq 0\Big) - R(y)\mathbf{P}\Big(\min_{j\leq n-\lceil n^r \rceil} S_j \geq 0\Big) \Big| \lesssim \frac{1+y}{\sqrt{n}}R(y)\mathbf{P}\Big(\min_{j\leq n-\lceil n^r \rceil} S_j \geq 0\Big).
\]
Combining Lemma \ref{lemma1} (i) with the facts that $n-\lceil n^r \rceil\asymp n$ and
$\left|R(y)-c^*y\right| \lesssim 1, y\geq 0$,  we get that for all $y\geq 0$ and large $n$,
\begin{align*}
		&\Big|\sqrt{n}\mathbf{P}_y\Big(\min_{j\leq n-\lceil n^r \rceil} S_j \geq 0\Big)- c^*y\sqrt{n}\mathbf{P}\Big(\min_{j\leq n-\lceil n^r \rceil} S_j \geq 0\Big)\Big|\\	&\lesssim \sqrt{n}\frac{\left(1+y\right)^2}{\sqrt{n-\lceil n^r \rceil}}\mathbf{P}\big(\min_{j\leq n-\lceil n^r \rceil} S_j \geq 0\big) + \sqrt{n}\mathbf{P}\big(\min_{j\leq n-\lceil n^r \rceil} S_j \geq 0\big)\\	& \lesssim\frac{(1+y)^2}{\sqrt{n}}+1.
\end{align*}
Recall that  $\tau_0^-:= \inf\left\{\ell \geq 0: S_\ell <0 \right\}$.  Then by \cite[(2.18)]{Hu},
\begin{align*}
	 	 &\Big|\mathbf{P}\Big(\min_{j\leq n-\lceil n^r \rceil} S_j \geq 0\Big) - \mathbf{P}\Big(\min_{j\leq n} S_j \geq 0\Big) \Big|= \mathbf{P}\Big(\min_{j\leq n-\lceil n^r \rceil} S_j \geq 0, \min_{n-\lceil n^r \rceil+1\leq j \leq n}S_j<0  \Big)\\
	 & = \sum_{\ell= n-\lceil n^r \rceil+1}^n \mathbf{P}\left(\tau_0^- = \ell\right)
	 \lesssim \sum_{\ell= n-\lceil n^r \rceil+1}^n \frac{1}{\ell^{3/2}}\leq \int_{n-\lceil n^r \rceil}^n \frac{1}{x^{3/2}}\mathrm{d}x\\
	& = \Big(\frac{1}{\sqrt{n-\lceil n^r \rceil}}-\frac{1}{\sqrt{n}} \Big) \lesssim \frac{n^r}{n^{3/2}}.
\end{align*}
Combining the two displays above with \eqref{Delta-n}, we get that for all $y\geq 0$ and large $n$,
\begin{align*}
		&\Big|\sqrt{n}\mathbf{P}_y	\Big(\min_{j\leq n-\lceil n^r \rceil} S_j \geq 0\Big)- \sqrt{\frac{2}{\pi \sigma^2}}\delta_n y\Big| \lesssim \frac{(1+y)^2}{\sqrt{n}}+ 1+ \frac{y}{n^{1-r}},
\end{align*}
which implies that
\begin{align*}
	& \Big|\sqrt{n}\mathbb{E}\Big(
	\widetilde{W}_{n}^{ \lceil n^{r}\rceil ,0}
	\Big| \mathcal{F}_{\lceil n^r \rceil}\Big) - \sqrt{\frac{2}{\pi \sigma^2}}\delta_n D_{\lceil n^r \rceil}\Big|\\
	& \lesssim 	\sum_{x\in\mathcal{N}(\lceil n^r \rceil)}
	e^{-V(x)}1_{\left\{V(x)\geq 0 \right\}}\Big(\frac{\left(1+V(x)\right)^2}{\sqrt{n}}+ 1 + \frac{V(x)}{n^{1-r}} \Big).
\end{align*}
By \eqref{Seneta-Heyde-scaling}, we have
\begin{align}\label{e:rs17}
	&\lim_{n\to\infty}\sqrt{\lceil n^r \rceil}	\sum_{x\in\mathcal{N}(\lceil n^r \rceil)}
	e^{-V(x)} \to \sqrt{\frac{2}{\pi \sigma^2}}D_\infty\quad \mbox{in probability}.
\end{align}
Since $\lim_{n\to\infty}D_n=D_\infty$ almost surely and $D_\infty$ is non-negative, we have
\begin{align}\label{e:rs18}
\lim_{n\to\infty} \sum_{x\in\mathcal{N}(\lceil n^r \rceil)} V(x)e^{-V(x)}1_{\{V(x)\geq 0\}} = D_\infty,\ \mbox{almost surely.}
\end{align}
Since $\beta\in (0, \frac{1}{2}((1-r)\wedge r))$,  using the two displays above, we get that for any $\varepsilon>0$,
\[
\lim_{n\to\infty} \mathbb{P}\bigg(n^\beta\bigg|\sum_{x\in\mathcal{N}(\lceil n^r \rceil)}
e^{-V(x)}1_{\left\{V(x)\geq 0 \right\}}\Big( 1 + \frac{V(x)}{n^{1-r}} \Big) \bigg| >\varepsilon \bigg)=0.
\]
For any $\delta>0$, let $L$ satisfy \eqref{Global-Minimum}. Then for any $\varepsilon>0$,
\begin{align*}
	&\mathbb{P}\bigg(n^\beta 	\sum_{x\in\mathcal{N}(\lceil n^r \rceil)}
	e^{-V(x)}1_{\left\{V(x)\geq 0 \right\}}\frac{\left(1+V(x)\right)^2}{\sqrt{n}} > \varepsilon\bigg)\\
	& \leq \delta + \mathbb{P}\bigg(n^\beta 	\sum_{x\in\mathcal{N}(\lceil n^r \rceil)}
	e^{-V(x)}1_{\left\{V(x)\geq 0 \right\}}\frac{\left(1+V(x)\right)^2}{\sqrt{n}}1_{\{\min_{j\leq \lceil n^r \rceil} V(x_j)\geq -L \}} > \varepsilon\bigg)\\
	& \leq \delta + \frac{1}{\varepsilon}n^{\beta - \frac{1}{2}}\mathbb{E}\bigg(
	\sum_{x\in\mathcal{N}(\lceil n^r \rceil)} \left(1+V(x)\right)^2e^{-V(x)}1_{\left\{V(x)\geq 0 \right\}}1_{\{\min_{j\leq \lceil n^r \rceil} V(x_j)\geq -L \}} \bigg)\\
	& = \delta + \frac{1}{\varepsilon}n^{\beta -\frac{1}{2}}
	\mathbf{E}
	\Big((1+S_{\lceil n^r \rceil})^2 1_{\{S_{\lceil n^r \rceil}\geq 0\}}1_{\left\{\min_{j\leq \lceil n^r \rceil} S_j \geq -L \right\}} \Big)\\
	& \lesssim \delta + \frac{1}{\varepsilon}n^{\beta -\frac{1}{2}}\Big( (1+L)^2+ (1+L)\sqrt{\lceil n^r \rceil}\Big),
\end{align*}
where in the last inequality we used Lemma \ref{lemma12}.  Letting $n\to\infty$ first, and then $\delta\to 0$, we get
\[
\lim_{n\to\infty} \mathbb{P}\bigg(n^\beta \sum_{x\in\mathcal{N}(\lceil n^r \rceil)}
e^{-V(x)}1_{\left\{V(x)\geq 0 \right\}}\frac{\left(1+V(x)\right)^2}{\sqrt{n}} > \varepsilon\bigg)=0.
\]
This completes the proof.
\hfill$\Box$

\bigskip

\begin{lemma}\label{proof-step48}
Let $\beta_+\in (0, \frac12((1-r)\wedge r))$ satisfy \eqref{Beta}. For any $\beta\in (0, \beta^+)$,
\begin{align}
	\lim_{n\to \infty} \mathbb{P}\left(\left|
	\widetilde{W}_{n}^{ \lceil n^{r}\rceil ,0}
	- \mathbb{E}\left(
	\widetilde{W}_{n}^{ \lceil n^{r}\rceil ,0}
	\Big| \mathcal{F}_{\lceil n^r \rceil}\right) \right|> 3 n^{-\beta-\frac{1}{2}}  \right)=0.
\end{align}
\end{lemma}
\textbf{Proof: } We only need to prove that
\begin{equation}\label{toprove}
\mathbb{P}\Big(\Big|
\widetilde{W}_{n}^{ \lceil n^{r}\rceil ,0}
- \mathbb{E}\left(
\widetilde{W}_{n}^{ \lceil n^{r}\rceil ,0}
\Big| \mathcal{F}_{\lceil n^r \rceil}\right) \Big|> 3 n^{-\beta-\frac{1}{2}} \big| \mathcal{F}_{\lceil n^r \rceil}\Big)\to 0\ \ \mbox{in probability}.
\end{equation}
For $m<n$, define
\begin{align*}
	&\mathcal{A}_{n,1}^{m, 0} :=
	\bigg\{x\in \mathcal{N}(n): 	\mbox{ for all }j=m+1, \dots, n,\\
	&\qquad \qquad\sum_{u\in \Omega(x_j)} \left(1+\left(V(u)-V(x_{j-1})\right)_+\right)e^{-\left(V(u)-V(x_{j-1})\right)} \leq  e^{V(x_{j-1})/2} \bigg\},\\
	& \widetilde{W}_{n,1}^{\lceil n^r \rceil,0}   =\sum_{x\in \mathcal{N}(n)}
	e^{-V(x)}1_{\left\{ \min_{j\in [n^{r},n]\cap \mathbb{Z}} V(x_j)\geq  0 \right\}}1_{\big\{x\in \mathcal{A}_{n,1}^{\lceil n^r \rceil, 0}\big\}}.
\end{align*}
Recall that $\alpha \in (0,1]$ is given in {\bf(A5)}. Then
\begin{align}\label{term1+term2}
	&\mathbb{P}\Big(\Big|
	\widetilde{W}_{n}^{ \lceil n^{r}\rceil ,0}
	- \mathbb{E}\left(
	\widetilde{W}_{n}^{ \lceil n^{r}\rceil ,0}
	 \Big| \mathcal{F}_{\lceil n^r \rceil}\right) \Big|> 3
	n^{-\beta-\frac{1}{2}} \big| \mathcal{F}_{\lceil n^r \rceil}\Big)\nonumber\\
	& \lesssim 2n^{\beta + \frac{1}{2}}\mathbb{E}\Big(
	\widetilde{W}_{n}^{ \lceil n^{r}\rceil ,0} - \widetilde{W}_{n,1}^{\lceil n^{r}\rceil ,0} \Big| \mathcal{F}_{\lceil n^r \rceil} \Big)\nonumber\\
 &\quad + n^{(\beta + \frac{1}{2})(1+\alpha)}\mathbb{E}\Big(\left|
 \widetilde{W}_{n,1}^{\lceil n^{r}\rceil ,0}-\mathbb{E}\left(\widetilde{W}_{n,1}^{\lceil n^{r}\rceil,0} \Big| \mathcal{F}_{\lceil n^r \rceil} \right)  \right|^{1+\alpha}\Big| \mathcal{F}_{\lceil n^r \rceil} \Big).
\end{align}
We estimate the two term on the right hand side of  the above display by the following two steps respectively.

{\bf Step 1}
 By \eqref{Many-to-one-1} and the branching property, we have
\begin{align}\label{cond-minus}
	&\mathbb{E}\left(\widetilde{W}_{n}^{\lceil n^{r}\rceil ,0}-\widetilde{W}_{n,1}^{\lceil n^{r}\rceil ,0} \Big| \mathcal{F}_{\lceil n^r \rceil} \right)\nonumber\\
	&= \sum_{u\in \mathcal{N}(\lceil n^r \rceil)}
	1_{\left\{V(u)\geq 0 \right\}}
\mathbb{E}_y\bigg( \sum_{x\in \mathcal{N}(n-\lceil n^r \rceil-1)}
	e^{-V(x)}1_{\left\{ \min_{j\leq n-\lceil n^r \rceil} V(x_j)\geq  0 \right\}}1_{\big\{x\notin \mathcal{A}_{n-\lceil n^r \rceil-1,1}^{0, 0}\big\}}\bigg)_{y= V(u)}\nonumber\\
	& \leq  \sum_{u\in \mathcal{N}(\lceil n^r \rceil)}	e^{-V(u)}1_{\{V(u)\geq 0\}} \sum_{\ell=1}^{n-\lceil n^r \rceil-1}  \mathbb{Q}\bigg(
	\min_{j\leq n-\lceil n^r \rceil-1} V(w_j)\geq -y, \nonumber\\	 &\quad\quad \sum_{u\in \Omega(w_\ell)}
	\left(1+\left(V(u)-V(w_{\ell-1})\right)_+\right)e^{-\left(V(u)-V(w_{\ell-1})\right)} >  e^{\left(V(w_{\ell-1})+ y\right)/2}  \bigg)\bigg|_{y= V(u)}.
\end{align}
Suppose $n$ is large so that $n-\lceil n^{r}\rceil >1$. For any positive integer $\ell \leq  n-\lceil n^r \rceil-1$, conditioned on $\mathcal{F}_{\ell}$, by Lemma \ref{lemma1} (i), we have
 \begin{align}
 &\mathbb{Q}\bigg(
	  \min_{j\leq n-\lceil n^r \rceil-1} V(w_j)\geq -y,	 \sum_{u\in \Omega(w_\ell)} \left(1+\left(V(u)-V(w_{\ell-1})\right)_+\right)e^{-\left(V(u)-V(w_{\ell-1})\right)} >  e^{\left(V(w_{\ell-1})+ y\right)/2}  \bigg)\nonumber\\
	& \lesssim \bigg(1\land \frac{1}{\sqrt{n-\lceil n^r \rceil-1-\ell}}\bigg)  \mathbb{E}_{\mathbb{Q}}\bigg(\left(1+V(w_\ell)+y\right)1_{\{\min_{j\leq \ell } V(w_j)\geq -y\}} \nonumber\\
	&\quad\quad\times  1_{\{\sum_{u\in \Omega(w_\ell)\cup\{w_\ell\}} \left(1+\left(V(u)-V(w_{\ell-1})\right)_+\right)e^{-\left(V(u)-V(w_{\ell-1})\right)} >  e^{\left(V(w_{\ell-1})+ y\right)/2} \}} \bigg)\nonumber\\
	& =: \bigg(1\land \frac{1}{\sqrt{n-\lceil n^r \rceil-1-\ell}}\bigg) \E_{\Q}(J_\ell )
	.\label{e:rs8}
 \end{align}
Conditioned on $\mathcal{F}_{\ell -1}$, we get that, given $V(w_{\ell -1})=z$,
\begin{align*}
	& \E_{\Q} (J_\ell\big| V(w_{\ell -1})=z)
	\\&=  \mathbb{E}_{\mathbb{Q}}\left(\left(1+ V(w_1)+z+y\right)1_{\{V(w_1)\geq -y-z\}} 1_{\{
		\sum_{x\in \mathcal{N}(1)}
		\left(1+\left(V(u)\right)_+\right)e^{-V(u)} >  e^{\left(z+ y\right)/2} \}} \right)\\
	& = \mathbb{E}\bigg(
	\sum_{x\in \mathcal{N}(1)}
	\left(1+V(x)+z+y \right)e^{-V(x)}1_{\{V(x)\geq -y-z\}}1_{\left\{z+y < 2\log_+ \left( W_1 + \widetilde{W}_1 \right) \right\}}  \bigg)\\
	& \leq \mathbb{E}\left(\left(W_1\left(1+ 2\log_+ \left( W_1 + \widetilde{W}_1 \right) \right)+ \widetilde{W}_1\right) 1_{\left\{z+y < 2\log_+ \left( W_1 + \widetilde{W}_1 \right) \right\}}   \right).
\end{align*}
Therefore,
\begin{align}
		& \E_{\Q}(J_\ell) \leq \mathbb{E}\bigg(\left(W_1\left(1+ 2\log_+ \left( W_1 + \widetilde{W}_1 \right) \right)+ \widetilde{W}_1\right) \nonumber\\
		& \quad\quad\times
				\mathbf{P}
		\left(\min_{j\leq \ell-1}S_j \geq -y, S_{\ell-1} +y < 2\log_+\left(m_1\right)\right)\bigg|_{m_1 = W_1 +\widetilde{W}_1} \bigg)\nonumber\\
		&=:  \mathbb{E}\bigg( W
					\mathbf{P}
		\bigg(\min_{j\leq \ell-1}S_j \geq -y, S_{\ell-1} +y < 2\log_+\left(m_1\right)\bigg)\bigg|_{m_1 = W_1 +\widetilde{W}_1} \bigg),\label{e:rs9}
\end{align}
where $W:= W_1\big(1+ 2\log_+ \big( W_1 + \widetilde{W}_1 \big) \big)+ \widetilde{W}_1.$
Thus, if $k_n$ is the integer such that $2k_n \leq n-\lceil n^r \rceil-1 < 2k_n +1$, then by Lemma \ref{lemma1} (ii)(iii), we have
\begin{align}
	&\sum_{\ell =1}^{n-\lceil n^r \rceil-1}\bigg(1\land \frac{1}{\sqrt{n-\lceil n^r \rceil-1-\ell}}\bigg)\E_{\Q}(J_\ell) \nonumber\\
	& \leq \frac{1}{\sqrt{n-\lceil n^r \rceil-1-k_n}}\mathbb{E}\bigg(W \sum_{\ell=1}^{k_n}
	\mathbf{P}
	\bigg(\min_{j\leq \ell-1}S_j \geq -y, S_{\ell-1} +y < 2\log_+\left(m_1\right)\bigg)\bigg|_{m_1 = W_1 +\widetilde{W}_1} \bigg)\nonumber\\
	& \quad + \sum_{\ell= k_n +1}^{n-\lceil n^r \rceil-1}\bigg(1\land \frac{1}{\sqrt{n-\lceil n^r \rceil-1-\ell}}\bigg)\mathbb{E}\bigg(W\frac{(y+1)\left(1+2\log_+\left(W_1 +\widetilde{W}_1\right)\right)^2}{(\ell-1)^{3/2}}\bigg)\nonumber\\
	& \lesssim \bigg(\frac{1}{\sqrt{n-\lceil n^r \rceil-1-k_n}} + (y+1)\sum_{\ell= k_n +1}^{n-\lceil n^r \rceil-1}\bigg(1\land \frac{1}{\sqrt{n-\lceil n^r \rceil-1-\ell}}\bigg) \frac{1}{(\ell-1)^{3/2}}  \bigg)F(1), \label{e:rs10}
 \end{align}
with $F(1)$ given by
\[
F(1):= \mathbb{E}\bigg(\Big(W_1\big(1+ 2\log_+ \big( W_1 + \widetilde{W}_1 \big) \big)+ \widetilde{W}_1\Big) \Big(1+2\log_+\big(W_1 +\widetilde{W}_1\big)\Big)^2\bigg)<\infty.
\]
If  $n$ is large enough, we have
\begin{align*}
	&\sum_{\ell= k_n +1}^{n-\lceil n^r \rceil-1}\bigg(1\land \frac{1}{\sqrt{n-\lceil n^r \rceil-1-\ell}}\bigg) \frac{1}{(\ell-1)^{3/2}}\\
	& \leq \frac{1}{(n-\lceil n^r \rceil-1)^{3/2}}+\sum_{\ell= k_n +1}^{n-\lceil n^r \rceil-2} \frac{1}{\sqrt{n-\lceil n^r \rceil-1-\ell}} \frac{1}{(\ell-1)^{3/2}}\\
	&\lesssim \frac{1}{n} +\sum_{\ell= k_n +1}^{n-\lceil n^r \rceil-2}  \int_{\ell}^{\ell+1}\frac{1}{\sqrt{(n-\lceil n^r \rceil-1-x)(x-2)^{3/2}}} \mathrm{d}x\\
	& =  \frac{1}{n} + \frac{1}{n}\int_{n^{-1}(k_n+1)}^{n^{-1}(n-\lceil n^r \rceil-2)}\frac{1}{\sqrt{(n^{-1}(n-\lceil n^r \rceil-1)-x)(x-2n^{-1})^{3/2}}} \mathrm{d}x\lesssim \frac{1}{n}.
\end{align*}
This implies that for $n$ large enough,
\[
\frac{1}{\sqrt{n-\lceil n^r \rceil-1-k_n}} + (y+1)\sum_{\ell= k_n +1}^{n-\lceil n^r \rceil-1}\bigg(1\land \frac{1}{\sqrt{n-\lceil n^r \rceil-1-\ell}}\bigg) \frac{1}{(\ell-1)^{3/2}}  \lesssim \frac{1}{\sqrt{n}}+\frac{y}{n}.
\]
Combining this with \eqref{cond-minus}, \eqref{e:rs8} and \eqref{e:rs10}, we get that for large $n$,
\begin{align}\label{estimate-term1}
	\mathbb{E}\left(\widetilde{W}_{n}^{\lceil n^{r}\rceil ,0}-\widetilde{W}_{n,1}^{\lceil n^{r}\rceil ,0} \Big| \mathcal{F}_{\lceil n^r \rceil} \right) \lesssim
	\sum_{u\in\mathcal{N}(\lceil n^r \rceil)}
	\Big(\frac{1}{\sqrt{n}}+\frac{V(u)}{n}\Big)e^{-V(u)}1_{\{V(u)\geq 0\}}.
\end{align}

{\bf Step 2}  By the branching property,
\begin{align*}
&\mathbb{E}\Big(\left|\widetilde{W}_{n,1}^{\lceil n^{r}\rceil ,0}- \mathbb{E} \left(\widetilde{W}_{n,1}^{\lceil n^{r}\rceil ,0} \Big| \mathcal{F}_{\lceil n^r \rceil} \right)  \right|^{1+\alpha}\Big| \mathcal{F}_{\lceil n^r \rceil} \Big) \\
	&\lesssim 	\sum_{u\in\mathcal{N}(\lceil n^r \rceil)}	1_{\{V(u)\geq 0\}} \mathbb{E}_y\bigg(\bigg(	\sum_{x\in\mathcal{N}(n-\lceil n^r \rceil-1)}
	 e^{-V(x)}1_{\left\{ \min_{j\leq n-\lceil n^r \rceil-1} V(x_j)\geq  0 \right\}}1_{\left\{x\in \mathcal{A}_{n-\lceil n^r \rceil-1,1}^{0, 0}\right\}} \bigg)^{1+\alpha} \bigg)\bigg|_{y= V(u)}
	 \nonumber\\ &=: \sum_{u\in\mathcal{N}(\lceil n^r \rceil)}
	 1_{\{V(u)\geq 0\}} \mathbb{E}_y\left( \Gamma^{1+\alpha}\right)\big|_{y=V(u)}.
\end{align*}
Recalling the definition of $\Q_y$ in beginning of Subsection \ref{ss:spine}, we similarly have
\begin{align*}
	&\mathbb{E}_y\left(\Gamma ^{1+\alpha} \right)=e^{-y}\mathbb{E}_{\mathbb{Q}_y}\Big(1_{\{\min_{j\leq n-\lceil n^r \rceil-1} V(w_j)\geq 0\}}1_{\big\{w\in \mathcal{A}_{n-\lceil n^r \rceil-1,1}^{0, 0} \big\}}\Gamma ^\alpha  \Big).
\end{align*}
For $u<x$,  define $V(x;u):= V(x)-V(u)$. Recall that $\mathcal{N}(u,m):= \left\{x\in \mathcal{N}(|u|+m): x>u\right\}$.  By the spine decomposition, we have
\begin{align*}
	& \Gamma = \sum_{x\in\mathcal{N}(n-\lceil n^r \rceil-1)}
	e^{-V(x)}1_{\left\{ \min_{j\leq n-\lceil n^r \rceil-1} V(x_j)\geq  0 \right\}}1_{\left\{x\in \mathcal{A}_{n-\lceil n^r \rceil-1}^{0, 0}\right\}}\\
	& \leq e^{-V(w_{n-\lceil n^r \rceil-1})}+ \sum_{k=1}^{n-\lceil n^r \rceil-1}\sum_{u\in \Omega(w_k)}e^{-V(u)}
         \sum_{x\in\mathcal{N}(u, n-\lceil n^r \rceil-1-k)}
	e^{-V(x;u)}1_{\{\min_{j\leq n-\lceil n^r \rceil-1-k} V(x_j;u)\geq - V(u) \}}\\
	&=:  e^{-V(w_{n-\lceil n^r \rceil-1})}+ \sum_{k=1}^{n-\lceil n^r \rceil-1}\sum_{u\in \Omega(w_k)}e^{-V(u)}H(u),
\end{align*}
where $H(u):=\sum_{x\in\mathcal{N}(u, n-\lceil n^r \rceil-1-k)}
e^{-V(x;u)}1_{\{\min_{j\leq n-\lceil n^r \rceil-1-k} V(x_j;k)\geq - V(u) \}}$.
Therefore, using the inequality $\E(|X|^\alpha)\leq (\E(|X|))^\alpha$, we get
\begin{align*}
	&\mathbb{E}_{\mathbb{Q}_y}\Big(1_{\{\min_{j\leq n-\lceil n^r \rceil-1} V(w_j)\geq 0\}}1_{\left\{w\in \mathcal{A}_{n-\lceil n^r \rceil-1,1}^{0, 0} \right\}}\Gamma^\alpha \Big)\\
	& \leq \mathbb{E}_{\mathbb{Q}_y}\left(e^{-\alpha V(w_{n-\lceil n^r \rceil-1})}1_{\left\{ \min_{j\leq n-\lceil n^r \rceil-1} V(x_j)\geq  0 \right\}} \right)\\
	& \quad +\bigg( \mathbb{E}_{\mathbb{Q}_y}\bigg(1_{\{\min_{j\leq n-\lceil n^r \rceil-1} V(w_j)\geq 0\}}1_{\left\{w\in \mathcal{A}_{n-\lceil n^r \rceil-1,1}^{0, 0} \right\}} \sum_{k=1}^{n-\lceil n^r \rceil-1}\sum_{u\in \Omega(w_k)}e^{-V(u)}H(u)\bigg)\bigg)^\alpha
\end{align*}
By the branching property and Lemma \ref{lemma1} (i), we have
\begin{align*}
	&\mathbb{E}_{\mathbb{Q}_y}\bigg(1_{\{\min_{j\leq n-\lceil n^r \rceil-1} V(w_j)\geq 0\}}1_{\big\{w\in \mathcal{A}_{n-\lceil n^r \rceil-1,1}^{0, 0} \big\}} \sum_{k=1}^{n-\lceil n^r \rceil-1}\sum_{u\in \Omega(w_k)}e^{-V(u)}H(u)\bigg)\\
	 & = \mathbb{E}_{\mathbb{Q}_y}\bigg(1_{\{\min_{j\leq n-\lceil n^r \rceil-1} V(w_j)\geq 0\}}1_{\big\{w\in \mathcal{A}_{n-\lceil n^r \rceil-1,1}^{0, 0} \big\}}\\
	 &\quad\quad \times \sum_{k=1}^{n-\lceil n^r \rceil-1}\sum_{u\in \Omega(w_k)}e^{-V(u)}
		 \mathbf{P}
	 \Big(\min_{j\leq n-\lceil n^r \rceil-1-k}S_j \geq -z \Big)\Big|_{z=V(u)}\bigg)\\
	& \lesssim \mathbb{E}_{\mathbb{Q}_y}\bigg(1_{\{\min_{j\leq n-\lceil n^r \rceil-1} V(w_j)\geq 0\}}1_{\big\{w\in \mathcal{A}_{n-\lceil n^r \rceil-1,1}^{0, 0} \big\}}\\
	&\quad\quad \times  \sum_{k=1}^{n-\lceil n^r \rceil-1}\sum_{u \in\Omega(w_k)} e^{-V(u)}\bigg(1\land \frac{(1+V(u))}{\sqrt{n-\lceil n^r \rceil-1-k}}\bigg)\bigg).
\end{align*}
Since
\begin{align*}
	&1\land \frac{(1+V(u))}{\sqrt{n-\lceil n^r \rceil-1-k}} \leq 1\land \frac{\left(1+\left(V(u)-V(w_{k-1})\right)_+ + V(w_{k-1})\right)}{\sqrt{n-\lceil n^r \rceil-1-k}}\\
	& \leq 1\land \frac{\left(1+\left(V(u)-V(w_{k-1})\right)_+ \right)\left(1+V(w_{k-1})\right)}{\sqrt{n-\lceil n^r \rceil-1-k}}\\
	& \leq \left(1+\left(V(u)-V(w_{k-1})\right)_+\right) \bigg(1\land \frac{(1+V(w_{k-1}))}{\sqrt{n-\lceil n^r \rceil-1-k}}\bigg),
\end{align*}
and, on the set $\big\{w\in \mathcal{A}_{n-\lceil n^r \rceil-1,1}^{0, 0} \big\}$, it holds that
\[
\sum_{u \in\Omega(w_k)} e^{-V(u)}\left(1+\left(V(u)-V(w_{k-1})\right)_+\right)\leq e^{-V(w_{k-1})/2}.
\]
We have
\begin{align*}
	&\mathbb{E}_{\mathbb{Q}_y}\bigg(1_{\{\min_{j\leq n-\lceil n^r \rceil-1} V(w_j)\geq 0\}}1_{\left\{w\in \mathcal{A}_{n-\lceil n^r \rceil-1,1}^{0, 0} \right\}} \\
	&\quad\quad\times \sum_{k=1}^{n-\lceil n^r \rceil-1}\sum_{u \in\Omega(w_k)} e^{-V(u)}\Big(1\land \frac{(1+V(u))}{\sqrt{n-\lceil n^r \rceil-1-k}}\Big)\bigg)\\
	& \lesssim \mathbb{E}_{\mathbb{Q}_y}\bigg(1_{\{\min_{j\leq n-\lceil n^r \rceil-1} V(w_j)\geq 0\}} \sum_{k=1}^{n-\lceil n^r \rceil-1} e^{-V(w_{k-1})/2} \bigg(1\land \frac{(1+V(w_{k-1}))}{\sqrt{n-\lceil n^r \rceil-1-k}}\bigg) \bigg).
\end{align*}
Therefore, by the arguments above, we conclude that
\begin{align}\label{step_30}
	&\mathbb{E}_y\left(\Gamma ^{1+\alpha} \right)\lesssim e^{-y}
		 \mathbf{E}_y
	\left(e^{-\alpha S_{n-\lceil n^r \rceil-1}} 1_{\{\min_{j\leq n-\lceil n^r \rceil-1} S_j \geq 0\}}\right) \nonumber\\
	& + e^{-y}\bigg( \sum_{k=1}^{n-\lceil n^r \rceil-1}
		 \mathbf{E}_y
	 \bigg(1_{\{\min_{j\leq n-\lceil n^r \rceil-1} S_j\geq 0\}} e^{-S_{k-1}/2} \bigg(1\land \frac{(1+S_{k-1})}{\sqrt{n-\lceil n^r \rceil-1-k}}\bigg) \bigg)	\bigg)^\alpha\nonumber\\
	&:= e^{-y}(I+II^\alpha).
\end{align}
By Lemma \ref{lemma1} (ii), for all $y\geq 0$ and large $n$, it holds that
\begin{align}
	I&
	\leq \frac{1}{n}+
		 \mathbf{P}_y
	\left(\min_{j\leq n-\lceil n^r \rceil-1} S_j \geq 0, S_{n-\lceil n^r \rceil-1}\leq \alpha^{-1}\log n \right)\nonumber\\
	&\lesssim \frac{1}{n}+  \frac{(1+y)(1+\alpha^{-1}\log n)^2}{(n-\lceil n^r \rceil-1)^{3/2}}\lesssim \frac{1}{n}+ y \frac{(\log n)^2}{n^{3/2}}.\label{e:rs12}
\end{align}
For $II$, we have
\begin{align}
	II& \leq \mathbf{E}_y \left(e^{-S_{n-\lceil n^r \rceil-1}/2} 1_{\{\min_{j\leq n-\lceil n^r \rceil-1} S_j \geq 0\}}\right) + \mathbf{E}_y \left(e^{-S_{n-\lceil n^r \rceil-2}/2} 1_{\{\min_{j\leq n-\lceil n^r \rceil-2} S_j \geq 0\}}\right)\nonumber\\
	& \quad +\sum_{k=1}^{n-\lceil n^r \rceil-3} \mathbf{E}_y\bigg(1_{\{\min_{j\leq n-\lceil n^r \rceil-1} S_j\geq 0\}} e^{-S_{k-1}/2}  \frac{(1+S_{k-1})}{\sqrt{n-\lceil n^r \rceil-1-k}}
	\bigg)\nonumber\\
	&:=II_1+II_2+II_3.\label{e:rs13}
\end{align}
Applying Lemma \ref{lemma1} (ii) with $k=n-\lceil n^r \rceil-1$ or $n-\lceil n^r \rceil-2$, we get
\[
\mathbf{E}_y
\left(e^{- S_{k}/2} 1_{\{\min_{j\leq k} S_j \geq 0\}}\right) \lesssim \frac{1}{n}+  \frac{(1+y)(1+2\log n)^2}{k^{3/2}}\lesssim\frac{1}{n}+ y\frac{(\log n)^2}{n^{3/2}}.
\]
Thus
\begin{equation}\label{e:rs14}
II_1\vee II_2\lesssim\frac{1}{n}+ y\frac{(\log n)^2}{n^{3/2}}.
\end{equation}
For $II_3$, we note that for, $k=1, \dots, n-\lceil n^r \rceil-3$,
\begin{align*}
&\mathbf{E}_y\bigg(e^{-S_{k-1}/2}  \frac{(1+S_{k-1})}{\sqrt{n-\lceil n^r \rceil-1-k}}1_{\{\min_{j\leq n-\lceil n^r \rceil-1} S_j\geq 0\}}  \bigg)\\
	& = \mathbf{E}_y\bigg(e^{-S_{k-1}/2}  \frac{(1+S_{k-1})}{\sqrt{n-\lceil n^r \rceil-1-k}}1_{\{\min_{j\leq k-1} S_j\geq 0\}}
	\mathbf{P}_{S_{k-1}}\bigg(\min_{j\leq n-\lceil n^r \rceil -k} S_j \geq 0 \bigg) \bigg)\\
	& \lesssim \mathbf{E}_y\bigg(e^{-S_{k-1}/2}  \frac{(1+S_{k-1})^2}{n-\lceil n^r \rceil-1-k}1_{\{\min_{j\leq k-1} S_j\geq 0\}} \bigg).
\end{align*}
 Recall that $k_n$  is the integer such that $2k_n \leq n-\lceil n^r \rceil-1 < 2k_n +1$. Combining Lemma \ref{lemma1} (iv) with $\lambda = \frac{1}{4}$ and the fact that $\sup_{x>0} \left(e^{-x/2}(1+x)^2\right)<\infty$,
\begin{align}
&\sum_{k=1}^{k_n}\mathbf{E}_y\bigg(e^{-S_{k-1}/2}
	\frac{(1+S_{k-1})}{\sqrt{n-\lceil n^r \rceil-1-k}}1_{\{\min_{j\leq n-\lceil n^r \rceil-1} S_j\geq 0\}}  \bigg)\nonumber\\
	& \lesssim \frac{1}{n-\lceil n^r \rceil-k_n-1}\sum_{k=0}^\infty \mathbf{E}_y \left(e^{-S_k/4}1_{\{\min_{j\leq k} S_j \geq 0\}}\right)\lesssim \frac{1}{n}.\label{e:rs11}
\end{align}
On the other hand, for $k_n+1 \leq k \leq n-\lceil n^r \rceil-3$, by Lemma \ref{lemma1} (ii), we have
\begin{align*}
		&\mathbf{E}_y\Big(e^{-S_{k-1}/2}  \frac{(1+S_{k-1})^2}{n-\lceil n^r \rceil-1-k}1_{\{\min_{j\leq k-1} S_j\geq 0\}} \Big)\\ &=\sum_{\ell=0}^\infty \mathbf{E}_y\Big(e^{-S_{k-1}/2}  \frac{(1+S_{k-1})^2}{n-\lceil n^r \rceil-1-k}1_{\{\min_{j\leq k-1} S_j\geq 0\}}1_{\{S_{k-1}\in [\ell, \ell+1) \}} \Big)\\
	&\leq \frac{1}{n-\lceil n^r \rceil-1-k}\sum_{\ell=0}^\infty e^{-\ell /2} (2+\ell)^2 \mathbf{P}_y\Big(\min_{j\leq k-1} S_j\geq 0, S_{k-1}\in [\ell, \ell+1)  \Big)\\
	& \lesssim  \frac{1}{n-\lceil n^r \rceil-1-k}\frac{(y+1)}{(k-1)^{3/2}}\sum_{\ell=0}^\infty e^{-\ell /2} (2+\ell)^2 (\ell +2)\lesssim  \frac{(1+y)}{\left(n-\lceil n^r \rceil-1-k\right)(k-1)^{3/2}}.
\end{align*}
This implies that
\begin{align*}
	&\sum_{k=k_n+1}^{n-\lceil n^r \rceil-3}
	\mathbf{E}_y
	\Big(e^{-S_{k-1}/2}  \frac{(1+S_{k-1})^2}{n-\lceil n^r \rceil-1-k}1_{\{\min_{j\leq k-1} S_j\geq 0\}} \Big)\\
	& \lesssim (1+y)\sum_{k=k_n+1}^{n-\lceil n^r \rceil-3}\frac{1}{\left(n-\lceil n^r \rceil-1-k\right)(k-1)^{3/2}}\\
	& \leq (1+y)\int_{k_n}^{n-\lceil n^r \rceil-3}\frac{\mathrm{d}x}{(n-\lceil n^r \rceil-2-x)(x-1)^{3/2}}\\
	& = (1+y) \left(\frac{\log \left(n-\lceil n^r \rceil-2-k_n \right)}{(k_n-1)^{3/2}} -\frac{3}{2} \int_{k_n}^{n-\lceil n^r \rceil-2} \frac{\log \left( n-\lceil n^r \rceil-2-x\right)\mathrm{d}x }{(x-1)^{5/2}}\right)\\
	& \leq (1+y)\frac{\log \left(n-\lceil n^r \rceil-2-k_n \right)}{(k_n-1)^{3/2}} \lesssim (1+y)\frac{(\log n)^2}{n^{3/2}}.
\end{align*}
Combining this with \eqref{e:rs11}, we get
\begin{equation}\label{e:rs15}
II_3 \lesssim \frac{1}{n}+ y\frac{(\log n)^2}{n^{3/2}}.
\end{equation}
Combining \eqref{step_30}, \eqref{e:rs12},\eqref{e:rs13}, \eqref{e:rs14} and \eqref{e:rs15}, we get
\begin{align}
	&\mathbb{E}_y\bigg(\bigg(
	\sum_{x\in \mathcal{N}(n-\lceil n^r \rceil-1)}
	e^{-V(x)}1_{\left\{ \min_{j\leq n-\lceil n^r \rceil-1} V(x_j)\geq  0 \right\}}1_{\left\{x\in \mathcal{A}_{n-\lceil n^r \rceil-1,1}^{0, 0}\right\}} \bigg)^{1+\alpha} \bigg)\nonumber\\
	& \lesssim e^{-y}\left( \left(\frac{1}{n}+ y \frac{(\log n)^2}{n^{3/2}}\right) +  \left(\frac{1}{n}+ y \frac{(\log n)^2}{n^{3/2}}\right)^\alpha \right).
\end{align}
Thus, when $n$ is large enough,
\begin{align}\label{estimate-term2}
&\mathbb{E}\Big(\left|\widetilde{W}_{n,1}^{\lceil n^{r}\rceil ,0}-\mathbb{E}\left(\widetilde{W}_{n,1}^{\lceil n^{r}\rceil ,0} \Big| \mathcal{F}_{\lceil n^r \rceil} \right)  \right|^{1+\alpha}\Big| \mathcal{F}_{\lceil n^r \rceil} \Big)\nonumber\\
	& \lesssim
	\sum_{u\in \mathcal{N}(\lceil n^r \rceil)}
	e^{-V(u)}\left(\left(\frac{1}{n}+ V(u) \frac{(\log n)^2}{n^{3/2}}\right)+\left(\frac{1}{n}+ V(u) \frac{(\log n)^2}{n^{3/2}}\right)^\alpha \right) 1_{\{V(u)\geq 0\}}\nonumber\\
	& \lesssim U_{\lceil n^r \rceil} + W_{\lceil n^r \rceil}^{1-\alpha}U_{\lceil n^r \rceil}^\alpha,
\end{align}
with
\[
U_{\lceil n^r \rceil}:=
	\sum_{u\in \mathcal{N}(\lceil n^r \rceil)}
e^{-V(u)}\Big(\frac{1}{n}+ V(u) \frac{(\log n)^2}{n^{3/2}}\Big)1_{\{V(u)\geq 0\}}.
\]

{\bf Step 3}
Now combining \eqref{term1+term2}, \eqref{estimate-term1} and \eqref{estimate-term2}, we get
\begin{align*}
	&\mathbb{P}\Big(\Big|\widetilde{W}_{n}^{\lceil n^{r}\rceil ,0} - \mathbb{E}\left(\widetilde{W}_{n}^{\lceil n^{r}\rceil ,0} \Big| \mathcal{F}_{\lceil n^r \rceil}\right) \Big|> 3	n^{-\beta-\frac{1}{2}} \big| \mathcal{F}_{\lceil n^r \rceil}\Big)\\
	& \lesssim n^\beta W_{\lceil n^r \rceil} +n^{\beta -\frac{1}{2}}
		\sum_{u\in \mathcal{N}(\lceil n^r \rceil)}
	V(u)e^{-V(u)}1_{\{V(u)\geq 0\}} +n^{(\beta + \frac{1}{2})(1+\alpha)} \left(U_{\lceil n^r \rceil} + W_{\lceil n^r \rceil}^{1-\alpha}U_{\lceil n^r \rceil}^\alpha \right).
\end{align*}
Using \eqref{e:rs17} and \eqref{e:rs18},
we get that for any $\varepsilon \in (0, (1+r)\land (\frac{3}{2}))$, $n^\varepsilon U_{n^{r}}\to 0$ in probability.
Therefore, we have  \eqref{toprove} and the proof is complete.
\hfill$\Box$


\begin{singlespace}

\end{singlespace}

\vskip 0.2truein
\vskip 0.2truein

\vskip 0.2truein
\vskip 0.2truein

\noindent{\bf Haojie Hou:}  School of Mathematical Sciences, Peking
University,   Beijing, 100871, P.R. China. Email: {\texttt
	houhaojie@pku.edu.cn}

\smallskip

\noindent{\bf Yan-Xia Ren:} LMAM School of Mathematical Sciences \& Center for
Statistical Science, Peking
University,  Beijing, 100871, P.R. China. Email: {\texttt
	yxren@math.pku.edu.cn}

\smallskip
\noindent {\bf Renming Song:} Department of Mathematics,
University of Illinois Urbana-Champaign,
Urbana, IL 61801, U.S.A.
Email: {\texttt rsong@illinois.edu}

\end{document}